\documentclass[11pt,a4paper]{article}
\usepackage{aligned-overset}
\usepackage[T1]{fontenc}
\usepackage[american]{babel}
\usepackage{comment}
\usepackage{mathtools,amssymb,amsthm}
\usepackage{authblk}
\usepackage{booktabs}
\usepackage{graphicx}
\usepackage{tikz}
\usetikzlibrary{arrows.meta}
\usepackage{subcaption}
\usepackage{placeins}
\usepackage{pgfplots}
\usepackage{tikz}
\usepgfplotslibrary{groupplots}
\pgfplotsset{compat=1.18}
\usepackage{geometry}
\usepackage{microtype}
\usepackage{xcolor}
\usepackage{hyperref}
\usepackage{csquotes}
\usepackage[
  backend=biber,
  style=numeric-comp,
  maxnames=99,
  giveninits=true,
  isbn=false
]{biblatex}
\usepackage[nameinlink,noabbrev,capitalize]{cleveref}
\usepackage{newtxtext,newtxmath}
\usepackage{xspace}

\hypersetup{
  colorlinks=true,
  linkcolor=blue!55!black,
  citecolor=blue!55!black,
  urlcolor=blue!55!black
}

\newcommand{\email}[1]{\href{mailto:#1}{#1}}

\numberwithin{equation}{section}

\crefname{remark}{Remark}{Remarks}
\Crefname{remark}{Remark}{Remarks}
\crefname{corollary}{Corollary}{Corollaries}
\Crefname{corollary}{Corollary}{Corollaries}
\AddToHook{env/remark/begin}{\crefalias{theorem}{remark}}
\AddToHook{env/lemma/begin}{\crefalias{theorem}{lemma}}
\AddToHook{env/proposition/begin}{\crefalias{theorem}{proposition}}
\AddToHook{env/corollary/begin}{\crefalias{theorem}{corollary}}

\newtheorem{theorem}{Theorem}
\newtheorem{proposition}[theorem]{Proposition}
\newtheorem{lemma}[theorem]{Lemma}
\newtheorem{corollary}[theorem]{Corollary}
\theoremstyle{definition}

\theoremstyle{remark}
\newtheorem{remark}[theorem]{Remark}

\newcommand{\Real}{\mathbb{R}}

\newcommand{\dd}{\operatorname{div}\operatorname{div}}
\newcommand{\Hess}{\nabla^2}

\newcommand{\polyrust}{\texttt{poly\_rust}\xspace}

\newcommand{\Mh}{\mathcal M_h}
\newcommand{\Th}{\mathcal T_h}
\newcommand{\Eh}{\mathcal E_h}
\newcommand{\Vh}{\mathcal V_h}
\newcommand{\ET}{\mathcal E_T}
\newcommand{\VE}{\mathcal V_E}
\newcommand{\VT}{\mathcal V_T}

\newcommand{\Poly}[1]{\mathcal{P}^{#1}}

\newcommand{\divdivKernelComplement}[1]{\mathcal{K}^{{\rm c}, #1}}

\newcommand{\lproj}[2]{\pi^{#1}_{#2}}

\newcommand{\Htwo}[2]{\underline{\mathcal H}_{2,#2}^{#1}}
\newcommand{\IHtwo}[2]{\underline{\mathcal I}_{2,#2}^{#1}}

\newcommand{\norm}[2]{\lVert #2\rVert_{#1}}
\newcommand{\seminorm}[2]{\lvert #2\rvert_{#1}}

\title{A Discrete de Rham method for linear and nonlinear fourth-order plate models}
\author[1]{Daniele A. Di Pietro}
\author[1,2]{J\'er\^ome Droniou}
\author[1]{Arax Leroy}
\author[1]{Vito Patierno}
\affil[1]{IMAG, Univ. Montpellier, CNRS, Montpellier, France, \email{daniele.di-pietro@umontpellier.fr}, \email{jerome.droniou@cnrs.fr}, \email{arax.leroy@umontpellier.fr}}
\affil[2]{School of Mathematics, Monash University, Melbourne, Australia}

\begin{document}

\maketitle

\begin{abstract}
  We develop and analyze a Discrete de Rham (DDR) method for the Kirchhoff--Love plate problem on general polygonal meshes and demonstrate its extension to geometrically nonlinear plates.
  The construction uses local degrees of freedom and polynomial reconstructions, avoiding the implementation difficulties associated with globally $C^1$-conforming approximation spaces.
  A new deflection reconstruction obtained from the discrete Hessian permits a uniform treatment of the load term at all polynomial degrees.
  The use of this reconstruction improves the convergence rates of previous DDR methods and enables the formulation of a lowest-order scheme with only vertex unknowns.
  We provide a stability and convergence analysis of the method and validate the theoretical results through a comprehensive set of numerical tests.
  The method is then extended to the coupled von K\'arm\'an equations, for which we demonstrate convergence numerically.
  We also study the flat-branch spectrum of a compressed perforated plate and its small-amplitude forced nonlinear response near the first critical load.
\end{abstract}

\begingroup
\small
\noindent\textbf{Keywords:} discrete de Rham methods; polygonal meshes; Kirchhoff--Love plates; Hessian reconstruction; von K\'arm\'an equations; buckling.
\par\smallskip
\noindent\textbf{Mathematics Subject Classification (2020):} 65N30, 65N12, 65N15, 74G60, 74K20.
\par
\endgroup
\medskip


\section{Introduction}
\label{sec:introduction}

In this work, we develop, analyze, and numerically investigate an arbitrary-order Discrete de Rham (DDR) method for fourth-order plate problems on general polygonal meshes.
Unlike traditional finite element discretizations, our method does not rely on spaces of globally regular functions, but rather on local reconstructions computed directly from the degrees of freedom (DOFs).
This considerably simplifies its implementation.
The main novelty is a new Hessian-based deflection reconstruction that yields a higher convergence rate than previous DDR methods.
\smallskip

The Kirchhoff--Love theory describes the deformation of thin elastic plates under forces and moments and leads to a fourth-order elliptic problem with a variational formulation in $H^2(\Omega)$; see \cite{Ciarlet:21}.
Its numerical approximation therefore requires either global $C^1$ continuity \cite{Argyris.Fried.ea:68,Clough.Tocher:65,Bell:69} or specific stabilization terms to control the bending energy, as in \cite{Brenner.Sung:05}.
The fourth-order structure of the Kirchhoff--Love problem and the related $C^1$ continuity requirement for conforming approximations also appear in geometrically nonlinear plate models, such as the von K\'arm\'an equations, whose finite element approximation has been considered, e.g., in \cite{Brezzi:78,Mallik.Nataraj:16,Carstensen.Nataraj:21,Miyoshi:76,Brenner.Neilan.ea:17,Carstensen.Mallik.ea:19}.
Constructing and transforming bases for $C^1$-conforming finite elements is, however, not straightforward, particularly at higher polynomial degrees \cite{Kirby:18,Chen.Chen.ea:25}.
\smallskip

In addition to supporting general meshes, polytopal methods avoid the need to construct discrete function spaces with global regularity by using reconstructions and projections computed directly from the local DOFs.
Virtual Element Methods retain global function spaces for the purpose of presentation and analysis; see \cite{Brezzi.Marini:13,Antonietti.Manzini.ea:18} concerning Kirchhoff--Love plates and \cite{Lovadina.Mora.ea:21,Shylaja.Kumar:24} concerning the von K\'arm\'an equations.
Fully discrete approaches such as the Weak Galerkin method \cite{Mu.Wang.ea:14} and the Hybrid High-Order (HHO) method \cite{Bonaldi.Di-Pietro.ea:18,Dong.Ern:22}, on the other hand, have formulations that more closely reflect their implementation, with discrete spaces spanned by vectors of local polynomials.
The more recent Discrete de Rham (DDR) methods follow a similar approach, but also organize discrete spaces and reconstructed differential operators to form compatible discrete complexes \cite{Di-Pietro.Droniou:23,Di-Pietro.Droniou.ea:25}.
\smallskip

The present work develops a DDR method for linear and nonlinear plate problems on polygonal meshes.
Its construction is closely related to the discrete Stokes--de Rham diagram and its Bernstein--Gelfand--Gelfand reduction proposed in \cite{Di-Pietro.Droniou.ea:25}.
In \cite{Di-Pietro.Droniou.ea:26}, this diagram is used to devise a Kirchhoff--Love discretization based on a degree-$k$ Hessian complex, for $k\ge1$ in our notation (which differs from that of \cite{Di-Pietro.Droniou.ea:25} by a shift of $+1$ in the polynomial degree).
The method converges at a rate of $h^k$ in the energy norm:
the limiting contribution arises from the adjoint consistency of a degree-$k$ potential reconstructed through the discrete gradient.

In this work, we show that relatively minor modifications of the scheme can increase this convergence rate by one order.
We achieve this by enriching the element DOFs (as already done in \cite{Leroy:26}, which applies the construction of \cite{Di-Pietro.Droniou.ea:25} to the Reissner--Mindlin plate problem) and, most importantly, by introducing a new potential reconstruction based on the reconstructed Hessian.
The use of this potential in the load term makes it possible to compare the discrete and continuous problems at the Hessian level, without composing the adjoint-consistency estimates of two first-order operators.
This results in an energy error estimate of order $h^{k+1}$, also covering the case $k = 0$, which was not considered in \cite{Di-Pietro.Droniou.ea:26}.
\smallskip

After matching polynomial degrees, our discrete space coincides with one of the choices contemplated in \cite{Brezzi.Marini:13}; see \cref{rem:comparison-bm}.
That method achieves convergence rates in the energy norm that match those obtained here, but requires special definitions of the load term at the lowest degrees.
In addition to introducing a uniform treatment of the load for all polynomial degrees through the new deflection reconstruction, the present work fills an important theoretical gap by avoiding the uniform norm-equivalence assumption \cite[Assumption~H2]{Brezzi.Marini:13} in the analysis.
As pointed out in \cite{Beirao-da-Veiga.Di-Pietro.ea:26}, this norm equivalence can be difficult to establish even for simpler problems; the corresponding assumption can be avoided using fully discrete analysis techniques in the spirit of the Third Strang Lemma \cite{Di-Pietro.Droniou:18}.
These techniques also naturally incorporate load discretizations that do not involve $L^2$-orthogonal projections of virtual functions; see \cite[Remark~3.18]{Beirao-da-Veiga.Di-Pietro.ea:26}.
\smallskip

To demonstrate the performance of the proposed method in more complex situations, we apply it to the coupled von K\'arm\'an equations.
We present a manufactured-solution convergence test, a study of the flat-branch spectrum for a compressed perforated plate, and a study of a small-amplitude forced nonlinear response near the first critical load.
These experiments illustrate the use of the method for both linearized buckling computations and the solution of the coupled nonlinear system.
\smallskip

The paper is organized as follows.
Section~\ref{sec:setting} introduces the continuous problem, the mesh, the discrete spaces, and the local reconstructions.
Section~\ref{sec:discrete-problem} defines the discrete problem, states the main error estimate, and presents numerical validation.
Section~\ref{sec:analysis} establishes the properties of the reconstructions and proves the stability and consistency estimates underlying the error analysis.
Finally, Section~\ref{sec:von-karman} applies the method to the von K\'arm\'an equations and presents numerical experiments on convergence, buckling, and nonlinear response.


\section{Setting}
\label{sec:setting}

\subsection{Continuous problem}

Let $\Omega\subset\Real^2$ be a bounded polygonal domain and let $f\in L^2(\Omega)$.
For simplicity, the flexural rigidity is normalized to one and the isotropic biharmonic model is considered:
\begin{equation}
  \Delta^2 u=f \quad\text{in }\Omega,
  \qquad
  u=\partial_n u=0 \quad\text{on }\partial\Omega.
  \label{eq:strong-problem}
\end{equation}
Writing $H_0^2(\Omega)$ for the subspace of $H^2(\Omega)$ with vanishing trace and normal derivative, the weak problem is: find $u\in H_0^2(\Omega)$ such that
\begin{equation}
  a(u,v)\coloneqq \int_\Omega \Hess u:\Hess v
  =\int_\Omega fv
  \qquad\forall v\in H_0^2(\Omega).
  \label{eq:continuous-weak}
\end{equation}

\subsection{Mesh}
\label{sec:setting:mesh}

Let $\mathcal M_h\coloneqq(\Th,\Eh,\Vh)$ be a conforming polygonal mesh of $\Omega$ in the usual sense; see, for instance, \cite[Definition~1.4]{Di-Pietro.Droniou:20}.
The diameter of an element $T$ is $h_T$, and $h\coloneqq\max_{T\in\Th}h_T$.
For the analysis, we assume that there is a mesh-independent parameter $\theta>0$ such that every element $T$ is star-shaped with respect to a ball of radius $\theta h_T$ centered at $x_T$.
This assumption provides discrete trace and inverse estimates as well as uniform polynomial approximation; see \cite{Droniou.Yemm:22} and \cite[Chapter 2]{Di-Pietro.Droniou.ea:27}.
Given a mesh element $T\in\Th$, the sets collecting the edges and vertices of $T$ are denoted by $\ET$ and $\VT$.

Each edge $E$ is equipped once and for all with a unit tangent $ t_E = (t_{E,1}, t_{E,2})$ and with the unit normal $ n_E\coloneqq(-t_{E,2},t_{E,1})$.
We denote by $\omega_{TE}\in\{-1,1\}$ the sign such that $ n_{TE}\coloneqq\omega_{TE} n_E$, with $ n_{TE}$ denoting the unit normal pointing out of $T$.
We also introduce the endpoint orientations $\omega_{EV}\in\{-1,1\}$ such that $\omega_{EV} t_E$ points toward $V\in\VE$.

With these orientation conventions, the following key integration-by-parts identity holds:
\begin{equation}
  \int_T \Hess v:\tau
  =\int_T v\,\dd\tau
  - \sum_{E \in \ET} \omega_{TE} \int_E v\,\operatorname{div}\tau\cdot n_E
  + \sum_{E \in \ET} \omega_{TE} \int_E \nabla v \cdot (\tau n_E),
  \label{eq:ibp-hessian}
\end{equation}
where the matrix divergence is taken row by row: $(\operatorname{div}\tau)_i=\sum_{j=1}^2\partial_j\tau_{ij}$.

For nonnegative quantities $A$ and $B$, we write $A\lesssim B$ if $A\le CB$ for a constant $C>0$ independent of the mesh size, the particular mesh entities, and the functions appearing in the estimate, but possibly depending on $\Omega$, $\theta$, and the polynomial degree of the method.

\subsection{Local polynomial spaces}

For a mesh entity $Y$, $\Poly\ell(Y)$ is the space of restrictions to $Y$ of polynomials in two variables of total degree at most $\ell$, with the convention $\Poly\ell(Y)\coloneqq\{0\}$ for $\ell<0$.
The $L^2$ projector onto this space is $\pi_Y^\ell$.
The same notation is used for componentwise vector- and matrix-valued projectors.
The zero-average polynomial subspace is denoted by $\Poly{\ell}_0(Y)$. 

\subsection{Discrete deflection space and interpolator}

Let $k \ge 0$ and define the global $H^2$-like deflection space
\begin{equation}
  \label{eq:global-space}
  \Htwo{k}{h} \coloneqq \Big\{
  \begin{aligned}[t]
    &\underline{v}_h = ( (v_T)_{T \in \Th}, (v_E, G_{v,E,n})_{E \in \Eh}, (v_V, G_{v,V})_{V \in \Vh}) \;:\;
    \\
    &\quad
    \text{$v_T \in \Poly{k-2}(T)$ for all $T \in \Th$},
    \\
    &\quad
    \text{$v_E \in \Poly{k-2}(E)$ and $G_{v,E,n} \in \Poly{k-1}(E)$ for all $E \in \Eh$},
    \\
    &\quad
    \text{$v_V \in \Real$ and $G_{v,V} \in \Real^2$ for all $V \in \Vh$},
    \Big\},
  \end{aligned}
\end{equation}
This space corresponds to the discrete $H^2$-like space of \cite{Di-Pietro.Droniou.ea:25,Di-Pietro.Droniou.ea:26}, with the index shift $k\to k-1$ and enriched element DOFs (since no serendipity reduction is performed).
The restrictions of the above space and of the corresponding global vectors of polynomials to a mesh element or edge $Y \in \Mh$ are denoted by replacing the subscript $h$ by $Y$.
The local DOFs in $\Htwo{k}{T}$ for a hexagonal element $T$ and $k \in \{0, 1, 2\}$ are illustrated in \cref{fig:deflection-dofs}, while \cref{fig:dof-comparison} contains a comparison with the lowest-order versions of common triangular finite elements for plate problems.

\begin{figure}[!htbp]
  \centering
\begin{tikzpicture}[
  line width=0.6pt,
  dof moment/.style={draw=blue!65!black, fill=blue!65!black},
  dof normal/.style={draw=orange!80!black, -{Stealth[length=3.5pt,width=3pt]}}
]
  \foreach \k/\ndofs in {0/18,1/24,2/37} {
    \begin{scope}[xshift={4.5*\k cm}]
      \foreach \i in {0,...,5} {
        \coordinate (V\i) at ({60*\i}:1.3);
      }
      \draw (V0) -- (V1) -- (V2) -- (V3) -- (V4) -- (V5) -- cycle;

      \foreach \i in {0,...,5} {
        \pgfmathsetmacro{\angle}{60*\i}
        \pgfmathsetmacro{\nextangle}{60*(\i+1)}
        \pgfmathsetmacro{\normalangle}{60*\i+30}
        \ifnum\k>0
          \foreach \j in {1,...,\k} {
            \pgfmathsetmacro{\fraction}{\j/(\k+1)}
            \pgfmathsetmacro{\px}{1.3*((1-\fraction)*cos(\angle)+\fraction*cos(\nextangle))}
            \pgfmathsetmacro{\py}{1.3*((1-\fraction)*sin(\angle)+\fraction*sin(\nextangle))}
            \draw[dof normal] (\px,\py) --
              ({\px+0.24*cos(\normalangle)},{\py+0.24*sin(\normalangle)});
          }
        \fi
        \ifnum\k>1
          \pgfmathtruncatemacro{\nscalar}{\k-1}
          \foreach \j in {1,...,\nscalar} {
            \pgfmathsetmacro{\fraction}{\j/\k}
            \pgfmathsetmacro{\px}{1.3*((1-\fraction)*cos(\angle)+\fraction*cos(\nextangle))}
            \pgfmathsetmacro{\py}{1.3*((1-\fraction)*sin(\angle)+\fraction*sin(\nextangle))}
            \path[dof moment] (\px,\py) ++(-1.6pt,-1.6pt) rectangle ++(3.2pt,3.2pt);
          }
        \fi
      }

      \foreach \i in {0,...,5} {
        \draw[fill=white] (V\i) circle[radius=3.4pt];
        \fill (V\i) circle[radius=1.3pt];
      }

      \ifnum\k=2
        \path[dof moment] (-1.6pt,-1.6pt) rectangle ++(3.2pt,3.2pt);
      \fi
      \node[font=\small] at (0,-1.65) {$k=\k$\,: \ndofs\ DOFs};
    \end{scope}
  }
\end{tikzpicture}
\par\smallskip
{\small
\begin{tabular}{@{}cl@{\qquad}cl@{}}
  \tikz[baseline=-0.5ex]{\fill (0,0) circle[radius=1.3pt];}
  & Vertex value $v_V$
  & \tikz[baseline=-0.5ex]{\draw[line width=0.6pt] (0,0) circle[radius=3.4pt];}
  & Vertex gradient $G_V$ (two DOFs)
  \\
  \tikz[baseline=-0.5ex]{\fill[blue!65!black] (-1.6pt,-1.6pt) rectangle ++(3.2pt,3.2pt);}
  & Scalar moments $v_E$, $v_T$
  & \tikz[baseline=-0.5ex]{\draw[orange!80!black,line width=0.6pt,-{Stealth[length=3.5pt,width=3pt]}] (0,-0.10) -- (0,0.14);}
  & Normal-gradient moments $G_{v,E,n}$
\end{tabular}}
  \caption{DOFs of $\Htwo{k}{T}$ on a hexagon.
    Each vertex carries one value and two gradient components.
    Each edge carries $\max\{k-1,0\}$ scalar moments and $k$ normal-gradient moments;
    the element carries $\frac{k(k-1)}{2}$ scalar moments.
    Nonvertex symbols indicate moments, not point evaluations.}
  \label{fig:deflection-dofs}
\end{figure}
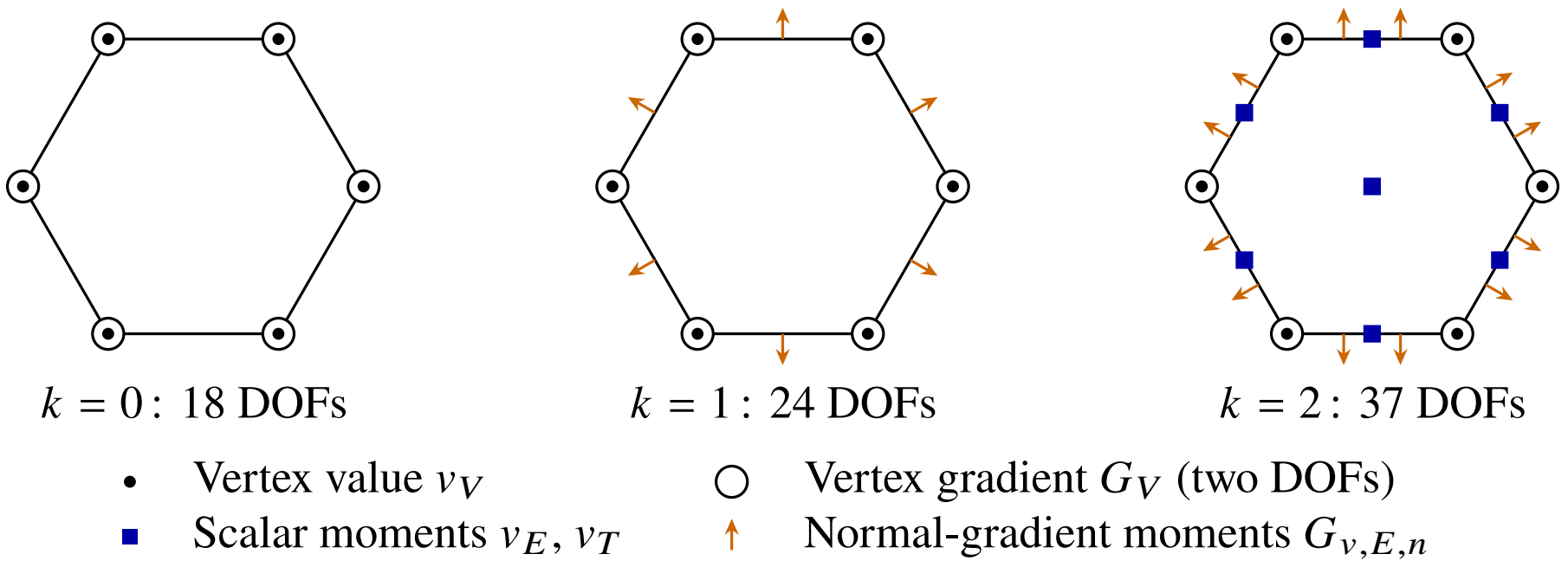

\begin{figure}[!htbp]
  \centering
\begingroup
\tikzset{
  dofcmp moment/.style={draw=blue!65!black,fill=blue!65!black},
  dofcmp normal/.style={
    draw=orange!80!black,-{Stealth[length=3.5pt,width=3pt]}}
}
\newcommand{\dofcmptriangle}{%
  \coordinate (A) at (-1.1,-0.65);
  \coordinate (B) at ( 1.1,-0.65);
  \coordinate (C) at ( 0,1.25526);
  \coordinate (O) at (0,-0.01491);
  \draw (A)--(B)--(C)--cycle;
}
\newcommand{\dofcmpvertices}[1]{%
  \foreach \v in {A,B,C}{%
    \ifnum#1=1
      \draw[fill=white] (\v) circle[radius=5.8pt];
    \fi
    \draw[fill=white] (\v) circle[radius=3.4pt];
    \fill (\v) circle[radius=1.3pt];
  }%
}
\newcommand{\dofcmpsquare}[1]{%
  \path[dofcmp moment] #1
    ++(-1.6pt,-1.6pt) rectangle ++(3.2pt,3.2pt);
}
\newcommand{\dofcmpddr}[1]{%
  \dofcmptriangle
  \foreach \a/\b/\normalangle in {A/B/-90,B/C/30,C/A/150}{%
    \ifnum#1>0
      \foreach \j in {1,...,#1}{%
        \pgfmathsetmacro{\fraction}{\j/(#1+1)}
        \path (\a)--(\b) coordinate[pos=\fraction] (M);
        \draw[dofcmp normal] (M)--++(\normalangle:0.24);
      }%
    \fi
    \ifnum#1>1
      \pgfmathtruncatemacro{\nscalar}{#1-1}
      \foreach \j in {1,...,\nscalar}{%
        \pgfmathsetmacro{\fraction}{\j/#1}
        \path (\a)--(\b) coordinate[pos=\fraction] (M);
        \dofcmpsquare{(M)}
      }%
    \fi
  }%
  \ifnum#1=2
    \dofcmpsquare{(O)}
  \fi
  \ifnum#1=3
    \dofcmpsquare{(-0.22,-0.16)}
    \dofcmpsquare{( 0.22,-0.16)}
    \dofcmpsquare{(0,0.22)}
  \fi
  \dofcmpvertices{0}
}
\newcommand{\dofcmpedgemoments}{%
  \foreach \a/\b/\normalangle in {A/B/-90,B/C/30,C/A/150}{%
    \path (\a)--(\b) coordinate[midway] (M);
    \draw[dofcmp normal] (M)--++(\normalangle:0.24);
  }%
}
\newcommand{\dofcmpddrpic}[1]{%
  \begin{tikzpicture}[line width=0.6pt,scale=0.9,baseline=-0.10cm]
    \dofcmpddr{#1}
  \end{tikzpicture}%
}
\newcommand{\dofcmpfempic}[1]{%
  \begin{tikzpicture}[line width=0.6pt,scale=0.9,baseline=-0.10cm]
    \dofcmptriangle
    \ifnum#1=1 \dofcmpedgemoments \fi
    \ifnum#1=3 \dofcmpedgemoments \fi
    \ifnum#1<2
      \dofcmpvertices{0}
    \else
      \dofcmpvertices{1}
    \fi
  \end{tikzpicture}%
}
\newcommand{\dofcmpvertexsymbol}[1]{%
  \tikz[baseline=-0.5ex,line width=0.6pt]{%
    \ifnum#1=1
      \draw[fill=white] (0,0) circle[radius=5.8pt];
    \fi
    \draw[fill=white] (0,0) circle[radius=3.4pt];
    \fill (0,0) circle[radius=1.3pt];
  }%
}
{\small
\setlength{\tabcolsep}{6pt}
\begin{tabular}{@{}clcrcr@{}}
  \toprule
  & \multicolumn{2}{c}{DDR} & \multicolumn{3}{c}{Classical} \\
  \cmidrule(lr){2-3}\cmidrule(l){4-6}
  $k$ & Local DOFs & Count & Name & Local DOFs & Count \\
  \midrule
  0 &  \dofcmpddrpic{0} & 9  & Reduced HCT \cite{Ciarlet:78} & \dofcmpfempic{0} & 9  \\
  \addlinespace[4pt]
  1 &  \dofcmpddrpic{1} & 12 & HCT \cite{Ciarlet:74}        &\dofcmpfempic{1} & 12 \\
  \addlinespace[4pt]
  2 & \dofcmpddrpic{2} & 19 & Bell  \cite{Bell:69}      & \dofcmpfempic{2} & 18 \\
  \addlinespace[4pt]
  3 &  \dofcmpddrpic{3} & 27 & Argyris  \cite{Kirby:18}   &\dofcmpfempic{3} & 21 \\
  \bottomrule
\end{tabular}}
\par\medskip
{\footnotesize
\renewcommand{\arraystretch}{1.45}
\begin{tabular}{@{}cl@{\qquad}cl@{}}
  \dofcmpvertexsymbol{0}
  & Vertex value and gradient (3 DOFs)
  & \dofcmpvertexsymbol{1}
  & Value, gradient and Hessian (6 DOFs)
  \\
  \tikz[baseline=-0.5ex]{%
    \path[dofcmp moment]
      (-1.6pt,-1.6pt) rectangle ++(3.2pt,3.2pt);}
  & Scalar moment (1 DOF)
  & \tikz[baseline=-0.5ex,line width=0.6pt]{%
      \draw[dofcmp normal] (0,-0.10)--(0,0.14);}
  & Normal-derivative moment (1 DOF)
\end{tabular}}
\endgroup

  \caption{Local DOFs and counts on a triangle for spaces that contain $\Poly{k+2}(T)$.
    HCT denotes Hsieh--Clough--Tocher.
    HCT and Argyris elements are represented using the equivalent moment-based choice of edge data.}
  \label{fig:dof-comparison}
\end{figure}

The meaning of the components of the discrete space is specified by the corresponding interpolator.
The interpolant of a sufficiently smooth function $v : \Omega \to \Real$ in $\Htwo{k}{h}$ is given by
\[
\IHtwo{k}{h} v \coloneqq \big(
(\lproj{k-2}{T} v)_{T \in \Th},
(\lproj{k-2}{E} v, \lproj{k-1}{E} (\nabla v \cdot n_E) )_{E \in \Eh},
(v(x_V), \nabla v(x_V))_{V \in \Vh}
\big),
\]
where $x_V \in \Real^2$ denotes the coordinate vector of the vertex $V$.

\subsection{Local reconstructions}
\label{sec:setting:local-reconstructions}

All reconstruction operators below are computed by solving small linear systems local to an edge or an element $Y\in\Eh\cup\Th$.
Their definitions mimic continuous integration-by-parts formulas, with polynomial projections and point values replaced by the corresponding DOFs or previously reconstructed quantities.
In what follows, $\underline{v}_Y$ denotes a generic element of $\Htwo{k}{Y}$.

\subsubsection{Edge operators}

A prime denotes differentiation in the direction $t_E$.
The edge Hessian
$H_E^k\underline v_E\in\Poly{k}(E;\Real^2)$
approximates $(\Hess v)_{|E}t_E$ and is defined by
\begin{equation}\label{eq:def.HE}
  \begin{aligned}
    \int_E H_E^k\underline v_E\cdot w
    &=
    \int_E v_E\,(w''\cdot t_E)
    -\int_E G_{v,E,n}\,(w'\cdot n_E)
    \\
    &\quad
    +\sum_{V\in\VE}\omega_{EV}
    \left[
      G_{v,V}\cdot w(x_V)
      -v_V\,(w'(x_V)\cdot t_E)
    \right]
    \quad\forall w\in\Poly{k}(E;\Real^2).
  \end{aligned}
\end{equation}

The edge gradient
$G_E^{k+1}\underline v_E\in\Poly{k+1}(E;\Real^2)$
is reconstructed from the edge Hessian by discrete integration
by parts:
\begin{equation}\label{eq:def.GE}
  \int_E G_E^{k+1}\underline v_E\cdot w'
  =
  -\int_E H_E^k\underline v_E\cdot w
  +\sum_{V\in\VE}\omega_{EV}\,G_{v,V}\cdot w(x_V)
  \quad\forall w\in\Poly{k+2}_0(E;\Real^2).
\end{equation}
The deflection trace
$\gamma_E^{k+2}\underline v_E\in\Poly{k+2}(E)$
is then defined by
\begin{equation}\label{eq:def.traceE}
  \int_E \gamma_E^{k+2}\underline v_E\,r'
  =
  -\int_E (G_E^{k+1}\underline v_E\cdot t_E)\,r
  +\sum_{V\in\VE}\omega_{EV}\,v_V\,r(x_V)
  \quad\forall r\in\Poly{k+3}_0(E).
\end{equation}

\begin{remark}[Validity of the reconstruction formulas]\label{rem:validity.formulae}
  Taking $w\in\Poly{0}(E;\Real^2)$ in \eqref{eq:def.HE} shows that the right-hand side of \eqref{eq:def.GE} vanishes for constant test functions $w$.
  Hence, \eqref{eq:def.GE} actually holds for all $w\in\Poly{k+2}(E;\Real^2)$.

  If $k\ge 1$, taking $r\in\Poly{0}(E)$ and choosing $w=r (x \cdot t_E) t_E\in\Poly{k}(E;\Real^2)$ (which satisfies $w''=  0$ and $w'\cdot n_E=0$) in \eqref{eq:def.HE} shows that
  \[
  -\int_E H_E^k\underline v_E\cdot w
  + \sum_{V\in\VE}\omega_{EV} G_{v,V} \cdot w(x_V)
  = \sum_{V\in\VE}\omega_{EV}v_V\,r.
  \]
  Substituting this into \eqref{eq:def.GE} gives
  $
  \int_E (G_E^{k+1}\underline v_E\cdot t_E)r=\sum_{V\in\VE}\omega_{EV}v_V\,r
  $,
  which, upon substitution into \eqref{eq:def.traceE}, shows that this formula also holds for constant $r$ and hence for every $r\in\Poly{k+3}(E)$.
  If $k=0$, this argument no longer applies, and \eqref{eq:def.traceE} is not guaranteed to hold for constant $r$.
\end{remark}

\subsubsection{Element operators}

The edge reconstructions provide the boundary data in the discrete counterparts of the integration-by-parts identity \eqref{eq:ibp-hessian}.
Specifically, the element Hessian $H_T^k\underline v_T\in\Poly{k}(T;\Real^{2\times2})$ is defined by
\begin{equation}
  \begin{aligned}
    \int_T H_T^k\underline v_T:\tau
    &=
    \int_T v_T\,\dd\tau
    -\sum_{E\in\ET}\omega_{TE}
      \int_E \gamma_E^{k+2}\underline v_E\,
      (\operatorname{div}\tau\cdot n_E)
    \\
    &\quad
    +\sum_{E\in\ET}\omega_{TE}
      \int_E G_E^{k+1}\underline v_E\cdot(\tau n_E)
    \quad\forall\tau\in\Poly{k}(T;\Real^{2\times2}).
  \end{aligned}
  \label{eq:element-hessian}
\end{equation}

\begin{remark}[Symmetry of the Hessian reconstruction]
  \label{rem:HT:symmetry}
  For $k\ge1$, selecting $\tau = J q$ above with $J\coloneqq \left(\begin{smallmatrix}0&-1\\1&0 \end{smallmatrix}\right)$ and $q \in \Poly{k}(T)$ and using the tangential derivative identity in \eqref{eq:edge-operators:characterization}, one can prove that $H_T^{k}\underline{v}_T$ is symmetric.
  For $k=0$, this is no longer the case: for any (possibly non-symmetric) matrix $J \in \Real^{2\times 2}$, vertex gradients $G_{v,V}=Jx_V$ give $G_E^1\underline v_E=Jx$ on every edge and hence $H_T^0\underline v_T=J$.
  
  Note, however, that symmetry of $H_T^{k}\underline{v}_T$ for $k=0$ could be imposed by defining this reconstruction with values in $\Poly{k}(T;\mathbb{S})$, with $\mathbb{S}$ denoting the space of symmetric matrices.
  This is not directly useful for our purposes.
\end{remark}

We next reconstruct a new element deflection by using \eqref{eq:ibp-hessian} in reverse, following an approach similar to that in \cite[Section 3.2]{Bonaldi.Di-Pietro.ea:18}.
For the point $x_T$ fixed in \cref{sec:setting:mesh}, introduce
\begin{equation*}
  \divdivKernelComplement{m}(T)
  \coloneqq
  (x-x_T)\otimes(x-x_T)\Poly{m-2}(T),
  \qquad m\ge2.
\end{equation*}
The Hessian-based deflection reconstruction
$p_T^{k+2}\underline v_T\in\Poly{k+2}(T)$
is defined by
\begin{equation}
  \begin{aligned}
    \int_T p_T^{k+2}\underline v_T\,\dd\tau
    &=
    \int_T H_T^k\underline v_T:\tau
    -\sum_{E\in\ET}\omega_{TE}
      \int_E G_E^{k+1}\underline v_E\cdot(\tau n_E)
    \\
    &\quad
    +\sum_{E\in\ET}\omega_{TE}
      \int_E \gamma_E^{k+2}\underline v_E\,
      (\operatorname{div}\tau\cdot n_E)
    \quad\forall\tau\in\divdivKernelComplement{k+4}(T).
  \end{aligned}
  \label{eq:hessian-potential}
\end{equation}
Well-posedness follows from the fact that the map $\dd:\divdivKernelComplement{m}(T)\to\Poly{m-2}(T)$ is an isomorphism; see \cite[Lemma~3.3]{Chen.Huang:22}.

The corresponding global piecewise polynomial reconstruction $p_h^{k+2}:\Htwo{k}{h}\to L^2(\Omega)$ is defined elementwise by
\begin{equation}
  (p_h^{k+2}\underline v_h)_{|T}
  \coloneqq p_T^{k+2}\underline v_T
  \quad\forall T\in\Th.
  \label{eq:global-hessian-potential}
\end{equation}

\begin{figure}[!htbp]
  \begin{tikzpicture}[
    x=1cm,
    y=1cm,
    >=stealth,
    every node/.style={
      font=\small,
      inner sep=3pt,
      outer sep=1pt
    }
  ]
    \node[anchor=east,align=right] at (-1.1,0)
      {Reconstruction};
    \node[anchor=east,align=right] at (-1.1,-1.1)
      {Continuous counterpart};


    \node (HE) at (0,0)
      {$H_E^k\underline v_E$};
    \node (GE) at (2.15,0)
      {$G_E^{k+1}\underline v_E$};
    \node (trace) at (4.3,0)
      {$\gamma_E^{k+2}\underline v_E$};
    \node (HT) at (6.45,0)
      {$H_T^k\underline v_T$};
    \node (potential) at (8.6,0)
      {$p_T^{k+2}\underline v_T$};

    \node at (0,-1.1)    {$(\Hess v)_{|E}\,t_E$};
    \node at (2.15,-1.1) {$(\nabla v)_{|E}$};
    \node at (4.3,-1.1)  {$v_{|E}$};
    \node at (6.45,-1.1) {$(\Hess v)_{|T}$};
    \node at (8.6,-1.1)  {$v_{|T}$};

    \draw[->] (HE.east) -- (GE.west);
    \draw[->] (GE.east) -- (trace.west);
    \draw[->] (trace.east) -- (HT.west);
    \draw[->] (HT.east) -- (potential.west);

    \draw[->] (GE.north)
      .. controls (2.15,0.7) and (6.45,0.7) ..
      (HT.north);

    \draw[->] (GE.north)
      .. controls (2.15,1.3) and (8.6,1.3) ..
      (potential.north);

    \draw[->] (trace.south)
      .. controls (4.3,-0.6) and (8.6,-0.6) ..
      (potential.south);
  \end{tikzpicture}
  \caption{Summary of the local reconstructions introduced in Section \ref{sec:setting:local-reconstructions} together with the respective continuous counterparts and dependencies.}
  \label{fig:local-reconstructions}
\end{figure}
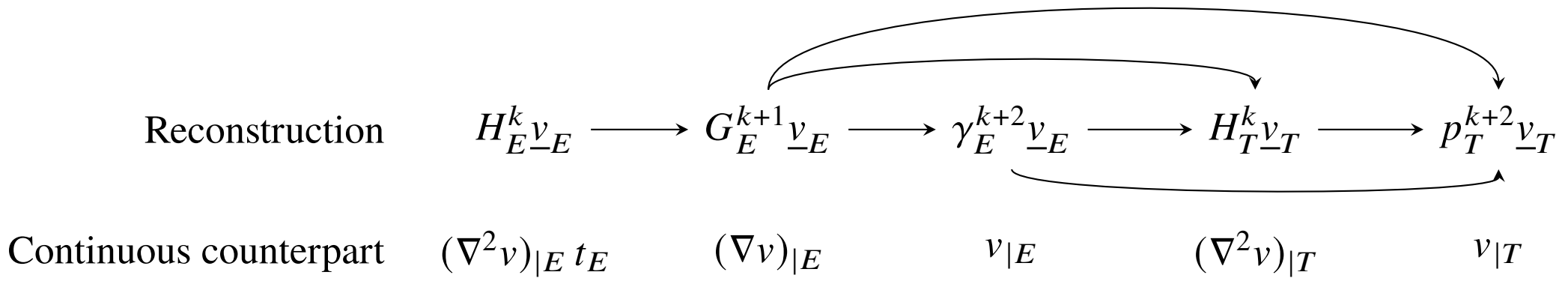

The local reconstructions are summarized in Figure \ref{fig:local-reconstructions} along with their continuous counterparts and the respective dependencies.

\section{Discrete problem, main results, and numerical validation}
\label{sec:discrete-problem}

\subsection{Discrete problem}

For all $T \in \Th$, the local bilinear form $a_T : \Htwo{k}{T} \times \Htwo{k}{T} \to \Real$ is defined by
\begin{equation*}
  a_T(\underline{u}_T,\underline{v}_T)
  \coloneqq\int_T H_T^{k}\underline{u}_T : H_T^{k}\underline{v}_T
  + s_T(\underline{u}_T, \underline{v}_T),
\end{equation*}
where $s_T : \Htwo{k}{T} \times \Htwo{k}{T} \to \Real$ is the following stabilization bilinear form:
\begin{equation}
  \begin{aligned}
    s_T(\underline{u}_T, \underline{v}_T)
    &\coloneqq
    h_T^{-1} \sum_{E \in \ET} \int_E (G_E^{k+1} \underline{u}_E - \nabla p_T^{k+2} \underline{u}_T) \cdot (G_E^{k+1} \underline{v}_E - \nabla p_T^{k+2} \underline{v}_T)
    \\
    &\quad
    + h_T^{-3} \sum_{E \in \ET} \int_E  (\gamma_E^{k+2} \underline{u}_E - p_T^{k+2} \underline{u}_T) (\gamma_E^{k+2} \underline{v}_E - p_T^{k+2} \underline{v}_T).
  \end{aligned}
  \label{eq:potential-stabilization}
\end{equation}

The global bilinear form $a_h : \Htwo{k}{h} \times \Htwo{k}{h} \to \Real$ and linear form $\ell_h : \Htwo{k}{h} \to \Real$ are given by
\begin{equation*}
  a_h(\underline{u}_h,\underline{v}_h)\coloneqq\sum_{T\in\Th}a_T(\underline{u}_T,\underline{v}_T),
  \qquad
  \ell_h(\underline{v}_h)\coloneqq\int_\Omega f\,p_h^{k+2}\underline{v}_h.
\end{equation*}

Denote by $\Htwo{k}{h,0}$ the subspace of $\Htwo{k}{h}$ where all boundary components are set to zero.
The discrete problem is:
find $\underline u_h\in\Htwo{k}{h,0}$ such that
\begin{equation}
  a_h(\underline u_h,\underline v_h)=\ell_h(\underline v_h)
  \qquad\forall\underline v_h\in\Htwo{k}{h,0}.
  \label{eq:discrete-problem}
\end{equation}

\subsection{Energy error estimate}

On $\Htwo{k}{h}$, define the seminorm induced by $a_h$:
\begin{equation}
  \norm{a,h}{\underline{v}_h}\coloneqq a_h(\underline{v}_h,\underline{v}_h)^{\frac12}\qquad
  \forall \underline{v}_h\in \Htwo{k}{h}.
  \label{eq:discrete-energy-seminorm}
\end{equation}
In what follows, $\seminorm{H^\ell(\Th)}{\cdot}$ denotes the standard broken $H^\ell$ seminorm.

\begin{theorem}[Energy error estimate]
  \label{thm:energy-error-estimate}
  Let $u\in H_0^2(\Omega)\cap H^{k+3}(\Th)$ solve the continuous weak problem \eqref{eq:continuous-weak}, and assume additionally $u\in H^4(\Th)$ if $k=0$.
  Then, the discrete problem \eqref{eq:discrete-problem} has a unique solution satisfying the following error estimate:
  \[
  \norm{a,h}{\underline u_h-\IHtwo{k}{h} u}
  \lesssim
  \begin{cases}
    \displaystyle
    h \seminorm{H^{k+3}(\Th)}{u} + h^2 \seminorm{H^{k+4}(\Th)}{u} & \text{if $k = 0$,}
    \\
    h^{k+1}\seminorm{H^{k+3}(\Th)}{u} & \text{if $k \ge 1$,}
  \end{cases}
  \]
\end{theorem}

\begin{proof}
  The bilinear form $a_h$ is coercive on
  $\Htwo{k}{h,0}$ by \cref{cor:norm-property} below,
  so the discrete problem has a unique solution.
  For $w\in H^2(\Omega)\cap H^4(\Th)$ such that
  $\Delta^2w\in L^2(\Omega)$ in the distributional sense,
  introduce the consistency functional
  \begin{equation}
    \label{eq:consistency-functional}
    \mathfrak E_h(w;\underline v_h)
    \coloneqq
    \int_\Omega \Delta^2w\,p_h^{k+2}\underline v_h
    -a_h(\IHtwo{k}{h}w,\underline v_h)
    \qquad
    \forall\underline v_h\in\Htwo{k}{h,0}.
  \end{equation}
  Its dual energy norm is denoted by
  \begin{equation}
    \label{eq:consistency-functional-norm}
    \norm{a,h,*}{\mathfrak E_h(w;\cdot)}
    \coloneqq
    \sup_{\underline v_h\in
      \Htwo{k}{h,0}\setminus\{\underline 0\}}
    \frac{|\mathfrak E_h(w;\underline v_h)|}
         {\norm{a,h}{\underline v_h}}.
  \end{equation}
  Set $\underline e_h\coloneqq
  \underline u_h-\IHtwo{k}{h}u\in\Htwo{k}{h,0}$.
  Since $\Delta^2u=f$, the discrete equation
  \eqref{eq:discrete-problem} gives
  \[
  \norm{a,h}{\underline e_h}^2
  =
  \ell_h(\underline e_h)
  -a_h(\IHtwo{k}{h}u,\underline e_h)
  =
  \mathfrak E_h(u;\underline e_h)
    \le
    \norm{a,h,*}{\mathfrak E_h(u;\cdot)}
    \norm{a,h}{\underline e_h}.
  \]
  Dividing by $\norm{a,h}{\underline e_h}$ when it is nonzero and applying \eqref{eq:adjoint-consistency:convergence} below proves the result.
  This is the coercive case of the Third Strang Lemma of \cite{Di-Pietro.Droniou:18}.
\end{proof}

\begin{remark}[Comparison with the method of \cite{Di-Pietro.Droniou.ea:26}]
  For $k\ge1$, aside from the fact that the present construction does not use serendipity to reduce by one the degree of the internal DOFs, the method differs in several respects from that of \cite{Di-Pietro.Droniou.ea:26}.

  The first difference is the treatment of the load.
  The distinction goes beyond the polynomial degree of the reconstruction.
  In \cite{Di-Pietro.Droniou.ea:26}, the identity that connects the source to the discrete Hessian is obtained by composing adjoint-consistency estimates for two first-order operators.
  The estimate for the first discrete gradient contains an $h^{k}$ contribution; this term limits the final energy rate, even though the Hessian reconstruction itself contains higher-order information.
  A closer look shows that the serendipity reduction in \cite{Di-Pietro.Droniou.ea:26} actually prevents the extension of the discrete integration-by-parts identity \eqref{eq:hessian-potential} to the full space $\Poly{k}(T;\Real^{2\times2})$.
  Retaining the full cell DOFs, as we do here, is therefore crucial to this end; see \eqref{eq:potential-low-degree-extension}.
  Thanks to this extension, the load can be compared to the continuous problem at the Hessian level (see \cref{lem:consistency-error} below), where the approximation error is $O(h^{k+1})$.

  The second difference lies in the choice of the stabilization.
  Adapted to the present context, the stabilization proposed in \cite{Di-Pietro.Droniou.ea:26} reads
  \begin{equation}
    s_T(\underline{u}_T, \underline{v}_T)
    = h_T \sum_{E \in \ET} \int_E (H_T^{k} \underline{u}_T t_E - H_E^{k} \underline{u}_E) \cdot (H_T^{k} \underline{v}_T t_E - H_E^{k} \underline{v}_E).
    \label{eq:hessian-stabilization}
  \end{equation}
  This stabilization is also admissible for the present method when $k\ge1$, but the corresponding analysis would be longer since it does not directly provide the boundary differences needed in the estimate of the adjoint consistency error; see \cref{lem:consistency-error} below.
  It is, however, not admissible for $k=0$: both $H_E^0$ and $H_T^0$ depend only on the vertex gradients, leaving arbitrary interior vertex values uncontrolled.
  The stabilization \eqref{eq:hessian-stabilization} is tested numerically for $k\ge1$ in \cref{sec:discrete-problem:numerical-validation}, while the potential-based stabilization \eqref{eq:potential-stabilization} is used in \cref{sec:von-karman} for the method with $k = 1$.

  For polytopal discretizations of the de Rham complex, schemes based on mixed formulations typically rely on stabilized discrete $L^2$ products on each space of the complex.
  Bilinear forms obtained using discrete differential operators as arguments of these $L^2$ products therefore do not directly depend on the associated potential reconstructions (see, e.g., \cite{Beirao-da-Veiga.Brezzi.ea:17,Beirao-da-Veiga.Dassi.ea:22}).
  The failure of the Hessian-based stabilization for $k = 0$ suggests that, in this case, a direct dependence on the reconstructed potentials is necessary for stability.
\end{remark}

\begin{remark}[Comparison with the method of \cite{Brezzi.Marini:13}]
  To compare polynomial degrees, denote by $\ell$ the accuracy parameter called $k$ in \cite{Brezzi.Marini:13}.
For $\ell\ge2$, setting $k=\ell-2$ gives exactly the data in the space \eqref{eq:global-space}, up to scaling and the representation of moments by $L^2$ projections.
For $\ell\ge4$, a simple load construction using only cell moments corresponds, in the present notation, to
  \begin{equation*}
    \ell_h^{\mathrm{cell}}(\underline{v}_h)
    =\sum_{T\in\Th}\int_T f\,v_T,
  \end{equation*}
  which vanishes identically for $k=0,1$, since $\Poly{k-2}(T)=\{0\}$ at these degrees.
  While \cite[Section~4.5]{Brezzi.Marini:13} supplies special loads for $\ell=2,3,4$, the present construction uses a single load construction based on the Hessian potential for all $k\ge0$, including the vertex-only case, with the potential-based stabilization.
  Additionally, as already mentioned, our analysis does not require the uniform norm-equivalence assumption made in \cite[Assumption~H2]{Brezzi.Marini:13}, whose verification can be technically demanding.
  \label{rem:comparison-bm}
\end{remark}

\subsection{Numerical validation}
\label{sec:discrete-problem:numerical-validation}

\begin{figure}[!htbp]
  \centering
  \subcaptionbox{Triangular mesh}{
    \includegraphics[height=4cm]{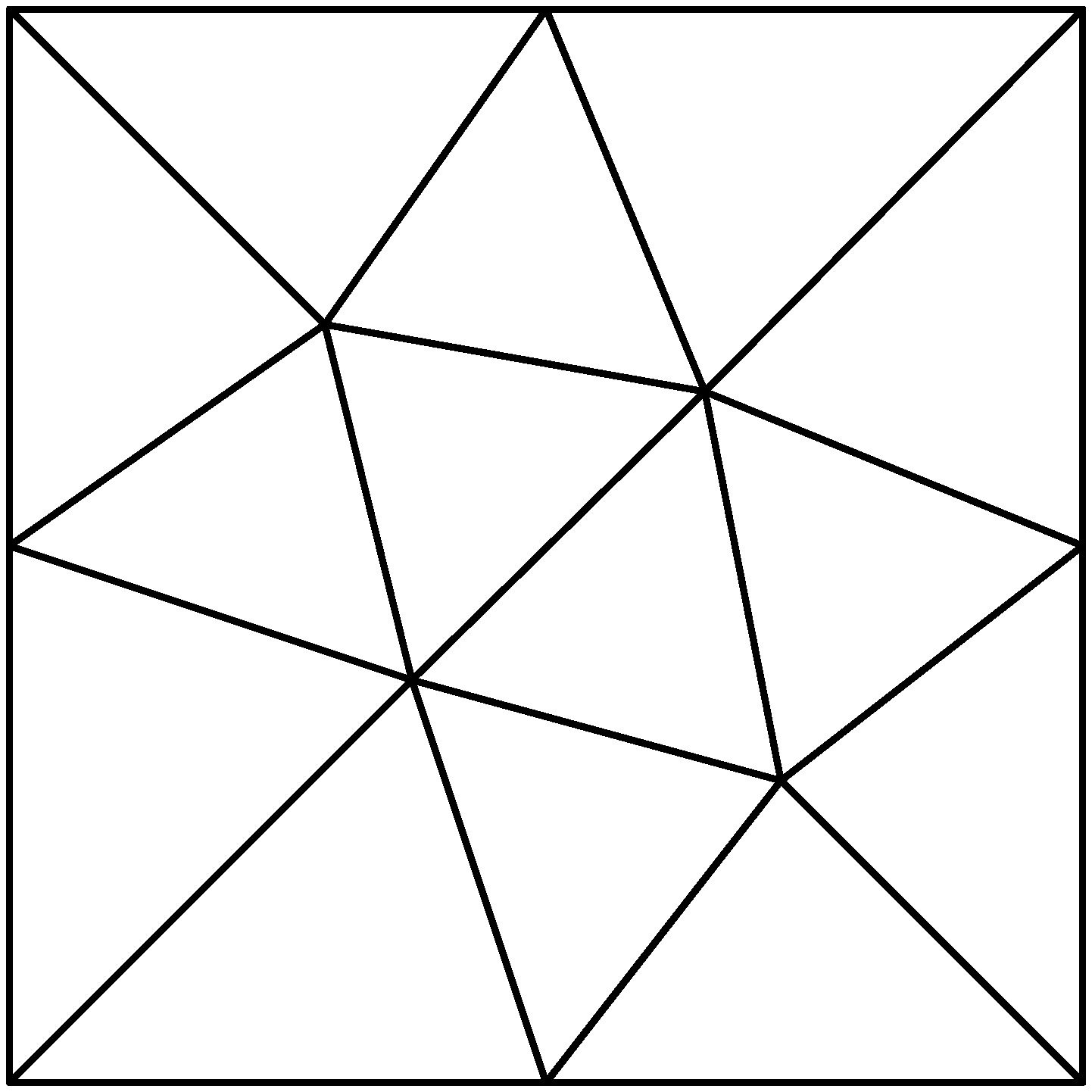}
  }
  \hspace{0.5cm}
  \subcaptionbox{Cartesian mesh}{
    \includegraphics[height=4cm]{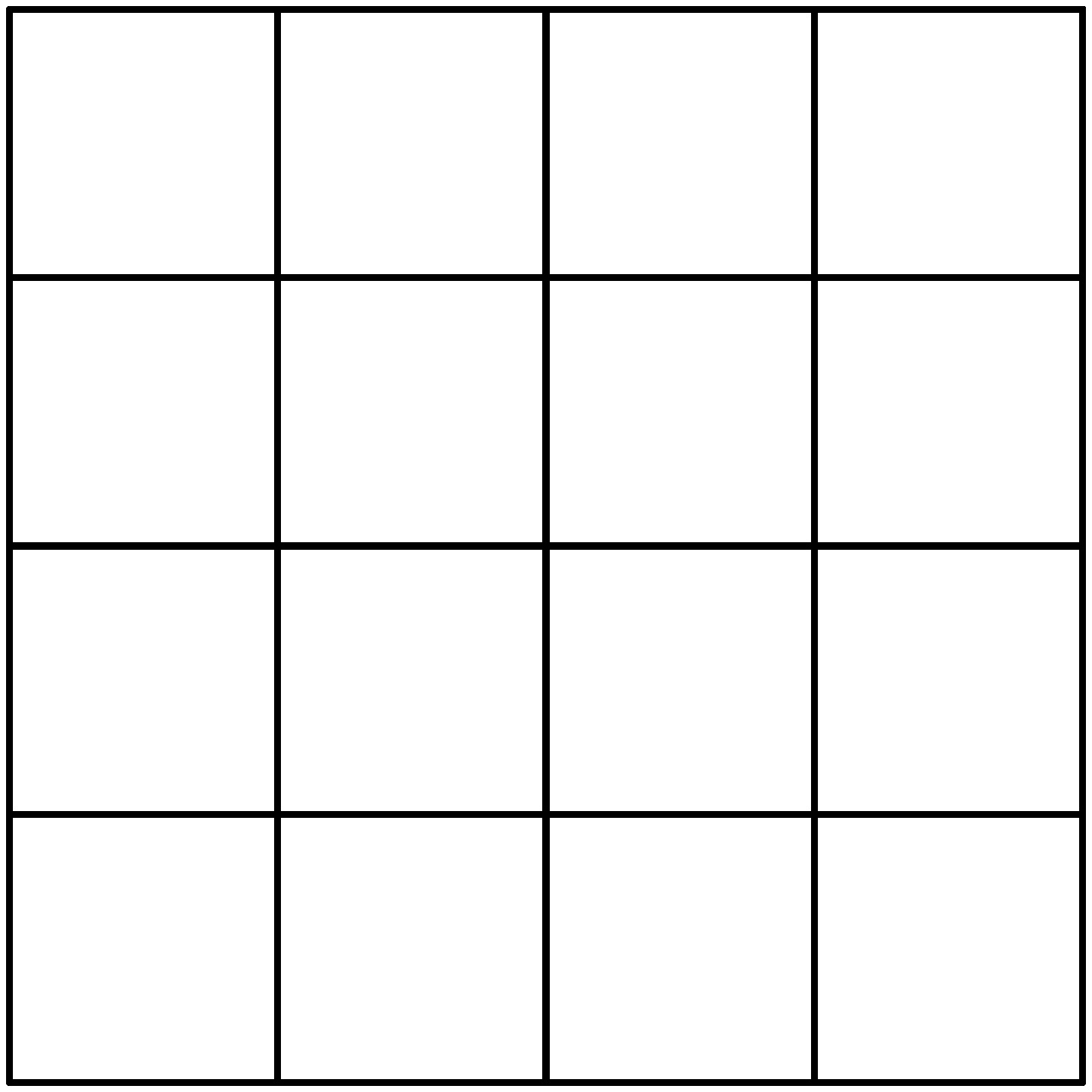}
  }
  \hspace{0.5cm}
  \subcaptionbox{Hexagonal mesh}{
    \includegraphics[height=4cm]{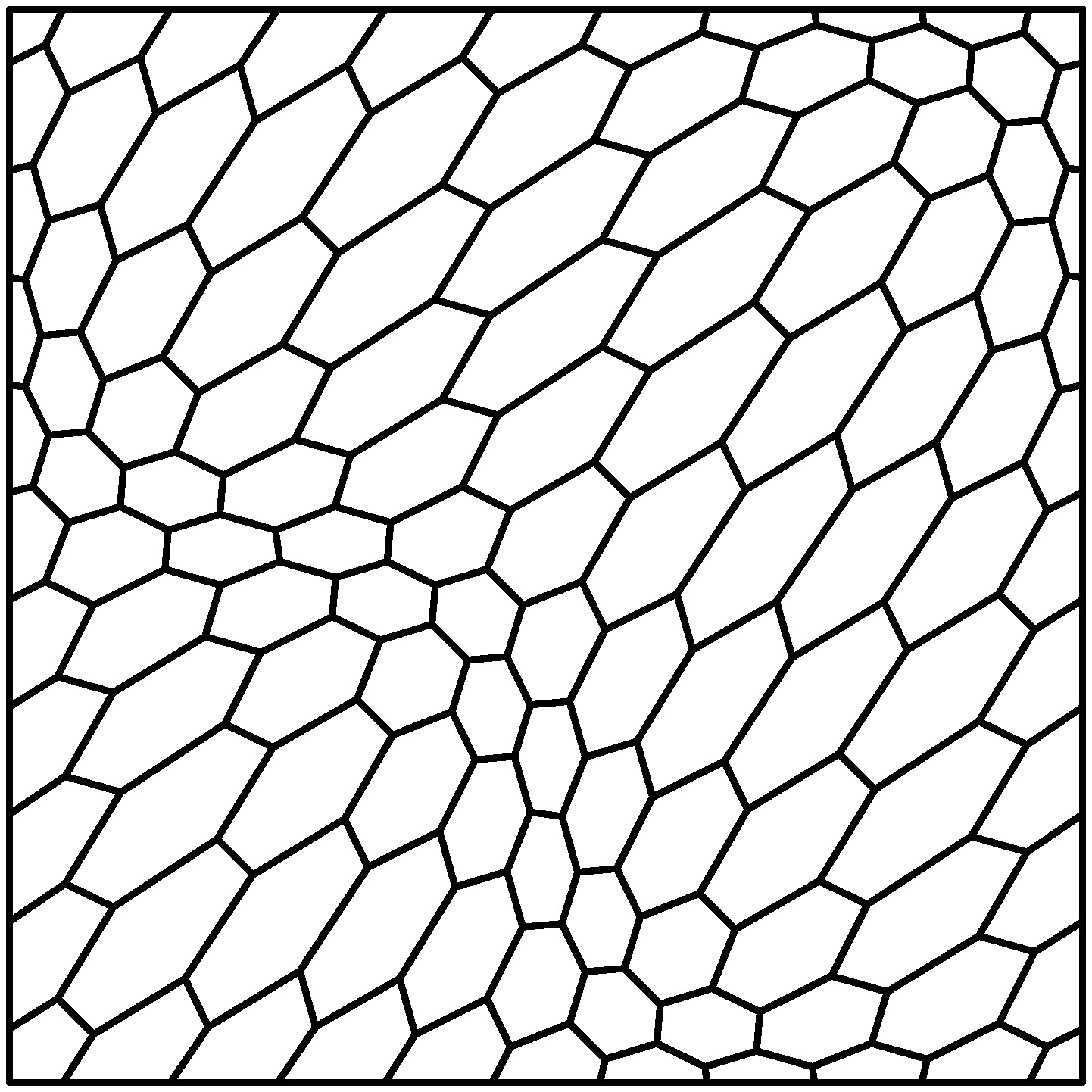}
  }
  \caption{Coarsest meshes used for the numerical validation of the error estimates in \cref{sec:discrete-problem:numerical-validation}.}
  \label{fig:convergence-test-meshes}
\end{figure}

To assess the convergence rate predicted by \cref{thm:energy-error-estimate}, we implemented the scheme \eqref{eq:discrete-problem} in the \polyrust library\footnote{\url{https://polyrust.org}} \cite{Di-Pietro:26} and ran the convergence test described in \cite[Section 5.3]{Di-Pietro.Droniou.ea:26} on mesh families obtained by refining the meshes depicted in \cref{fig:convergence-test-meshes}.
On the unit square $\Omega=(0,1)^2$, the manufactured solution and source are
\begin{equation}
  u(x,y)=\sin(\pi x)\sin(\pi y),
  \qquad
  f(x,y)=\Delta^2u(x,y)=4\pi^4\sin(\pi x)\sin(\pi y).
  \label{eq:manufactured-kirchhoff-love}
\end{equation}
We prescribe the exact boundary traces: $u=0$ on $\partial\Omega$, while the outward normal derivative is
\[
  \partial_n u=
  \begin{cases}
    -\pi\sin(\pi y) & \text{if $x=0$ or $x=1$}, \\
    -\pi\sin(\pi x) & \text{if $y=0$ or $y=1$}.
  \end{cases}
\]
Thus, this test involves a slight generalization of \eqref{eq:strong-problem} to accommodate nonhomogeneous normal-derivative data.
The boundary DOFs are prescribed as those of $\IHtwo{k}{h} u$, using exact vertex values and gradients and the edge $L^2$ projections specified by the interpolator.
We consider degrees $k=0,1,2,3$ with the potential-based stabilization \eqref{eq:potential-stabilization} and $k=1,2,3$ with the Hessian-based stabilization \eqref{eq:hessian-stabilization} on the triangular, Cartesian, and hexagonal mesh families depicted in \cref{fig:convergence-test-meshes}.

\begingroup
\tikzset{every node/.append style={scale=0.85}}
\begin{figure}[!htbp]
  \centering
  \subcaptionbox{Triangular mesh family}{
  \begin{minipage}{0.48\textwidth}
\noindent\begin{tikzpicture}
\begin{loglogaxis}[
  ymin=1e-8,
  ymax=1e1,
  enlarge y limits=false,
  width=\linewidth,
  height=0.78\linewidth,
  title={energy norm vs. $h$},
  title style={font=\small},
  label style={font=\scriptsize},
  tick label style={font=\scriptsize},
  grid=both,
  major grid style={gray!40},
  minor grid style={gray!20},
  legend pos=south east,
  legend style={font=\scriptsize, inner xsep=2pt, inner ysep=1pt},
  legend cell align={left}
]
\addplot+[blue, mark=*, thick] coordinates { (5.3033008588991082e-1,9.5477347523426093e-1) (2.6516504294495546e-1,4.6300854590201629e-1) (1.3258252147247773e-1,2.2735605397114023e-1) (6.6291260736238949e-2,1.1563466232087682e-1) (3.3145630368119572e-2,5.8677860792313678e-2) };
\addlegendentry{$k = 0$}
\draw[draw=none] (axis cs:5.3033008588991082e-1,9.5477347523426093e-1) -- (axis cs:2.6516504294495546e-1,4.6300854590201629e-1) node[midway, above, text=blue, fill=white, fill opacity=0.8, text opacity=1, inner sep=1pt, font=\scriptsize] {1.04};
\draw[draw=none] (axis cs:2.6516504294495546e-1,4.6300854590201629e-1) -- (axis cs:1.3258252147247773e-1,2.2735605397114023e-1) node[midway, above, text=blue, fill=white, fill opacity=0.8, text opacity=1, inner sep=1pt, font=\scriptsize] {1.03};
\draw[draw=none] (axis cs:1.3258252147247773e-1,2.2735605397114023e-1) -- (axis cs:6.6291260736238949e-2,1.1563466232087682e-1) node[midway, above, text=blue, fill=white, fill opacity=0.8, text opacity=1, inner sep=1pt, font=\scriptsize] {0.98};
\draw[draw=none] (axis cs:6.6291260736238949e-2,1.1563466232087682e-1) -- (axis cs:3.3145630368119572e-2,5.8677860792313678e-2) node[midway, above, text=blue, fill=white, fill opacity=0.8, text opacity=1, inner sep=1pt, font=\scriptsize] {0.98};
\addplot+[red, mark=square*, thick] coordinates { (5.3033008588991082e-1,7.4989452218480301e-1) (2.6516504294495546e-1,1.6612824218735489e-1) (1.3258252147247773e-1,3.0975829062538710e-2) (6.6291260736238949e-2,5.2204935175618221e-3) (3.3145630368119572e-2,8.5645739416316178e-4) };
\addlegendentry{$k = 1$}
\draw[draw=none] (axis cs:5.3033008588991082e-1,7.4989452218480301e-1) -- (axis cs:2.6516504294495546e-1,1.6612824218735489e-1) node[midway, above, text=red, fill=white, fill opacity=0.8, text opacity=1, inner sep=1pt, font=\scriptsize] {2.17};
\draw[draw=none] (axis cs:2.6516504294495546e-1,1.6612824218735489e-1) -- (axis cs:1.3258252147247773e-1,3.0975829062538710e-2) node[midway, above, text=red, fill=white, fill opacity=0.8, text opacity=1, inner sep=1pt, font=\scriptsize] {2.42};
\draw[draw=none] (axis cs:1.3258252147247773e-1,3.0975829062538710e-2) -- (axis cs:6.6291260736238949e-2,5.2204935175618221e-3) node[midway, above, text=red, fill=white, fill opacity=0.8, text opacity=1, inner sep=1pt, font=\scriptsize] {2.57};
\draw[draw=none] (axis cs:6.6291260736238949e-2,5.2204935175618221e-3) -- (axis cs:3.3145630368119572e-2,8.5645739416316178e-4) node[midway, above, text=red, fill=white, fill opacity=0.8, text opacity=1, inner sep=1pt, font=\scriptsize] {2.61};
\addplot+[brown, mark=o, thick] coordinates { (5.3033008588991082e-1,1.5398330634983948e-1) (2.6516504294495546e-1,1.9281632222919814e-2) (1.3258252147247773e-1,2.4613100766488626e-3) (6.6291260736238949e-2,3.1227522045136169e-4) (3.3145630368119572e-2,3.9333115504591302e-5) };
\addlegendentry{$k = 2$}
\draw[draw=none] (axis cs:5.3033008588991082e-1,1.5398330634983948e-1) -- (axis cs:2.6516504294495546e-1,1.9281632222919814e-2) node[midway, above, text=brown, fill=white, fill opacity=0.8, text opacity=1, inner sep=1pt, font=\scriptsize] {3.00};
\draw[draw=none] (axis cs:2.6516504294495546e-1,1.9281632222919814e-2) -- (axis cs:1.3258252147247773e-1,2.4613100766488626e-3) node[midway, above, text=brown, fill=white, fill opacity=0.8, text opacity=1, inner sep=1pt, font=\scriptsize] {2.97};
\draw[draw=none] (axis cs:1.3258252147247773e-1,2.4613100766488626e-3) -- (axis cs:6.6291260736238949e-2,3.1227522045136169e-4) node[midway, above, text=brown, fill=white, fill opacity=0.8, text opacity=1, inner sep=1pt, font=\scriptsize] {2.98};
\draw[draw=none] (axis cs:6.6291260736238949e-2,3.1227522045136169e-4) -- (axis cs:3.3145630368119572e-2,3.9333115504591302e-5) node[midway, above, text=brown, fill=white, fill opacity=0.8, text opacity=1, inner sep=1pt, font=\scriptsize] {2.99};
\addplot+[teal, mark=triangle*, thick] coordinates { (5.3033008588991082e-1,1.3008077595569279e-2) (2.6516504294495546e-1,7.9009504126308007e-4) (1.3258252147247773e-1,4.9934730231416032e-5) (6.6291260736238949e-2,3.1272436769566881e-6) (3.3145630368119572e-2,2.0283405487819562e-7) };
\addlegendentry{$k = 3$}
\draw[draw=none] (axis cs:5.3033008588991082e-1,1.3008077595569279e-2) -- (axis cs:2.6516504294495546e-1,7.9009504126308007e-4) node[midway, above, text=teal, fill=white, fill opacity=0.8, text opacity=1, inner sep=1pt, font=\scriptsize] {4.04};
\draw[draw=none] (axis cs:2.6516504294495546e-1,7.9009504126308007e-4) -- (axis cs:1.3258252147247773e-1,4.9934730231416032e-5) node[midway, above, text=teal, fill=white, fill opacity=0.8, text opacity=1, inner sep=1pt, font=\scriptsize] {3.98};
\draw[draw=none] (axis cs:1.3258252147247773e-1,4.9934730231416032e-5) -- (axis cs:6.6291260736238949e-2,3.1272436769566881e-6) node[midway, above, text=teal, fill=white, fill opacity=0.8, text opacity=1, inner sep=1pt, font=\scriptsize] {4.00};
\draw[draw=none] (axis cs:6.6291260736238949e-2,3.1272436769566881e-6) -- (axis cs:3.3145630368119572e-2,2.0283405487819562e-7) node[midway, above, text=teal, fill=white, fill opacity=0.8, text opacity=1, inner sep=1pt, font=\scriptsize] {3.95};
\end{loglogaxis}
\end{tikzpicture}
    \end{minipage}
    \begin{minipage}{0.48\textwidth}
\noindent\begin{tikzpicture}
\begin{loglogaxis}[
  ymin=1e-11,
  ymax=1e-1,
  enlarge y limits=false,
  width=\linewidth,
  height=0.78\linewidth,
  title={l2 potential norm vs. $h$},
  title style={font=\small},
  label style={font=\scriptsize},
  tick label style={font=\scriptsize},
  grid=both,
  major grid style={gray!40},
  minor grid style={gray!20},
  legend pos=south east,
  legend style={font=\scriptsize, inner xsep=2pt, inner ysep=1pt},
  legend cell align={left}
]
\addplot+[blue, mark=*, thick] coordinates { (5.3033008588991082e-1,2.4838491750328939e-2) (2.6516504294495546e-1,4.7184526622268633e-3) (1.3258252147247773e-1,1.1526691413898534e-3) (6.6291260736238949e-2,2.9000874080330392e-4) (3.3145630368119572e-2,7.3118745079422042e-5) };
\addlegendentry{$k = 0$}
\draw[draw=none] (axis cs:5.3033008588991082e-1,2.4838491750328939e-2) -- (axis cs:2.6516504294495546e-1,4.7184526622268633e-3) node[midway, above, text=blue, fill=white, fill opacity=0.8, text opacity=1, inner sep=1pt, font=\scriptsize] {2.40};
\draw[draw=none] (axis cs:2.6516504294495546e-1,4.7184526622268633e-3) -- (axis cs:1.3258252147247773e-1,1.1526691413898534e-3) node[midway, above, text=blue, fill=white, fill opacity=0.8, text opacity=1, inner sep=1pt, font=\scriptsize] {2.03};
\draw[draw=none] (axis cs:1.3258252147247773e-1,1.1526691413898534e-3) -- (axis cs:6.6291260736238949e-2,2.9000874080330392e-4) node[midway, above, text=blue, fill=white, fill opacity=0.8, text opacity=1, inner sep=1pt, font=\scriptsize] {1.99};
\draw[draw=none] (axis cs:6.6291260736238949e-2,2.9000874080330392e-4) -- (axis cs:3.3145630368119572e-2,7.3118745079422042e-5) node[midway, above, text=blue, fill=white, fill opacity=0.8, text opacity=1, inner sep=1pt, font=\scriptsize] {1.99};
\addplot+[red, mark=square*, thick] coordinates { (5.3033008588991082e-1,8.2758615159385027e-3) (2.6516504294495546e-1,5.9236446584087473e-4) (1.3258252147247773e-1,4.4896652342583270e-5) (6.6291260736238949e-2,3.1284704653863018e-6) (3.3145630368119572e-2,2.0317907149514628e-7) };
\addlegendentry{$k = 1$}
\draw[draw=none] (axis cs:5.3033008588991082e-1,8.2758615159385027e-3) -- (axis cs:2.6516504294495546e-1,5.9236446584087473e-4) node[midway, above, text=red, fill=white, fill opacity=0.8, text opacity=1, inner sep=1pt, font=\scriptsize] {3.80};
\draw[draw=none] (axis cs:2.6516504294495546e-1,5.9236446584087473e-4) -- (axis cs:1.3258252147247773e-1,4.4896652342583270e-5) node[midway, above, text=red, fill=white, fill opacity=0.8, text opacity=1, inner sep=1pt, font=\scriptsize] {3.72};
\draw[draw=none] (axis cs:1.3258252147247773e-1,4.4896652342583270e-5) -- (axis cs:6.6291260736238949e-2,3.1284704653863018e-6) node[midway, above, text=red, fill=white, fill opacity=0.8, text opacity=1, inner sep=1pt, font=\scriptsize] {3.84};
\draw[draw=none] (axis cs:6.6291260736238949e-2,3.1284704653863018e-6) -- (axis cs:3.3145630368119572e-2,2.0317907149514628e-7) node[midway, above, text=red, fill=white, fill opacity=0.8, text opacity=1, inner sep=1pt, font=\scriptsize] {3.94};
\addplot+[brown, mark=o, thick] coordinates { (5.3033008588991082e-1,4.4960652800998910e-4) (2.6516504294495546e-1,1.0997294305524605e-5) (1.3258252147247773e-1,2.8786608737118227e-7) (6.6291260736238949e-2,8.3686377797291738e-9) (3.3145630368119572e-2,2.9274958307872642e-10) };
\addlegendentry{$k = 2$}
\draw[draw=none] (axis cs:5.3033008588991082e-1,4.4960652800998910e-4) -- (axis cs:2.6516504294495546e-1,1.0997294305524605e-5) node[midway, above, text=brown, fill=white, fill opacity=0.8, text opacity=1, inner sep=1pt, font=\scriptsize] {5.35};
\draw[draw=none] (axis cs:2.6516504294495546e-1,1.0997294305524605e-5) -- (axis cs:1.3258252147247773e-1,2.8786608737118227e-7) node[midway, above, text=brown, fill=white, fill opacity=0.8, text opacity=1, inner sep=1pt, font=\scriptsize] {5.26};
\draw[draw=none] (axis cs:1.3258252147247773e-1,2.8786608737118227e-7) -- (axis cs:6.6291260736238949e-2,8.3686377797291738e-9) node[midway, above, text=brown, fill=white, fill opacity=0.8, text opacity=1, inner sep=1pt, font=\scriptsize] {5.10};
\draw[draw=none] (axis cs:6.6291260736238949e-2,8.3686377797291738e-9) -- (axis cs:3.3145630368119572e-2,2.9274958307872642e-10) node[midway, above, text=brown, fill=white, fill opacity=0.8, text opacity=1, inner sep=1pt, font=\scriptsize] {4.84};
\addplot+[teal, mark=triangle*, thick] coordinates { (5.3033008588991082e-1,1.7957127768833336e-5) (2.6516504294495546e-1,2.5831729921190953e-7) (1.3258252147247773e-1,3.3447767355543770e-9) (6.6291260736238949e-2,7.4750684558457156e-11) (3.3145630368119572e-2,1.3926758765044121e-9) };
\addlegendentry{$k = 3$}
\draw[draw=none] (axis cs:5.3033008588991082e-1,1.7957127768833336e-5) -- (axis cs:2.6516504294495546e-1,2.5831729921190953e-7) node[midway, above, text=teal, fill=white, fill opacity=0.8, text opacity=1, inner sep=1pt, font=\scriptsize] {6.12};
\draw[draw=none] (axis cs:2.6516504294495546e-1,2.5831729921190953e-7) -- (axis cs:1.3258252147247773e-1,3.3447767355543770e-9) node[midway, above, text=teal, fill=white, fill opacity=0.8, text opacity=1, inner sep=1pt, font=\scriptsize] {6.27};
\draw[draw=none] (axis cs:1.3258252147247773e-1,3.3447767355543770e-9) -- (axis cs:6.6291260736238949e-2,7.4750684558457156e-11) node[midway, above, text=teal, fill=white, fill opacity=0.8, text opacity=1, inner sep=1pt, font=\scriptsize] {5.48};
\draw[draw=none] (axis cs:6.6291260736238949e-2,7.4750684558457156e-11) -- (axis cs:3.3145630368119572e-2,1.3926758765044121e-9) node[midway, above, text=teal, fill=white, fill opacity=0.8, text opacity=1, inner sep=1pt, font=\scriptsize] {-4.22};
\end{loglogaxis}
\end{tikzpicture}
    \end{minipage}
  }
  \\
  \subcaptionbox{Cartesian mesh family}{
    \begin{minipage}{0.48\textwidth}
\noindent\begin{tikzpicture}
\begin{loglogaxis}[
  ymin=1e-8,
  ymax=1e1,
  enlarge y limits=false,
  width=\linewidth,
  height=0.78\linewidth,
  title={energy norm vs. $h$},
  title style={font=\small},
  label style={font=\scriptsize},
  tick label style={font=\scriptsize},
  grid=both,
  major grid style={gray!40},
  minor grid style={gray!20},
  legend pos=south east,
  legend style={font=\scriptsize, inner xsep=2pt, inner ysep=1pt},
  legend cell align={left}
]
\addplot+[blue, mark=*, thick] coordinates { (3.5355339059421481e-1,2.7181091901154888e0) (1.7677669529722873e-1,1.3491509191074791e0) (8.8388347648650084e-2,6.7265720686553121e-1) (4.4194173824334396e-2,3.3606908599282209e-1) (2.2097086912169481e-2,1.6800152749385627e-1) };
\addlegendentry{$k = 0$}
\draw[draw=none] (axis cs:3.5355339059421481e-1,2.7181091901154888e0) -- (axis cs:1.7677669529722873e-1,1.3491509191074791e0) node[midway, above, text=blue, fill=white, fill opacity=0.8, text opacity=1, inner sep=1pt, font=\scriptsize] {1.01};
\draw[draw=none] (axis cs:1.7677669529722873e-1,1.3491509191074791e0) -- (axis cs:8.8388347648650084e-2,6.7265720686553121e-1) node[midway, above, text=blue, fill=white, fill opacity=0.8, text opacity=1, inner sep=1pt, font=\scriptsize] {1.00};
\draw[draw=none] (axis cs:8.8388347648650084e-2,6.7265720686553121e-1) -- (axis cs:4.4194173824334396e-2,3.3606908599282209e-1) node[midway, above, text=blue, fill=white, fill opacity=0.8, text opacity=1, inner sep=1pt, font=\scriptsize] {1.00};
\draw[draw=none] (axis cs:4.4194173824334396e-2,3.3606908599282209e-1) -- (axis cs:2.2097086912169481e-2,1.6800152749385627e-1) node[midway, above, text=blue, fill=white, fill opacity=0.8, text opacity=1, inner sep=1pt, font=\scriptsize] {1.00};
\addplot+[red, mark=square*, thick] coordinates { (3.5355339059421481e-1,5.8974755866329769e-2) (1.7677669529722873e-1,8.6135185055726000e-3) (8.8388347648650084e-2,1.1496942705828651e-3) (4.4194173824334396e-2,1.4778275903455721e-4) (2.2097086912169481e-2,1.8707391757827739e-5) };
\addlegendentry{$k = 1$}
\draw[draw=none] (axis cs:3.5355339059421481e-1,5.8974755866329769e-2) -- (axis cs:1.7677669529722873e-1,8.6135185055726000e-3) node[midway, above, text=red, fill=white, fill opacity=0.8, text opacity=1, inner sep=1pt, font=\scriptsize] {2.78};
\draw[draw=none] (axis cs:1.7677669529722873e-1,8.6135185055726000e-3) -- (axis cs:8.8388347648650084e-2,1.1496942705828651e-3) node[midway, above, text=red, fill=white, fill opacity=0.8, text opacity=1, inner sep=1pt, font=\scriptsize] {2.91};
\draw[draw=none] (axis cs:8.8388347648650084e-2,1.1496942705828651e-3) -- (axis cs:4.4194173824334396e-2,1.4778275903455721e-4) node[midway, above, text=red, fill=white, fill opacity=0.8, text opacity=1, inner sep=1pt, font=\scriptsize] {2.96};
\draw[draw=none] (axis cs:4.4194173824334396e-2,1.4778275903455721e-4) -- (axis cs:2.2097086912169481e-2,1.8707391757827739e-5) node[midway, above, text=red, fill=white, fill opacity=0.8, text opacity=1, inner sep=1pt, font=\scriptsize] {2.98};
\addplot+[brown, mark=o, thick] coordinates { (3.5355339059421481e-1,1.9424005175958083e-2) (1.7677669529722873e-1,2.6714696956190668e-3) (8.8388347648650084e-2,3.5245655806468579e-4) (4.4194173824334396e-2,4.5269114078361764e-5) (2.2097086912169481e-2,5.7349995849077317e-6) };
\addlegendentry{$k = 2$}
\draw[draw=none] (axis cs:3.5355339059421481e-1,1.9424005175958083e-2) -- (axis cs:1.7677669529722873e-1,2.6714696956190668e-3) node[midway, above, text=brown, fill=white, fill opacity=0.8, text opacity=1, inner sep=1pt, font=\scriptsize] {2.86};
\draw[draw=none] (axis cs:1.7677669529722873e-1,2.6714696956190668e-3) -- (axis cs:8.8388347648650084e-2,3.5245655806468579e-4) node[midway, above, text=brown, fill=white, fill opacity=0.8, text opacity=1, inner sep=1pt, font=\scriptsize] {2.92};
\draw[draw=none] (axis cs:8.8388347648650084e-2,3.5245655806468579e-4) -- (axis cs:4.4194173824334396e-2,4.5269114078361764e-5) node[midway, above, text=brown, fill=white, fill opacity=0.8, text opacity=1, inner sep=1pt, font=\scriptsize] {2.96};
\draw[draw=none] (axis cs:4.4194173824334396e-2,4.5269114078361764e-5) -- (axis cs:2.2097086912169481e-2,5.7349995849077317e-6) node[midway, above, text=brown, fill=white, fill opacity=0.8, text opacity=1, inner sep=1pt, font=\scriptsize] {2.98};
\addplot+[teal, mark=triangle*, thick] coordinates { (3.5355339059421481e-1,4.0157463136445502e-3) (1.7677669529722873e-1,2.6032295505902189e-4) (8.8388347648650084e-2,1.6403687394129399e-5) (4.4194173824334396e-2,1.0274757289943739e-6) (2.2097086912169481e-2,6.5411700737728746e-8) };
\addlegendentry{$k = 3$}
\draw[draw=none] (axis cs:3.5355339059421481e-1,4.0157463136445502e-3) -- (axis cs:1.7677669529722873e-1,2.6032295505902189e-4) node[midway, above, text=teal, fill=white, fill opacity=0.8, text opacity=1, inner sep=1pt, font=\scriptsize] {3.95};
\draw[draw=none] (axis cs:1.7677669529722873e-1,2.6032295505902189e-4) -- (axis cs:8.8388347648650084e-2,1.6403687394129399e-5) node[midway, above, text=teal, fill=white, fill opacity=0.8, text opacity=1, inner sep=1pt, font=\scriptsize] {3.99};
\draw[draw=none] (axis cs:8.8388347648650084e-2,1.6403687394129399e-5) -- (axis cs:4.4194173824334396e-2,1.0274757289943739e-6) node[midway, above, text=teal, fill=white, fill opacity=0.8, text opacity=1, inner sep=1pt, font=\scriptsize] {4.00};
\draw[draw=none] (axis cs:4.4194173824334396e-2,1.0274757289943739e-6) -- (axis cs:2.2097086912169481e-2,6.5411700737728746e-8) node[midway, above, text=teal, fill=white, fill opacity=0.8, text opacity=1, inner sep=1pt, font=\scriptsize] {3.97};
\end{loglogaxis}
\end{tikzpicture}
    \end{minipage}
    \begin{minipage}{0.48\textwidth}
\noindent\begin{tikzpicture}
\begin{loglogaxis}[
  ymin=1e-11,
  ymax=1e0,
  enlarge y limits=false,
  width=\linewidth,
  height=0.78\linewidth,
  title={l2 potential norm vs. $h$},
  title style={font=\small},
  label style={font=\scriptsize},
  tick label style={font=\scriptsize},
  grid=both,
  major grid style={gray!40},
  minor grid style={gray!20},
  legend pos=south east,
  legend style={font=\scriptsize, inner xsep=2pt, inner ysep=1pt},
  legend cell align={left}
]
\addplot+[blue, mark=*, thick] coordinates { (3.5355339059421481e-1,7.4983275578525119e-2) (1.7677669529722873e-1,1.9440352494051669e-2) (8.8388347648650084e-2,4.8954858572773366e-3) (4.4194173824334396e-2,1.2257832269347458e-3) (2.2097086912169481e-2,3.0655640099700284e-4) };
\addlegendentry{$k = 0$}
\draw[draw=none] (axis cs:3.5355339059421481e-1,7.4983275578525119e-2) -- (axis cs:1.7677669529722873e-1,1.9440352494051669e-2) node[midway, above, text=blue, fill=white, fill opacity=0.8, text opacity=1, inner sep=1pt, font=\scriptsize] {1.95};
\draw[draw=none] (axis cs:1.7677669529722873e-1,1.9440352494051669e-2) -- (axis cs:8.8388347648650084e-2,4.8954858572773366e-3) node[midway, above, text=blue, fill=white, fill opacity=0.8, text opacity=1, inner sep=1pt, font=\scriptsize] {1.99};
\draw[draw=none] (axis cs:8.8388347648650084e-2,4.8954858572773366e-3) -- (axis cs:4.4194173824334396e-2,1.2257832269347458e-3) node[midway, above, text=blue, fill=white, fill opacity=0.8, text opacity=1, inner sep=1pt, font=\scriptsize] {2.00};
\draw[draw=none] (axis cs:4.4194173824334396e-2,1.2257832269347458e-3) -- (axis cs:2.2097086912169481e-2,3.0655640099700284e-4) node[midway, above, text=blue, fill=white, fill opacity=0.8, text opacity=1, inner sep=1pt, font=\scriptsize] {2.00};
\addplot+[red, mark=square*, thick] coordinates { (3.5355339059421481e-1,7.8929721534067498e-4) (1.7677669529722873e-1,4.9688618205499002e-5) (8.8388347648650084e-2,3.0874792424075346e-6) (4.4194173824334396e-2,1.9246600365712817e-7) (2.2097086912169481e-2,1.2016809762822955e-8) };
\addlegendentry{$k = 1$}
\draw[draw=none] (axis cs:3.5355339059421481e-1,7.8929721534067498e-4) -- (axis cs:1.7677669529722873e-1,4.9688618205499002e-5) node[midway, above, text=red, fill=white, fill opacity=0.8, text opacity=1, inner sep=1pt, font=\scriptsize] {3.99};
\draw[draw=none] (axis cs:1.7677669529722873e-1,4.9688618205499002e-5) -- (axis cs:8.8388347648650084e-2,3.0874792424075346e-6) node[midway, above, text=red, fill=white, fill opacity=0.8, text opacity=1, inner sep=1pt, font=\scriptsize] {4.01};
\draw[draw=none] (axis cs:8.8388347648650084e-2,3.0874792424075346e-6) -- (axis cs:4.4194173824334396e-2,1.9246600365712817e-7) node[midway, above, text=red, fill=white, fill opacity=0.8, text opacity=1, inner sep=1pt, font=\scriptsize] {4.00};
\draw[draw=none] (axis cs:4.4194173824334396e-2,1.9246600365712817e-7) -- (axis cs:2.2097086912169481e-2,1.2016809762822955e-8) node[midway, above, text=red, fill=white, fill opacity=0.8, text opacity=1, inner sep=1pt, font=\scriptsize] {4.00};
\addplot+[brown, mark=o, thick] coordinates { (3.5355339059421481e-1,2.5846973126558257e-5) (1.7677669529722873e-1,8.3944018784207701e-7) (8.8388347648650084e-2,2.7753727899386165e-8) (4.4194173824334396e-2,8.9503259485079093e-10) (2.2097086912169481e-2,6.8325873463579545e-11) };
\addlegendentry{$k = 2$}
\draw[draw=none] (axis cs:3.5355339059421481e-1,2.5846973126558257e-5) -- (axis cs:1.7677669529722873e-1,8.3944018784207701e-7) node[midway, above, text=brown, fill=white, fill opacity=0.8, text opacity=1, inner sep=1pt, font=\scriptsize] {4.94};
\draw[draw=none] (axis cs:1.7677669529722873e-1,8.3944018784207701e-7) -- (axis cs:8.8388347648650084e-2,2.7753727899386165e-8) node[midway, above, text=brown, fill=white, fill opacity=0.8, text opacity=1, inner sep=1pt, font=\scriptsize] {4.92};
\draw[draw=none] (axis cs:8.8388347648650084e-2,2.7753727899386165e-8) -- (axis cs:4.4194173824334396e-2,8.9503259485079093e-10) node[midway, above, text=brown, fill=white, fill opacity=0.8, text opacity=1, inner sep=1pt, font=\scriptsize] {4.95};
\draw[draw=none] (axis cs:4.4194173824334396e-2,8.9503259485079093e-10) -- (axis cs:2.2097086912169481e-2,6.8325873463579545e-11) node[midway, above, text=brown, fill=white, fill opacity=0.8, text opacity=1, inner sep=1pt, font=\scriptsize] {3.71};
\addplot+[teal, mark=triangle*, thick] coordinates { (3.5355339059421481e-1,2.9248734305263561e-6) (1.7677669529722873e-1,4.8420341347917955e-8) (8.8388347648650084e-2,7.6949707874422123e-10) (4.4194173824334396e-2,4.9376471358508072e-11) (2.2097086912169481e-2,9.5407182445197509e-11) };
\addlegendentry{$k = 3$}
\draw[draw=none] (axis cs:3.5355339059421481e-1,2.9248734305263561e-6) -- (axis cs:1.7677669529722873e-1,4.8420341347917955e-8) node[midway, above, text=teal, fill=white, fill opacity=0.8, text opacity=1, inner sep=1pt, font=\scriptsize] {5.92};
\draw[draw=none] (axis cs:1.7677669529722873e-1,4.8420341347917955e-8) -- (axis cs:8.8388347648650084e-2,7.6949707874422123e-10) node[midway, above, text=teal, fill=white, fill opacity=0.8, text opacity=1, inner sep=1pt, font=\scriptsize] {5.98};
\draw[draw=none] (axis cs:8.8388347648650084e-2,7.6949707874422123e-10) -- (axis cs:4.4194173824334396e-2,4.9376471358508072e-11) node[midway, above, text=teal, fill=white, fill opacity=0.8, text opacity=1, inner sep=1pt, font=\scriptsize] {3.96};
\draw[draw=none] (axis cs:4.4194173824334396e-2,4.9376471358508072e-11) -- (axis cs:2.2097086912169481e-2,9.5407182445197509e-11) node[midway, above, text=teal, fill=white, fill opacity=0.8, text opacity=1, inner sep=1pt, font=\scriptsize] {-0.95};
\end{loglogaxis}
\end{tikzpicture}
    \end{minipage}
  }
  \\
  \subcaptionbox{Hexagonal mesh family}{
    \begin{minipage}{0.48\textwidth}
\noindent\begin{tikzpicture}
\begin{loglogaxis}[
  ymin=1e-7,
  ymax=1e1,
  enlarge y limits=false,
  width=\linewidth,
  height=0.78\linewidth,
  title={energy norm vs. $h$},
  title style={font=\small},
  label style={font=\scriptsize},
  tick label style={font=\scriptsize},
  grid=both,
  major grid style={gray!40},
  minor grid style={gray!20},
  legend pos=south east,
  legend style={font=\scriptsize, inner xsep=2pt, inner ysep=1pt},
  legend cell align={left}
]
\addplot+[blue, mark=*, thick] coordinates { (2.4141220176769076e-1,1.4545190449197398e0) (1.2971299742290118e-1,7.7576934311085222e-1) (6.5736358782959317e-2,3.9791698040671458e-1) (3.2979939347191366e-2,2.0101142814883552e-1) (1.6503995713160725e-2,1.0095258124158042e-1) };
\addlegendentry{$k = 0$}
\draw[draw=none] (axis cs:2.4141220176769076e-1,1.4545190449197398e0) -- (axis cs:1.2971299742290118e-1,7.7576934311085222e-1) node[midway, above, text=blue, fill=white, fill opacity=0.8, text opacity=1, inner sep=1pt, font=\scriptsize] {1.01};
\draw[draw=none] (axis cs:1.2971299742290118e-1,7.7576934311085222e-1) -- (axis cs:6.5736358782959317e-2,3.9791698040671458e-1) node[midway, above, text=blue, fill=white, fill opacity=0.8, text opacity=1, inner sep=1pt, font=\scriptsize] {0.98};
\draw[draw=none] (axis cs:6.5736358782959317e-2,3.9791698040671458e-1) -- (axis cs:3.2979939347191366e-2,2.0101142814883552e-1) node[midway, above, text=blue, fill=white, fill opacity=0.8, text opacity=1, inner sep=1pt, font=\scriptsize] {0.99};
\draw[draw=none] (axis cs:3.2979939347191366e-2,2.0101142814883552e-1) -- (axis cs:1.6503995713160725e-2,1.0095258124158042e-1) node[midway, above, text=blue, fill=white, fill opacity=0.8, text opacity=1, inner sep=1pt, font=\scriptsize] {0.99};
\addplot+[red, mark=square*, thick] coordinates { (2.4141220176769076e-1,1.3407673725136143e-1) (1.2971299742290118e-1,4.2135073033526899e-2) (6.5736358782959317e-2,1.1254303252917449e-2) (3.2979939347191366e-2,2.8722033703083003e-3) (1.6503995713160725e-2,7.2322939244990475e-4) };
\addlegendentry{$k = 1$}
\draw[draw=none] (axis cs:2.4141220176769076e-1,1.3407673725136143e-1) -- (axis cs:1.2971299742290118e-1,4.2135073033526899e-2) node[midway, above, text=red, fill=white, fill opacity=0.8, text opacity=1, inner sep=1pt, font=\scriptsize] {1.86};
\draw[draw=none] (axis cs:1.2971299742290118e-1,4.2135073033526899e-2) -- (axis cs:6.5736358782959317e-2,1.1254303252917449e-2) node[midway, above, text=red, fill=white, fill opacity=0.8, text opacity=1, inner sep=1pt, font=\scriptsize] {1.94};
\draw[draw=none] (axis cs:6.5736358782959317e-2,1.1254303252917449e-2) -- (axis cs:3.2979939347191366e-2,2.8722033703083003e-3) node[midway, above, text=red, fill=white, fill opacity=0.8, text opacity=1, inner sep=1pt, font=\scriptsize] {1.98};
\draw[draw=none] (axis cs:3.2979939347191366e-2,2.8722033703083003e-3) -- (axis cs:1.6503995713160725e-2,7.2322939244990475e-4) node[midway, above, text=red, fill=white, fill opacity=0.8, text opacity=1, inner sep=1pt, font=\scriptsize] {1.99};
\addplot+[brown, mark=o, thick] coordinates { (2.4141220176769076e-1,1.0173065834689955e-2) (1.2971299742290118e-1,1.7434316477190016e-3) (6.5736358782959317e-2,2.3862940865974010e-4) (3.2979939347191366e-2,3.0695011997040260e-5) (1.6503995713160725e-2,3.8802312493991125e-6) };
\addlegendentry{$k = 2$}
\draw[draw=none] (axis cs:2.4141220176769076e-1,1.0173065834689955e-2) -- (axis cs:1.2971299742290118e-1,1.7434316477190016e-3) node[midway, above, text=brown, fill=white, fill opacity=0.8, text opacity=1, inner sep=1pt, font=\scriptsize] {2.84};
\draw[draw=none] (axis cs:1.2971299742290118e-1,1.7434316477190016e-3) -- (axis cs:6.5736358782959317e-2,2.3862940865974010e-4) node[midway, above, text=brown, fill=white, fill opacity=0.8, text opacity=1, inner sep=1pt, font=\scriptsize] {2.93};
\draw[draw=none] (axis cs:6.5736358782959317e-2,2.3862940865974010e-4) -- (axis cs:3.2979939347191366e-2,3.0695011997040260e-5) node[midway, above, text=brown, fill=white, fill opacity=0.8, text opacity=1, inner sep=1pt, font=\scriptsize] {2.97};
\draw[draw=none] (axis cs:3.2979939347191366e-2,3.0695011997040260e-5) -- (axis cs:1.6503995713160725e-2,3.8802312493991125e-6) node[midway, above, text=brown, fill=white, fill opacity=0.8, text opacity=1, inner sep=1pt, font=\scriptsize] {2.99};
\addplot+[teal, mark=triangle*, thick] coordinates { (2.4141220176769076e-1,6.9654080108716963e-4) (1.2971299742290118e-1,6.4515170791025082e-5) (6.5736358782959317e-2,4.5622544741291515e-6) (3.2979939347191366e-2,3.0345684770948574e-7) (1.6503995713160725e-2,3.2097252134112667e-7) };
\addlegendentry{$k = 3$}
\draw[draw=none] (axis cs:2.4141220176769076e-1,6.9654080108716963e-4) -- (axis cs:1.2971299742290118e-1,6.4515170791025082e-5) node[midway, above, text=teal, fill=white, fill opacity=0.8, text opacity=1, inner sep=1pt, font=\scriptsize] {3.83};
\draw[draw=none] (axis cs:1.2971299742290118e-1,6.4515170791025082e-5) -- (axis cs:6.5736358782959317e-2,4.5622544741291515e-6) node[midway, above, text=teal, fill=white, fill opacity=0.8, text opacity=1, inner sep=1pt, font=\scriptsize] {3.90};
\draw[draw=none] (axis cs:6.5736358782959317e-2,4.5622544741291515e-6) -- (axis cs:3.2979939347191366e-2,3.0345684770948574e-7) node[midway, above, text=teal, fill=white, fill opacity=0.8, text opacity=1, inner sep=1pt, font=\scriptsize] {3.93};
\draw[draw=none] (axis cs:3.2979939347191366e-2,3.0345684770948574e-7) -- (axis cs:1.6503995713160725e-2,3.2097252134112667e-7) node[midway, above, text=teal, fill=white, fill opacity=0.8, text opacity=1, inner sep=1pt, font=\scriptsize] {-0.08};
\end{loglogaxis}
\end{tikzpicture}
    \end{minipage}
    \begin{minipage}{0.48\textwidth}
\noindent\begin{tikzpicture}
\begin{loglogaxis}[
  ymin=1e-11,
  ymax=1e-1,
  enlarge y limits=false,
  width=\linewidth,
  height=0.78\linewidth,
  title={l2 potential norm vs. $h$},
  title style={font=\small},
  label style={font=\scriptsize},
  tick label style={font=\scriptsize},
  grid=both,
  major grid style={gray!40},
  minor grid style={gray!20},
  legend pos=south east,
  legend style={font=\scriptsize, inner xsep=2pt, inner ysep=1pt},
  legend cell align={left}
]
\addplot+[blue, mark=*, thick] coordinates { (2.4141220176769076e-1,1.5929717274728847e-2) (1.2971299742290118e-1,4.5381743717795428e-3) (6.5736358782959317e-2,1.1973071718306186e-3) (3.2979939347191366e-2,3.0598912511518980e-4) (1.6503995713160725e-2,7.7233397394379607e-5) };
\addlegendentry{$k = 0$}
\draw[draw=none] (axis cs:2.4141220176769076e-1,1.5929717274728847e-2) -- (axis cs:1.2971299742290118e-1,4.5381743717795428e-3) node[midway, above, text=blue, fill=white, fill opacity=0.8, text opacity=1, inner sep=1pt, font=\scriptsize] {2.02};
\draw[draw=none] (axis cs:1.2971299742290118e-1,4.5381743717795428e-3) -- (axis cs:6.5736358782959317e-2,1.1973071718306186e-3) node[midway, above, text=blue, fill=white, fill opacity=0.8, text opacity=1, inner sep=1pt, font=\scriptsize] {1.96};
\draw[draw=none] (axis cs:6.5736358782959317e-2,1.1973071718306186e-3) -- (axis cs:3.2979939347191366e-2,3.0598912511518980e-4) node[midway, above, text=blue, fill=white, fill opacity=0.8, text opacity=1, inner sep=1pt, font=\scriptsize] {1.98};
\draw[draw=none] (axis cs:3.2979939347191366e-2,3.0598912511518980e-4) -- (axis cs:1.6503995713160725e-2,7.7233397394379607e-5) node[midway, above, text=blue, fill=white, fill opacity=0.8, text opacity=1, inner sep=1pt, font=\scriptsize] {1.99};
\addplot+[red, mark=square*, thick] coordinates { (2.4141220176769076e-1,1.0284534564954570e-4) (1.2971299742290118e-1,1.1333175808499249e-5) (6.5736358782959317e-2,8.7215795998732866e-7) (3.2979939347191366e-2,5.8879138372214777e-8) (1.6503995713160725e-2,3.8299621857750724e-9) };
\addlegendentry{$k = 1$}
\draw[draw=none] (axis cs:2.4141220176769076e-1,1.0284534564954570e-4) -- (axis cs:1.2971299742290118e-1,1.1333175808499249e-5) node[midway, above, text=red, fill=white, fill opacity=0.8, text opacity=1, inner sep=1pt, font=\scriptsize] {3.55};
\draw[draw=none] (axis cs:1.2971299742290118e-1,1.1333175808499249e-5) -- (axis cs:6.5736358782959317e-2,8.7215795998732866e-7) node[midway, above, text=red, fill=white, fill opacity=0.8, text opacity=1, inner sep=1pt, font=\scriptsize] {3.77};
\draw[draw=none] (axis cs:6.5736358782959317e-2,8.7215795998732866e-7) -- (axis cs:3.2979939347191366e-2,5.8879138372214777e-8) node[midway, above, text=red, fill=white, fill opacity=0.8, text opacity=1, inner sep=1pt, font=\scriptsize] {3.91};
\draw[draw=none] (axis cs:3.2979939347191366e-2,5.8879138372214777e-8) -- (axis cs:1.6503995713160725e-2,3.8299621857750724e-9) node[midway, above, text=red, fill=white, fill opacity=0.8, text opacity=1, inner sep=1pt, font=\scriptsize] {3.95};
\addplot+[brown, mark=o, thick] coordinates { (2.4141220176769076e-1,2.7849940120454784e-6) (1.2971299742290118e-1,1.2659653586202679e-7) (6.5736358782959317e-2,4.3719480671991266e-9) (3.2979939347191366e-2,3.6978575314898654e-10) (1.6503995713160725e-2,3.9931262294886599e-9) };
\addlegendentry{$k = 2$}
\draw[draw=none] (axis cs:2.4141220176769076e-1,2.7849940120454784e-6) -- (axis cs:1.2971299742290118e-1,1.2659653586202679e-7) node[midway, above, text=brown, fill=white, fill opacity=0.8, text opacity=1, inner sep=1pt, font=\scriptsize] {4.98};
\draw[draw=none] (axis cs:1.2971299742290118e-1,1.2659653586202679e-7) -- (axis cs:6.5736358782959317e-2,4.3719480671991266e-9) node[midway, above, text=brown, fill=white, fill opacity=0.8, text opacity=1, inner sep=1pt, font=\scriptsize] {4.95};
\draw[draw=none] (axis cs:6.5736358782959317e-2,4.3719480671991266e-9) -- (axis cs:3.2979939347191366e-2,3.6978575314898654e-10) node[midway, above, text=brown, fill=white, fill opacity=0.8, text opacity=1, inner sep=1pt, font=\scriptsize] {3.58};
\draw[draw=none] (axis cs:3.2979939347191366e-2,3.6978575314898654e-10) -- (axis cs:1.6503995713160725e-2,3.9931262294886599e-9) node[midway, above, text=brown, fill=white, fill opacity=0.8, text opacity=1, inner sep=1pt, font=\scriptsize] {-3.44};
\addplot+[teal, mark=triangle*, thick] coordinates { (2.4141220176769076e-1,9.7524046654869218e-8) (1.2971299742290118e-1,2.2879697357593217e-9) (6.5736358782959317e-2,7.4049309187907299e-11) (3.2979939347191366e-2,1.2454368991622481e-9) (1.6503995713160725e-2,4.6280098552463533e-9) };
\addlegendentry{$k = 3$}
\draw[draw=none] (axis cs:2.4141220176769076e-1,9.7524046654869218e-8) -- (axis cs:1.2971299742290118e-1,2.2879697357593217e-9) node[midway, above, text=teal, fill=white, fill opacity=0.8, text opacity=1, inner sep=1pt, font=\scriptsize] {6.04};
\draw[draw=none] (axis cs:1.2971299742290118e-1,2.2879697357593217e-9) -- (axis cs:6.5736358782959317e-2,7.4049309187907299e-11) node[midway, above, text=teal, fill=white, fill opacity=0.8, text opacity=1, inner sep=1pt, font=\scriptsize] {5.05};
\draw[draw=none] (axis cs:6.5736358782959317e-2,7.4049309187907299e-11) -- (axis cs:3.2979939347191366e-2,1.2454368991622481e-9) node[midway, above, text=teal, fill=white, fill opacity=0.8, text opacity=1, inner sep=1pt, font=\scriptsize] {-4.09};
\draw[draw=none] (axis cs:3.2979939347191366e-2,1.2454368991622481e-9) -- (axis cs:1.6503995713160725e-2,4.6280098552463533e-9) node[midway, above, text=teal, fill=white, fill opacity=0.8, text opacity=1, inner sep=1pt, font=\scriptsize] {-1.90};
\end{loglogaxis}
\end{tikzpicture}
    \end{minipage}
  }
  \caption{Energy error (left) and reconstructed potential error in $L^2(\Omega)$ (right) for the manufactured solution \eqref{eq:manufactured-kirchhoff-love} using the potential-based stabilization \eqref{eq:potential-stabilization}.
    The numbers along the curves are the observed orders between successive refinements.}
  \label{fig:convergence-plots:potential-stabilization}
\end{figure}
\endgroup

\begingroup
\tikzset{every node/.append style={scale=0.85}}
\begin{figure}[!htbp]
  \centering
  \subcaptionbox{Triangular mesh family}{
    \begin{minipage}{0.48\textwidth}
\noindent\begin{tikzpicture}
\begin{loglogaxis}[
  cycle list shift=1,
  ymin=1e-8,
  ymax=1e1,
  enlarge y limits=false,
  width=\linewidth,
  height=0.78\linewidth,
  title={energy norm vs. $h$},
  title style={font=\small},
  label style={font=\scriptsize},
  tick label style={font=\scriptsize},
  grid=both,
  major grid style={gray!40},
  minor grid style={gray!20},
  legend pos=south east,
  legend style={font=\scriptsize, inner xsep=2pt, inner ysep=1pt},
  legend cell align={left}
]
\addplot+[red, mark=square*, thick] coordinates { (5.3033008588991082e-1,1.1366102063750496e0) (2.6516504294495546e-1,2.8875066287120621e-1) (1.3258252147247773e-1,6.1120555071753901e-2) (6.6291260736238949e-2,1.1561236468728605e-2) (3.3145630368119572e-2,2.0591770899256383e-3) };
\addlegendentry{$k = 1$}
\draw[draw=none] (axis cs:5.3033008588991082e-1,1.1366102063750496e0) -- (axis cs:2.6516504294495546e-1,2.8875066287120621e-1) node[midway, above, text=red, fill=white, fill opacity=0.8, text opacity=1, inner sep=1pt, font=\scriptsize] {1.98};
\draw[draw=none] (axis cs:2.6516504294495546e-1,2.8875066287120621e-1) -- (axis cs:1.3258252147247773e-1,6.1120555071753901e-2) node[midway, above, text=red, fill=white, fill opacity=0.8, text opacity=1, inner sep=1pt, font=\scriptsize] {2.24};
\draw[draw=none] (axis cs:1.3258252147247773e-1,6.1120555071753901e-2) -- (axis cs:6.6291260736238949e-2,1.1561236468728605e-2) node[midway, above, text=red, fill=white, fill opacity=0.8, text opacity=1, inner sep=1pt, font=\scriptsize] {2.40};
\draw[draw=none] (axis cs:6.6291260736238949e-2,1.1561236468728605e-2) -- (axis cs:3.3145630368119572e-2,2.0591770899256383e-3) node[midway, above, text=red, fill=white, fill opacity=0.8, text opacity=1, inner sep=1pt, font=\scriptsize] {2.49};
\addplot+[brown, mark=o, thick] coordinates { (5.3033008588991082e-1,4.4037687060438357e-1) (2.6516504294495546e-1,5.3887748627849993e-2) (1.3258252147247773e-1,6.5785624123912466e-3) (6.6291260736238949e-2,8.1208948262706049e-4) (3.3145630368119572e-2,1.0087916334986233e-4) };
\addlegendentry{$k = 2$}
\draw[draw=none] (axis cs:5.3033008588991082e-1,4.4037687060438357e-1) -- (axis cs:2.6516504294495546e-1,5.3887748627849993e-2) node[midway, above, text=brown, fill=white, fill opacity=0.8, text opacity=1, inner sep=1pt, font=\scriptsize] {3.03};
\draw[draw=none] (axis cs:2.6516504294495546e-1,5.3887748627849993e-2) -- (axis cs:1.3258252147247773e-1,6.5785624123912466e-3) node[midway, above, text=brown, fill=white, fill opacity=0.8, text opacity=1, inner sep=1pt, font=\scriptsize] {3.03};
\draw[draw=none] (axis cs:1.3258252147247773e-1,6.5785624123912466e-3) -- (axis cs:6.6291260736238949e-2,8.1208948262706049e-4) node[midway, above, text=brown, fill=white, fill opacity=0.8, text opacity=1, inner sep=1pt, font=\scriptsize] {3.02};
\draw[draw=none] (axis cs:6.6291260736238949e-2,8.1208948262706049e-4) -- (axis cs:3.3145630368119572e-2,1.0087916334986233e-4) node[midway, above, text=brown, fill=white, fill opacity=0.8, text opacity=1, inner sep=1pt, font=\scriptsize] {3.01};
\addplot+[teal, mark=triangle*, thick] coordinates { (5.3033008588991082e-1,3.8779703400024575e-2) (2.6516504294495546e-1,2.4358066473014416e-3) (1.3258252147247773e-1,1.5645897646795196e-4) (6.6291260736238949e-2,9.8801227558885538e-6) (3.3145630368119572e-2,6.2743586168995003e-7) };
\addlegendentry{$k = 3$}
\draw[draw=none] (axis cs:5.3033008588991082e-1,3.8779703400024575e-2) -- (axis cs:2.6516504294495546e-1,2.4358066473014416e-3) node[midway, above, text=teal, fill=white, fill opacity=0.8, text opacity=1, inner sep=1pt, font=\scriptsize] {3.99};
\draw[draw=none] (axis cs:2.6516504294495546e-1,2.4358066473014416e-3) -- (axis cs:1.3258252147247773e-1,1.5645897646795196e-4) node[midway, above, text=teal, fill=white, fill opacity=0.8, text opacity=1, inner sep=1pt, font=\scriptsize] {3.96};
\draw[draw=none] (axis cs:1.3258252147247773e-1,1.5645897646795196e-4) -- (axis cs:6.6291260736238949e-2,9.8801227558885538e-6) node[midway, above, text=teal, fill=white, fill opacity=0.8, text opacity=1, inner sep=1pt, font=\scriptsize] {3.99};
\draw[draw=none] (axis cs:6.6291260736238949e-2,9.8801227558885538e-6) -- (axis cs:3.3145630368119572e-2,6.2743586168995003e-7) node[midway, above, text=teal, fill=white, fill opacity=0.8, text opacity=1, inner sep=1pt, font=\scriptsize] {3.98};
\end{loglogaxis}
\end{tikzpicture}
    \end{minipage}
    \begin{minipage}{0.48\textwidth}
\noindent\begin{tikzpicture}
\begin{loglogaxis}[
  cycle list shift=1,
  ymin=1e-11,
  ymax=1e-1,
  enlarge y limits=false,
  width=\linewidth,
  height=0.78\linewidth,
  title={l2 potential norm vs. $h$},
  title style={font=\small},
  label style={font=\scriptsize},
  tick label style={font=\scriptsize},
  grid=both,
  major grid style={gray!40},
  minor grid style={gray!20},
  legend pos=south east,
  legend style={font=\scriptsize, inner xsep=2pt, inner ysep=1pt},
  legend cell align={left}
]
\addplot+[red, mark=square*, thick] coordinates { (5.3033008588991082e-1,9.9778127047748955e-3) (2.6516504294495546e-1,7.9987600297565795e-4) (1.3258252147247773e-1,7.4690118166765088e-5) (6.6291260736238949e-2,5.9069510305123007e-6) (3.3145630368119572e-2,4.0454470747568716e-7) };
\addlegendentry{$k = 1$}
\draw[draw=none] (axis cs:5.3033008588991082e-1,9.9778127047748955e-3) -- (axis cs:2.6516504294495546e-1,7.9987600297565795e-4) node[midway, above, text=red, fill=white, fill opacity=0.8, text opacity=1, inner sep=1pt, font=\scriptsize] {3.64};
\draw[draw=none] (axis cs:2.6516504294495546e-1,7.9987600297565795e-4) -- (axis cs:1.3258252147247773e-1,7.4690118166765088e-5) node[midway, above, text=red, fill=white, fill opacity=0.8, text opacity=1, inner sep=1pt, font=\scriptsize] {3.42};
\draw[draw=none] (axis cs:1.3258252147247773e-1,7.4690118166765088e-5) -- (axis cs:6.6291260736238949e-2,5.9069510305123007e-6) node[midway, above, text=red, fill=white, fill opacity=0.8, text opacity=1, inner sep=1pt, font=\scriptsize] {3.66};
\draw[draw=none] (axis cs:6.6291260736238949e-2,5.9069510305123007e-6) -- (axis cs:3.3145630368119572e-2,4.0454470747568716e-7) node[midway, above, text=red, fill=white, fill opacity=0.8, text opacity=1, inner sep=1pt, font=\scriptsize] {3.87};
\addplot+[brown, mark=o, thick] coordinates { (5.3033008588991082e-1,9.0416994016830369e-4) (2.6516504294495546e-1,2.6983539642897052e-5) (1.3258252147247773e-1,6.2391256215866707e-7) (6.6291260736238949e-2,1.4627731356693631e-8) (3.3145630368119572e-2,7.2264285247531604e-10) };
\addlegendentry{$k = 2$}
\draw[draw=none] (axis cs:5.3033008588991082e-1,9.0416994016830369e-4) -- (axis cs:2.6516504294495546e-1,2.6983539642897052e-5) node[midway, above, text=brown, fill=white, fill opacity=0.8, text opacity=1, inner sep=1pt, font=\scriptsize] {5.07};
\draw[draw=none] (axis cs:2.6516504294495546e-1,2.6983539642897052e-5) -- (axis cs:1.3258252147247773e-1,6.2391256215866707e-7) node[midway, above, text=brown, fill=white, fill opacity=0.8, text opacity=1, inner sep=1pt, font=\scriptsize] {5.43};
\draw[draw=none] (axis cs:1.3258252147247773e-1,6.2391256215866707e-7) -- (axis cs:6.6291260736238949e-2,1.4627731356693631e-8) node[midway, above, text=brown, fill=white, fill opacity=0.8, text opacity=1, inner sep=1pt, font=\scriptsize] {5.41};
\draw[draw=none] (axis cs:6.6291260736238949e-2,1.4627731356693631e-8) -- (axis cs:3.3145630368119572e-2,7.2264285247531604e-10) node[midway, above, text=brown, fill=white, fill opacity=0.8, text opacity=1, inner sep=1pt, font=\scriptsize] {4.34};
\addplot+[teal, mark=triangle*, thick] coordinates { (5.3033008588991082e-1,3.7531383884137935e-5) (2.6516504294495546e-1,5.7597467285680484e-7) (1.3258252147247773e-1,7.1491477451534321e-9) (6.6291260736238949e-2,2.4681923684533964e-10) (3.3145630368119572e-2,1.3712541555034167e-9) };
\addlegendentry{$k = 3$}
\draw[draw=none] (axis cs:5.3033008588991082e-1,3.7531383884137935e-5) -- (axis cs:2.6516504294495546e-1,5.7597467285680484e-7) node[midway, above, text=teal, fill=white, fill opacity=0.8, text opacity=1, inner sep=1pt, font=\scriptsize] {6.03};
\draw[draw=none] (axis cs:2.6516504294495546e-1,5.7597467285680484e-7) -- (axis cs:1.3258252147247773e-1,7.1491477451534321e-9) node[midway, above, text=teal, fill=white, fill opacity=0.8, text opacity=1, inner sep=1pt, font=\scriptsize] {6.33};
\draw[draw=none] (axis cs:1.3258252147247773e-1,7.1491477451534321e-9) -- (axis cs:6.6291260736238949e-2,2.4681923684533964e-10) node[midway, above, text=teal, fill=white, fill opacity=0.8, text opacity=1, inner sep=1pt, font=\scriptsize] {4.86};
\draw[draw=none] (axis cs:6.6291260736238949e-2,2.4681923684533964e-10) -- (axis cs:3.3145630368119572e-2,1.3712541555034167e-9) node[midway, above, text=teal, fill=white, fill opacity=0.8, text opacity=1, inner sep=1pt, font=\scriptsize] {-2.47};
\end{loglogaxis}
\end{tikzpicture}
    \end{minipage}
  }
  \\
  \subcaptionbox{Cartesian mesh family}{
    \begin{minipage}{0.48\textwidth}
\noindent\begin{tikzpicture}
\begin{loglogaxis}[
  cycle list shift=1,
  ymin=1e-8,
  ymax=1e1,
  enlarge y limits=false,
  width=\linewidth,
  height=0.78\linewidth,
  title={energy norm vs. $h$},
  title style={font=\small},
  label style={font=\scriptsize},
  tick label style={font=\scriptsize},
  grid=both,
  major grid style={gray!40},
  minor grid style={gray!20},
  legend pos=south east,
  legend style={font=\scriptsize, inner xsep=2pt, inner ysep=1pt},
  legend cell align={left}
]
\addplot+[red, mark=square*, thick] coordinates { (3.5355339059421481e-1,6.7111571846710849e-1) (1.7677669529722873e-1,1.0604383451064117e-1) (8.8388347648650084e-2,1.4342517619994758e-2) (4.4194173824334396e-2,1.8386532405428198e-3) (2.2097086912169481e-2,2.3181698568148845e-4) };
\addlegendentry{$k = 1$}
\draw[draw=none] (axis cs:3.5355339059421481e-1,6.7111571846710849e-1) -- (axis cs:1.7677669529722873e-1,1.0604383451064117e-1) node[midway, above, text=red, fill=white, fill opacity=0.8, text opacity=1, inner sep=1pt, font=\scriptsize] {2.66};
\draw[draw=none] (axis cs:1.7677669529722873e-1,1.0604383451064117e-1) -- (axis cs:8.8388347648650084e-2,1.4342517619994758e-2) node[midway, above, text=red, fill=white, fill opacity=0.8, text opacity=1, inner sep=1pt, font=\scriptsize] {2.89};
\draw[draw=none] (axis cs:8.8388347648650084e-2,1.4342517619994758e-2) -- (axis cs:4.4194173824334396e-2,1.8386532405428198e-3) node[midway, above, text=red, fill=white, fill opacity=0.8, text opacity=1, inner sep=1pt, font=\scriptsize] {2.96};
\draw[draw=none] (axis cs:4.4194173824334396e-2,1.8386532405428198e-3) -- (axis cs:2.2097086912169481e-2,2.3181698568148845e-4) node[midway, above, text=red, fill=white, fill opacity=0.8, text opacity=1, inner sep=1pt, font=\scriptsize] {2.99};
\addplot+[brown, mark=o, thick] coordinates { (3.5355339059421481e-1,1.4903655393162810e-1) (1.7677669529722873e-1,2.0290458709576412e-2) (8.8388347648650084e-2,2.6629585332263051e-3) (4.4194173824334396e-2,3.4175167780804748e-4) (2.2097086912169481e-2,4.3298481918428494e-5) };
\addlegendentry{$k = 2$}
\draw[draw=none] (axis cs:3.5355339059421481e-1,1.4903655393162810e-1) -- (axis cs:1.7677669529722873e-1,2.0290458709576412e-2) node[midway, above, text=brown, fill=white, fill opacity=0.8, text opacity=1, inner sep=1pt, font=\scriptsize] {2.88};
\draw[draw=none] (axis cs:1.7677669529722873e-1,2.0290458709576412e-2) -- (axis cs:8.8388347648650084e-2,2.6629585332263051e-3) node[midway, above, text=brown, fill=white, fill opacity=0.8, text opacity=1, inner sep=1pt, font=\scriptsize] {2.93};
\draw[draw=none] (axis cs:8.8388347648650084e-2,2.6629585332263051e-3) -- (axis cs:4.4194173824334396e-2,3.4175167780804748e-4) node[midway, above, text=brown, fill=white, fill opacity=0.8, text opacity=1, inner sep=1pt, font=\scriptsize] {2.96};
\draw[draw=none] (axis cs:4.4194173824334396e-2,3.4175167780804748e-4) -- (axis cs:2.2097086912169481e-2,4.3298481918428494e-5) node[midway, above, text=brown, fill=white, fill opacity=0.8, text opacity=1, inner sep=1pt, font=\scriptsize] {2.98};
\addplot+[teal, mark=triangle*, thick] coordinates { (3.5355339059421481e-1,2.0757354836175503e-2) (1.7677669529722873e-1,1.3898842610265525e-3) (8.8388347648650084e-2,8.9093508112645745e-5) (4.4194173824334396e-2,5.6305868667554623e-6) (2.2097086912169481e-2,3.5963769680466136e-7) };
\addlegendentry{$k = 3$}
\draw[draw=none] (axis cs:3.5355339059421481e-1,2.0757354836175503e-2) -- (axis cs:1.7677669529722873e-1,1.3898842610265525e-3) node[midway, above, text=teal, fill=white, fill opacity=0.8, text opacity=1, inner sep=1pt, font=\scriptsize] {3.90};
\draw[draw=none] (axis cs:1.7677669529722873e-1,1.3898842610265525e-3) -- (axis cs:8.8388347648650084e-2,8.9093508112645745e-5) node[midway, above, text=teal, fill=white, fill opacity=0.8, text opacity=1, inner sep=1pt, font=\scriptsize] {3.96};
\draw[draw=none] (axis cs:8.8388347648650084e-2,8.9093508112645745e-5) -- (axis cs:4.4194173824334396e-2,5.6305868667554623e-6) node[midway, above, text=teal, fill=white, fill opacity=0.8, text opacity=1, inner sep=1pt, font=\scriptsize] {3.98};
\draw[draw=none] (axis cs:4.4194173824334396e-2,5.6305868667554623e-6) -- (axis cs:2.2097086912169481e-2,3.5963769680466136e-7) node[midway, above, text=teal, fill=white, fill opacity=0.8, text opacity=1, inner sep=1pt, font=\scriptsize] {3.97};
\end{loglogaxis}
\end{tikzpicture}
    \end{minipage}
    \begin{minipage}{0.48\textwidth}
\noindent\begin{tikzpicture}
\begin{loglogaxis}[
  cycle list shift=1,
  ymin=1e-11,
  ymax=1e0,
  enlarge y limits=false,
  width=\linewidth,
  height=0.78\linewidth,
  title={l2 potential norm vs. $h$},
  title style={font=\small},
  label style={font=\scriptsize},
  tick label style={font=\scriptsize},
  grid=both,
  major grid style={gray!40},
  minor grid style={gray!20},
  legend pos=south east,
  legend style={font=\scriptsize, inner xsep=2pt, inner ysep=1pt},
  legend cell align={left}
]
\addplot+[red, mark=square*, thick] coordinates { (3.5355339059421481e-1,7.3881428078425380e-3) (1.7677669529722873e-1,6.7109898883444810e-4) (8.8388347648650084e-2,4.7380832319068326e-5) (4.4194173824334396e-2,3.0633382432351684e-6) (2.2097086912169481e-2,1.9314893909841945e-7) };
\addlegendentry{$k = 1$}
\draw[draw=none] (axis cs:3.5355339059421481e-1,7.3881428078425380e-3) -- (axis cs:1.7677669529722873e-1,6.7109898883444810e-4) node[midway, above, text=red, fill=white, fill opacity=0.8, text opacity=1, inner sep=1pt, font=\scriptsize] {3.46};
\draw[draw=none] (axis cs:1.7677669529722873e-1,6.7109898883444810e-4) -- (axis cs:8.8388347648650084e-2,4.7380832319068326e-5) node[midway, above, text=red, fill=white, fill opacity=0.8, text opacity=1, inner sep=1pt, font=\scriptsize] {3.82};
\draw[draw=none] (axis cs:8.8388347648650084e-2,4.7380832319068326e-5) -- (axis cs:4.4194173824334396e-2,3.0633382432351684e-6) node[midway, above, text=red, fill=white, fill opacity=0.8, text opacity=1, inner sep=1pt, font=\scriptsize] {3.95};
\draw[draw=none] (axis cs:4.4194173824334396e-2,3.0633382432351684e-6) -- (axis cs:2.2097086912169481e-2,1.9314893909841945e-7) node[midway, above, text=red, fill=white, fill opacity=0.8, text opacity=1, inner sep=1pt, font=\scriptsize] {3.99};
\addplot+[brown, mark=o, thick] coordinates { (3.5355339059421481e-1,2.7080687961442476e-4) (1.7677669529722873e-1,5.6159361796988817e-6) (8.8388347648650084e-2,1.0720282321670221e-7) (4.4194173824334396e-2,2.5474858869546837e-9) (2.2097086912169481e-2,7.3079248045894059e-10) };
\addlegendentry{$k = 2$}
\draw[draw=none] (axis cs:3.5355339059421481e-1,2.7080687961442476e-4) -- (axis cs:1.7677669529722873e-1,5.6159361796988817e-6) node[midway, above, text=brown, fill=white, fill opacity=0.8, text opacity=1, inner sep=1pt, font=\scriptsize] {5.59};
\draw[draw=none] (axis cs:1.7677669529722873e-1,5.6159361796988817e-6) -- (axis cs:8.8388347648650084e-2,1.0720282321670221e-7) node[midway, above, text=brown, fill=white, fill opacity=0.8, text opacity=1, inner sep=1pt, font=\scriptsize] {5.71};
\draw[draw=none] (axis cs:8.8388347648650084e-2,1.0720282321670221e-7) -- (axis cs:4.4194173824334396e-2,2.5474858869546837e-9) node[midway, above, text=brown, fill=white, fill opacity=0.8, text opacity=1, inner sep=1pt, font=\scriptsize] {5.40};
\draw[draw=none] (axis cs:4.4194173824334396e-2,2.5474858869546837e-9) -- (axis cs:2.2097086912169481e-2,7.3079248045894059e-10) node[midway, above, text=brown, fill=white, fill opacity=0.8, text opacity=1, inner sep=1pt, font=\scriptsize] {1.80};
\addplot+[teal, mark=triangle*, thick] coordinates { (3.5355339059421481e-1,1.5502600452426904e-5) (1.7677669529722873e-1,1.5918200831379220e-7) (8.8388347648650084e-2,1.9864383807743298e-9) (4.4194173824334396e-2,6.5983077237132235e-11) (2.2097086912169481e-2,1.5137027218903972e-9) };
\addlegendentry{$k = 3$}
\draw[draw=none] (axis cs:3.5355339059421481e-1,1.5502600452426904e-5) -- (axis cs:1.7677669529722873e-1,1.5918200831379220e-7) node[midway, above, text=teal, fill=white, fill opacity=0.8, text opacity=1, inner sep=1pt, font=\scriptsize] {6.61};
\draw[draw=none] (axis cs:1.7677669529722873e-1,1.5918200831379220e-7) -- (axis cs:8.8388347648650084e-2,1.9864383807743298e-9) node[midway, above, text=teal, fill=white, fill opacity=0.8, text opacity=1, inner sep=1pt, font=\scriptsize] {6.32};
\draw[draw=none] (axis cs:8.8388347648650084e-2,1.9864383807743298e-9) -- (axis cs:4.4194173824334396e-2,6.5983077237132235e-11) node[midway, above, text=teal, fill=white, fill opacity=0.8, text opacity=1, inner sep=1pt, font=\scriptsize] {4.91};
\draw[draw=none] (axis cs:4.4194173824334396e-2,6.5983077237132235e-11) -- (axis cs:2.2097086912169481e-2,1.5137027218903972e-9) node[midway, above, text=teal, fill=white, fill opacity=0.8, text opacity=1, inner sep=1pt, font=\scriptsize] {-4.52};
\end{loglogaxis}
\end{tikzpicture}
    \end{minipage}
  }
  \\
  \subcaptionbox{Hexagonal mesh family}{
    \begin{minipage}{0.48\textwidth}
\noindent\begin{tikzpicture}
\begin{loglogaxis}[
  cycle list shift=1,
  ymin=1e-7,
  ymax=1e1,
  enlarge y limits=false,
  width=\linewidth,
  height=0.78\linewidth,
  title={energy norm vs. $h$},
  title style={font=\small},
  label style={font=\scriptsize},
  tick label style={font=\scriptsize},
  grid=both,
  major grid style={gray!40},
  minor grid style={gray!20},
  legend pos=south east,
  legend style={font=\scriptsize, inner xsep=2pt, inner ysep=1pt},
  legend cell align={left}
]
\addplot+[red, mark=square*, thick] coordinates { (2.4141220176769076e-1,2.3171658482280380e-1) (1.2971299742290118e-1,6.7810706707061810e-2) (6.5736358782959317e-2,1.8043902251068403e-2) (3.2979939347191366e-2,4.6283124323850070e-3) (1.6503995713160725e-2,1.1696859971684821e-3) };
\addlegendentry{$k = 1$}
\draw[draw=none] (axis cs:2.4141220176769076e-1,2.3171658482280380e-1) -- (axis cs:1.2971299742290118e-1,6.7810706707061810e-2) node[midway, above, text=red, fill=white, fill opacity=0.8, text opacity=1, inner sep=1pt, font=\scriptsize] {1.98};
\draw[draw=none] (axis cs:1.2971299742290118e-1,6.7810706707061810e-2) -- (axis cs:6.5736358782959317e-2,1.8043902251068403e-2) node[midway, above, text=red, fill=white, fill opacity=0.8, text opacity=1, inner sep=1pt, font=\scriptsize] {1.95};
\draw[draw=none] (axis cs:6.5736358782959317e-2,1.8043902251068403e-2) -- (axis cs:3.2979939347191366e-2,4.6283124323850070e-3) node[midway, above, text=red, fill=white, fill opacity=0.8, text opacity=1, inner sep=1pt, font=\scriptsize] {1.97};
\draw[draw=none] (axis cs:3.2979939347191366e-2,4.6283124323850070e-3) -- (axis cs:1.6503995713160725e-2,1.1696859971684821e-3) node[midway, above, text=red, fill=white, fill opacity=0.8, text opacity=1, inner sep=1pt, font=\scriptsize] {1.99};
\addplot+[brown, mark=o, thick] coordinates { (2.4141220176769076e-1,2.4764808096042333e-2) (1.2971299742290118e-1,4.0198245377986719e-3) (6.5736358782959317e-2,5.4417199229171308e-4) (3.2979939347191366e-2,6.9895056157929267e-5) (1.6503995713160725e-2,8.8786444384559414e-6) };
\addlegendentry{$k = 2$}
\draw[draw=none] (axis cs:2.4141220176769076e-1,2.4764808096042333e-2) -- (axis cs:1.2971299742290118e-1,4.0198245377986719e-3) node[midway, above, text=brown, fill=white, fill opacity=0.8, text opacity=1, inner sep=1pt, font=\scriptsize] {2.93};
\draw[draw=none] (axis cs:1.2971299742290118e-1,4.0198245377986719e-3) -- (axis cs:6.5736358782959317e-2,5.4417199229171308e-4) node[midway, above, text=brown, fill=white, fill opacity=0.8, text opacity=1, inner sep=1pt, font=\scriptsize] {2.94};
\draw[draw=none] (axis cs:6.5736358782959317e-2,5.4417199229171308e-4) -- (axis cs:3.2979939347191366e-2,6.9895056157929267e-5) node[midway, above, text=brown, fill=white, fill opacity=0.8, text opacity=1, inner sep=1pt, font=\scriptsize] {2.98};
\draw[draw=none] (axis cs:3.2979939347191366e-2,6.9895056157929267e-5) -- (axis cs:1.6503995713160725e-2,8.8786444384559414e-6) node[midway, above, text=brown, fill=white, fill opacity=0.8, text opacity=1, inner sep=1pt, font=\scriptsize] {2.98};
\addplot+[teal, mark=triangle*, thick] coordinates { (2.4141220176769076e-1,1.7163574322476862e-3) (1.2971299742290118e-1,1.4436838039981088e-4) (6.5736358782959317e-2,9.7794792152773599e-6) (3.2979939347191366e-2,7.0684446245630955e-7) (1.6503995713160725e-2,2.6027793254736623e-6) };
\addlegendentry{$k = 3$}
\draw[draw=none] (axis cs:2.4141220176769076e-1,1.7163574322476862e-3) -- (axis cs:1.2971299742290118e-1,1.4436838039981088e-4) node[midway, above, text=teal, fill=white, fill opacity=0.8, text opacity=1, inner sep=1pt, font=\scriptsize] {3.99};
\draw[draw=none] (axis cs:1.2971299742290118e-1,1.4436838039981088e-4) -- (axis cs:6.5736358782959317e-2,9.7794792152773599e-6) node[midway, above, text=teal, fill=white, fill opacity=0.8, text opacity=1, inner sep=1pt, font=\scriptsize] {3.96};
\draw[draw=none] (axis cs:6.5736358782959317e-2,9.7794792152773599e-6) -- (axis cs:3.2979939347191366e-2,7.0684446245630955e-7) node[midway, above, text=teal, fill=white, fill opacity=0.8, text opacity=1, inner sep=1pt, font=\scriptsize] {3.81};
\draw[draw=none] (axis cs:3.2979939347191366e-2,7.0684446245630955e-7) -- (axis cs:1.6503995713160725e-2,2.6027793254736623e-6) node[midway, above, text=teal, fill=white, fill opacity=0.8, text opacity=1, inner sep=1pt, font=\scriptsize] {-1.88};
\end{loglogaxis}
\end{tikzpicture}
    \end{minipage}
    \begin{minipage}{0.48\textwidth}
\noindent\begin{tikzpicture}
\begin{loglogaxis}[
  cycle list shift=1,
  ymin=1e-11,
  ymax=1e-1,
  enlarge y limits=false,
  width=\linewidth,
  height=0.78\linewidth,
  title={l2 potential norm vs. $h$},
  title style={font=\small},
  label style={font=\scriptsize},
  tick label style={font=\scriptsize},
  grid=both,
  major grid style={gray!40},
  minor grid style={gray!20},
  legend pos=south east,
  legend style={font=\scriptsize, inner xsep=2pt, inner ysep=1pt},
  legend cell align={left}
]
\addplot+[red, mark=square*, thick] coordinates { (2.4141220176769076e-1,6.6608029349912528e-4) (1.2971299742290118e-1,7.1037386503621769e-5) (6.5736358782959317e-2,5.2165005136104410e-6) (3.2979939347191366e-2,3.4352269653414489e-7) (1.6503995713160725e-2,2.1982189065695533e-8) };
\addlegendentry{$k = 1$}
\draw[draw=none] (axis cs:2.4141220176769076e-1,6.6608029349912528e-4) -- (axis cs:1.2971299742290118e-1,7.1037386503621769e-5) node[midway, above, text=red, fill=white, fill opacity=0.8, text opacity=1, inner sep=1pt, font=\scriptsize] {3.60};
\draw[draw=none] (axis cs:1.2971299742290118e-1,7.1037386503621769e-5) -- (axis cs:6.5736358782959317e-2,5.2165005136104410e-6) node[midway, above, text=red, fill=white, fill opacity=0.8, text opacity=1, inner sep=1pt, font=\scriptsize] {3.84};
\draw[draw=none] (axis cs:6.5736358782959317e-2,5.2165005136104410e-6) -- (axis cs:3.2979939347191366e-2,3.4352269653414489e-7) node[midway, above, text=red, fill=white, fill opacity=0.8, text opacity=1, inner sep=1pt, font=\scriptsize] {3.94};
\draw[draw=none] (axis cs:3.2979939347191366e-2,3.4352269653414489e-7) -- (axis cs:1.6503995713160725e-2,2.1982189065695533e-8) node[midway, above, text=red, fill=white, fill opacity=0.8, text opacity=1, inner sep=1pt, font=\scriptsize] {3.97};
\addplot+[brown, mark=o, thick] coordinates { (2.4141220176769076e-1,9.7391503086200963e-6) (1.2971299742290118e-1,3.2811142374312410e-7) (6.5736358782959317e-2,8.3736033382070245e-9) (3.2979939347191366e-2,5.9826238869921495e-9) (1.6503995713160725e-2,8.5957220732557805e-9) };
\addlegendentry{$k = 2$}
\draw[draw=none] (axis cs:2.4141220176769076e-1,9.7391503086200963e-6) -- (axis cs:1.2971299742290118e-1,3.2811142374312410e-7) node[midway, above, text=brown, fill=white, fill opacity=0.8, text opacity=1, inner sep=1pt, font=\scriptsize] {5.46};
\draw[draw=none] (axis cs:1.2971299742290118e-1,3.2811142374312410e-7) -- (axis cs:6.5736358782959317e-2,8.3736033382070245e-9) node[midway, above, text=brown, fill=white, fill opacity=0.8, text opacity=1, inner sep=1pt, font=\scriptsize] {5.40};
\draw[draw=none] (axis cs:6.5736358782959317e-2,8.3736033382070245e-9) -- (axis cs:3.2979939347191366e-2,5.9826238869921495e-9) node[midway, above, text=brown, fill=white, fill opacity=0.8, text opacity=1, inner sep=1pt, font=\scriptsize] {0.49};
\draw[draw=none] (axis cs:3.2979939347191366e-2,5.9826238869921495e-9) -- (axis cs:1.6503995713160725e-2,8.5957220732557805e-9) node[midway, above, text=brown, fill=white, fill opacity=0.8, text opacity=1, inner sep=1pt, font=\scriptsize] {-0.52};
\addplot+[teal, mark=triangle*, thick] coordinates { (2.4141220176769076e-1,2.8236941333878056e-7) (1.2971299742290118e-1,5.4253117329541393e-9) (6.5736358782959317e-2,9.5212230275598809e-10) (3.2979939347191366e-2,3.7343531752342603e-9) (1.6503995713160725e-2,3.1277861298671562e-8) };
\addlegendentry{$k = 3$}
\draw[draw=none] (axis cs:2.4141220176769076e-1,2.8236941333878056e-7) -- (axis cs:1.2971299742290118e-1,5.4253117329541393e-9) node[midway, above, text=teal, fill=white, fill opacity=0.8, text opacity=1, inner sep=1pt, font=\scriptsize] {6.36};
\draw[draw=none] (axis cs:1.2971299742290118e-1,5.4253117329541393e-9) -- (axis cs:6.5736358782959317e-2,9.5212230275598809e-10) node[midway, above, text=teal, fill=white, fill opacity=0.8, text opacity=1, inner sep=1pt, font=\scriptsize] {2.56};
\draw[draw=none] (axis cs:6.5736358782959317e-2,9.5212230275598809e-10) -- (axis cs:3.2979939347191366e-2,3.7343531752342603e-9) node[midway, above, text=teal, fill=white, fill opacity=0.8, text opacity=1, inner sep=1pt, font=\scriptsize] {-1.98};
\draw[draw=none] (axis cs:3.2979939347191366e-2,3.7343531752342603e-9) -- (axis cs:1.6503995713160725e-2,3.1277861298671562e-8) node[midway, above, text=teal, fill=white, fill opacity=0.8, text opacity=1, inner sep=1pt, font=\scriptsize] {-3.07};
\end{loglogaxis}
\end{tikzpicture}
    \end{minipage}
  }
  \caption{Energy error (left) and reconstructed potential error in $L^2(\Omega)$ (right) for the manufactured solution \eqref{eq:manufactured-kirchhoff-love} using the Hessian-based stabilization \eqref{eq:hessian-stabilization}.
    The numbers along the curves are the observed orders between successive refinements.}
  \label{fig:convergence-plots:hessian-stabilization}
\end{figure}
\endgroup

For each stabilization, we measure the discrete error $\underline e_h\coloneqq\underline u_h-\IHtwo{k}{h}u$ in the energy seminorm induced by the corresponding bilinear form $a_h$.
We also measure the reconstructed $L^2$ potential error $\norm{L^2(\Omega)}{p_h^{k+2}\underline e_h}$.
\Cref{fig:convergence-plots:potential-stabilization,fig:convergence-plots:hessian-stabilization} show these errors for the potential- and Hessian-based stabilizations \eqref{eq:potential-stabilization} and \eqref{eq:hessian-stabilization}, respectively.
The energy-error magnitudes are measured in different but comparable norms (due to the different choices of stabilization), whereas the reconstructed potential errors provide a common measure for comparing the two methods.

The potential-based stabilization \eqref{eq:potential-stabilization} includes the case $k=0$, for which the energy errors converge at first order on all mesh families considered.
This agrees with the lowest-order estimate of \cref{thm:energy-error-estimate}, whose additional regularity assumption is satisfied by the smooth manufactured solution.
For $k\ge1$, both stabilizations exhibit energy convergence consistent with $h^{k+1}$ before numerical saturation.
For $k=2$ and $k=3$, the finer-mesh rates approach $3$ and $4$, respectively, on all three mesh families before saturation.
For $k=1$, both stabilizations exhibit an additional order on Cartesian meshes: the rates on the finest refinements are $2.98$ and $2.99$ for the potential- and Hessian-based stabilizations, respectively.
The corresponding triangular rates are $2.61$ and $2.49$, whereas both hexagonal rates remain close to $2$. 
These enhanced rates at $k=1$ are observed for the present solution and mesh families; their full theoretical justification is left for future work.

The reconstructed potential errors converge faster than the energy errors.
At $k=0$, the finest rates are close to $2$ on all three families.
For $k\ge1$, the rates before saturation are broadly consistent with $h^{k+3}$, approaching $4$ at $k=1$ and reaching values near $5$ and $6$ at $k=2$ and $k=3$, respectively.

On the finest refinements, several higher-degree reconstructed potential-error curves flatten or increase, particularly on the hexagonal family. At $k=3$, the potential-based errors reach approximately $5.8\times10^{-11}$ on triangular meshes and $4.3\times10^{-11}$ on Cartesian meshes before increasing to $2.1\times10^{-9}$ and $1.3\times10^{-10}$, respectively. On the hexagonal family, the potential-based energy error at $k=3$ also nearly stops decreasing on the last refinement.
These departures are consistent with finite-precision effects becoming comparable to the discretization error.
The potential-based stabilization generally sustains convergence down to smaller error levels before saturation.
This is likely a consequence of better conditioning of the linear system associated with the potential-based stabilization.
\Cref{fig:condition-numbers} compares the spectral-norm condition numbers of the statically condensed systems on triangular and Cartesian meshes, including $k=0$ only for the potential-based stabilization.
The condition numbers appear to grow like $h^{-4}$, as expected. For each degree $k=1,2,3$ shared by the two methods, the potential-based stabilization yields smaller condition numbers on every tested mesh, with a more pronounced difference at higher polynomial degrees.

\begingroup
\tikzset{every node/.append style={scale=0.85}}
\pgfplotsset{every axis post/.append style={xticklabel style={rotate=25,anchor=north east}}}
\begin{figure}[!htbp]
  \centering
  \subcaptionbox{Triangular mesh family}{
    \begin{minipage}{0.48\textwidth}
\noindent\begin{tikzpicture}
\begin{loglogaxis}[
  ymin=1e1,
  ymax=1e10,
  enlarge y limits=false,
  width=\linewidth, height=0.78\linewidth,
  xtick=data, xticklabel={\begingroup\pgfkeys{/pgf/fpu=true}\pgfmathparse{exp(\tick)}\pgfmathprintnumber[sci,precision=2]{\pgfmathresult}\endgroup},
  title style={font=\small}, label style={font=\scriptsize},
  tick label style={font=\scriptsize}, grid=both,
  major grid style={gray!40}, minor grid style={gray!20},
  legend pos=north east, legend style={font=\scriptsize}, legend cell align={left}]
\addplot+[blue, mark=*, thick] coordinates { (5.3033008588991082e-01,9.4090276324246133e+01) (2.6516504294495546e-01,1.3158124223207767e+03) (1.3258252147247773e-01,2.0317515096252322e+04) (6.6291260736238949e-02,3.2582838008703478e+05) (3.3145630368119572e-02,5.2455397538726470e+06) };
\addlegendentry{$k = 0$}
\draw[draw=none] (axis cs:5.3033008588991082e-01,9.4090276324246133e+01) -- (axis cs:2.6516504294495546e-01,1.3158124223207767e+03) node[midway, above, text=blue, fill=white, fill opacity=0.8, text opacity=1, inner sep=1pt, font=\scriptsize] {3.81};
\draw[draw=none] (axis cs:2.6516504294495546e-01,1.3158124223207767e+03) -- (axis cs:1.3258252147247773e-01,2.0317515096252322e+04) node[midway, above, text=blue, fill=white, fill opacity=0.8, text opacity=1, inner sep=1pt, font=\scriptsize] {3.95};
\draw[draw=none] (axis cs:1.3258252147247773e-01,2.0317515096252322e+04) -- (axis cs:6.6291260736238949e-02,3.2582838008703478e+05) node[midway, above, text=blue, fill=white, fill opacity=0.8, text opacity=1, inner sep=1pt, font=\scriptsize] {4.00};
\draw[draw=none] (axis cs:6.6291260736238949e-02,3.2582838008703478e+05) -- (axis cs:3.3145630368119572e-02,5.2455397538726470e+06) node[midway, above, text=blue, fill=white, fill opacity=0.8, text opacity=1, inner sep=1pt, font=\scriptsize] {4.01};
\addplot+[red, mark=square*, thick] coordinates { (5.3033008588991082e-01,4.0073883560354273e+02) (2.6516504294495546e-01,4.8780195389052997e+03) (1.3258252147247773e-01,7.7990710890873219e+04) (6.6291260736238949e-02,1.2460892225498431e+06) (3.3145630368119572e-02,1.9653206766864445e+07) };
\addlegendentry{$k = 1$}
\draw[draw=none] (axis cs:5.3033008588991082e-01,4.0073883560354273e+02) -- (axis cs:2.6516504294495546e-01,4.8780195389052997e+03) node[midway, above, text=red, fill=white, fill opacity=0.8, text opacity=1, inner sep=1pt, font=\scriptsize] {3.61};
\draw[draw=none] (axis cs:2.6516504294495546e-01,4.8780195389052997e+03) -- (axis cs:1.3258252147247773e-01,7.7990710890873219e+04) node[midway, above, text=red, fill=white, fill opacity=0.8, text opacity=1, inner sep=1pt, font=\scriptsize] {4.00};
\draw[draw=none] (axis cs:1.3258252147247773e-01,7.7990710890873219e+04) -- (axis cs:6.6291260736238949e-02,1.2460892225498431e+06) node[midway, above, text=red, fill=white, fill opacity=0.8, text opacity=1, inner sep=1pt, font=\scriptsize] {4.00};
\draw[draw=none] (axis cs:6.6291260736238949e-02,1.2460892225498431e+06) -- (axis cs:3.3145630368119572e-02,1.9653206766864445e+07) node[midway, above, text=red, fill=white, fill opacity=0.8, text opacity=1, inner sep=1pt, font=\scriptsize] {3.98};
\addplot+[brown, mark=o, thick] coordinates { (5.3033008588991082e-01,7.1418043933785375e+03) (2.6516504294495546e-01,6.9028625199411210e+04) (1.3258252147247773e-01,9.0258682446224603e+05) (6.6291260736238949e-02,1.3465286594508654e+07) (3.3145630368119572e-02,2.0927758299045131e+08) };
\addlegendentry{$k = 2$}
\draw[draw=none] (axis cs:5.3033008588991082e-01,7.1418043933785375e+03) -- (axis cs:2.6516504294495546e-01,6.9028625199411210e+04) node[midway, above, text=brown, fill=white, fill opacity=0.8, text opacity=1, inner sep=1pt, font=\scriptsize] {3.27};
\draw[draw=none] (axis cs:2.6516504294495546e-01,6.9028625199411210e+04) -- (axis cs:1.3258252147247773e-01,9.0258682446224603e+05) node[midway, above, text=brown, fill=white, fill opacity=0.8, text opacity=1, inner sep=1pt, font=\scriptsize] {3.71};
\draw[draw=none] (axis cs:1.3258252147247773e-01,9.0258682446224603e+05) -- (axis cs:6.6291260736238949e-02,1.3465286594508654e+07) node[midway, above, text=brown, fill=white, fill opacity=0.8, text opacity=1, inner sep=1pt, font=\scriptsize] {3.90};
\draw[draw=none] (axis cs:6.6291260736238949e-02,1.3465286594508654e+07) -- (axis cs:3.3145630368119572e-02,2.0927758299045131e+08) node[midway, above, text=brown, fill=white, fill opacity=0.8, text opacity=1, inner sep=1pt, font=\scriptsize] {3.96};
\addplot+[teal, mark=triangle*, thick] coordinates { (5.3033008588991082e-01,3.6901978186702392e+04) (2.6516504294495546e-01,3.5109830784454505e+05) (1.3258252147247773e-01,4.4283780197430244e+06) (6.6291260736238949e-02,6.4451254050881594e+07) (3.3145630368119572e-02,9.9191725073382521e+08) };
\addlegendentry{$k = 3$}
\draw[draw=none] (axis cs:5.3033008588991082e-01,3.6901978186702392e+04) -- (axis cs:2.6516504294495546e-01,3.5109830784454505e+05) node[midway, above, text=teal, fill=white, fill opacity=0.8, text opacity=1, inner sep=1pt, font=\scriptsize] {3.25};
\draw[draw=none] (axis cs:2.6516504294495546e-01,3.5109830784454505e+05) -- (axis cs:1.3258252147247773e-01,4.4283780197430244e+06) node[midway, above, text=teal, fill=white, fill opacity=0.8, text opacity=1, inner sep=1pt, font=\scriptsize] {3.66};
\draw[draw=none] (axis cs:1.3258252147247773e-01,4.4283780197430244e+06) -- (axis cs:6.6291260736238949e-02,6.4451254050881594e+07) node[midway, above, text=teal, fill=white, fill opacity=0.8, text opacity=1, inner sep=1pt, font=\scriptsize] {3.86};
\draw[draw=none] (axis cs:6.6291260736238949e-02,6.4451254050881594e+07) -- (axis cs:3.3145630368119572e-02,9.9191725073382521e+08) node[midway, above, text=teal, fill=white, fill opacity=0.8, text opacity=1, inner sep=1pt, font=\scriptsize] {3.94};
\end{loglogaxis}
\end{tikzpicture}
    \end{minipage}
    \begin{minipage}{0.48\textwidth}
\noindent\begin{tikzpicture}
\begin{loglogaxis}[
  cycle list shift=1,
  ymin=1e1,
  ymax=1e10,
  enlarge y limits=false,
  width=\linewidth, height=0.78\linewidth,
  xtick=data, xticklabel={\begingroup\pgfkeys{/pgf/fpu=true}\pgfmathparse{exp(\tick)}\pgfmathprintnumber[sci,precision=2]{\pgfmathresult}\endgroup},
  title style={font=\small}, label style={font=\scriptsize},
  tick label style={font=\scriptsize}, grid=both,
  major grid style={gray!40}, minor grid style={gray!20},
  legend pos=north east, legend style={font=\scriptsize}, legend cell align={left}]
\addplot+[red, mark=square*, thick] coordinates { (5.3033008588991082e-01,7.4716031426418294e+02) (2.6516504294495546e-01,1.0860739544670139e+04) (1.3258252147247773e-01,1.9513782259948092e+05) (6.6291260736238949e-02,3.2182952426945521e+06) (3.3145630368119572e-02,5.1151962060059026e+07) };
\addlegendentry{$k = 1$}
\draw[draw=none] (axis cs:5.3033008588991082e-01,7.4716031426418294e+02) -- (axis cs:2.6516504294495546e-01,1.0860739544670139e+04) node[midway, above, text=red, fill=white, fill opacity=0.8, text opacity=1, inner sep=1pt, font=\scriptsize] {3.86};
\draw[draw=none] (axis cs:2.6516504294495546e-01,1.0860739544670139e+04) -- (axis cs:1.3258252147247773e-01,1.9513782259948092e+05) node[midway, above, text=red, fill=white, fill opacity=0.8, text opacity=1, inner sep=1pt, font=\scriptsize] {4.17};
\draw[draw=none] (axis cs:1.3258252147247773e-01,1.9513782259948092e+05) -- (axis cs:6.6291260736238949e-02,3.2182952426945521e+06) node[midway, above, text=red, fill=white, fill opacity=0.8, text opacity=1, inner sep=1pt, font=\scriptsize] {4.04};
\draw[draw=none] (axis cs:6.6291260736238949e-02,3.2182952426945521e+06) -- (axis cs:3.3145630368119572e-02,5.1151962060059026e+07) node[midway, above, text=red, fill=white, fill opacity=0.8, text opacity=1, inner sep=1pt, font=\scriptsize] {3.99};
\addplot+[brown, mark=o, thick] coordinates { (5.3033008588991082e-01,1.7697412465529349e+04) (2.6516504294495546e-01,1.9470012334354900e+05) (1.3258252147247773e-01,2.7540698166728425e+06) (6.6291260736238949e-02,4.2190850501034625e+07) (3.3145630368119572e-02,6.5944233458370364e+08) };
\addlegendentry{$k = 2$}
\draw[draw=none] (axis cs:5.3033008588991082e-01,1.7697412465529349e+04) -- (axis cs:2.6516504294495546e-01,1.9470012334354900e+05) node[midway, above, text=brown, fill=white, fill opacity=0.8, text opacity=1, inner sep=1pt, font=\scriptsize] {3.46};
\draw[draw=none] (axis cs:2.6516504294495546e-01,1.9470012334354900e+05) -- (axis cs:1.3258252147247773e-01,2.7540698166728425e+06) node[midway, above, text=brown, fill=white, fill opacity=0.8, text opacity=1, inner sep=1pt, font=\scriptsize] {3.82};
\draw[draw=none] (axis cs:1.3258252147247773e-01,2.7540698166728425e+06) -- (axis cs:6.6291260736238949e-02,4.2190850501034625e+07) node[midway, above, text=brown, fill=white, fill opacity=0.8, text opacity=1, inner sep=1pt, font=\scriptsize] {3.94};
\draw[draw=none] (axis cs:6.6291260736238949e-02,4.2190850501034625e+07) -- (axis cs:3.3145630368119572e-02,6.5944233458370364e+08) node[midway, above, text=brown, fill=white, fill opacity=0.8, text opacity=1, inner sep=1pt, font=\scriptsize] {3.97};
\addplot+[teal, mark=triangle*, thick] coordinates { (5.3033008588991082e-01,1.2940886737528276e+05) (2.6516504294495546e-01,1.3787574677030982e+06) (1.3258252147247773e-01,1.7972631141772863e+07) (6.6291260736238949e-02,2.6378094633753133e+08) (3.3145630368119572e-02,4.0628205602055840e+09) };
\addlegendentry{$k = 3$}
\draw[draw=none] (axis cs:5.3033008588991082e-01,1.2940886737528276e+05) -- (axis cs:2.6516504294495546e-01,1.3787574677030982e+06) node[midway, above, text=teal, fill=white, fill opacity=0.8, text opacity=1, inner sep=1pt, font=\scriptsize] {3.41};
\draw[draw=none] (axis cs:2.6516504294495546e-01,1.3787574677030982e+06) -- (axis cs:1.3258252147247773e-01,1.7972631141772863e+07) node[midway, above, text=teal, fill=white, fill opacity=0.8, text opacity=1, inner sep=1pt, font=\scriptsize] {3.70};
\draw[draw=none] (axis cs:1.3258252147247773e-01,1.7972631141772863e+07) -- (axis cs:6.6291260736238949e-02,2.6378094633753133e+08) node[midway, above, text=teal, fill=white, fill opacity=0.8, text opacity=1, inner sep=1pt, font=\scriptsize] {3.88};
\draw[draw=none] (axis cs:6.6291260736238949e-02,2.6378094633753133e+08) -- (axis cs:3.3145630368119572e-02,4.0628205602055840e+09) node[midway, above, text=teal, fill=white, fill opacity=0.8, text opacity=1, inner sep=1pt, font=\scriptsize] {3.95};
\end{loglogaxis}
\end{tikzpicture}
    \end{minipage}
  }
  \medskip\\
  \subcaptionbox{Cartesian mesh family}{
    \begin{minipage}{0.48\textwidth}
\noindent\begin{tikzpicture}
\begin{loglogaxis}[
  ymin=1e1,
  ymax=1e11,
  enlarge y limits=false,
  width=\linewidth, height=0.78\linewidth,
  xtick=data, xticklabel={\begingroup\pgfkeys{/pgf/fpu=true}\pgfmathparse{exp(\tick)}\pgfmathprintnumber[sci,precision=2]{\pgfmathresult}\endgroup},
  title style={font=\small}, label style={font=\scriptsize},
  tick label style={font=\scriptsize}, grid=both,
  major grid style={gray!40}, minor grid style={gray!20},
  legend pos=north east, legend style={font=\scriptsize}, legend cell align={left}]
\addplot+[blue, mark=*, thick] coordinates { (3.5355339059421481e-01,2.2942484142099394e+02) (1.7677669529722873e-01,2.9944770991856217e+03) (8.8388347648650084e-02,4.5237314724283773e+04) (4.4194173824334396e-02,7.1268088530088181e+05) (2.2097086912169481e-02,1.1357391425691146e+07) };
\addlegendentry{$k = 0$}
\draw[draw=none] (axis cs:3.5355339059421481e-01,2.2942484142099394e+02) -- (axis cs:1.7677669529722873e-01,2.9944770991856217e+03) node[midway, above, text=blue, fill=white, fill opacity=0.8, text opacity=1, inner sep=1pt, font=\scriptsize] {3.71};
\draw[draw=none] (axis cs:1.7677669529722873e-01,2.9944770991856217e+03) -- (axis cs:8.8388347648650084e-02,4.5237314724283773e+04) node[midway, above, text=blue, fill=white, fill opacity=0.8, text opacity=1, inner sep=1pt, font=\scriptsize] {3.92};
\draw[draw=none] (axis cs:8.8388347648650084e-02,4.5237314724283773e+04) -- (axis cs:4.4194173824334396e-02,7.1268088530088181e+05) node[midway, above, text=blue, fill=white, fill opacity=0.8, text opacity=1, inner sep=1pt, font=\scriptsize] {3.98};
\draw[draw=none] (axis cs:4.4194173824334396e-02,7.1268088530088181e+05) -- (axis cs:2.2097086912169481e-02,1.1357391425691146e+07) node[midway, above, text=blue, fill=white, fill opacity=0.8, text opacity=1, inner sep=1pt, font=\scriptsize] {3.99};
\addplot+[red, mark=square*, thick] coordinates { (3.5355339059421481e-01,6.2778766495655850e+02) (1.7677669529722873e-01,7.4832087077656770e+03) (8.8388347648650084e-02,1.0303827673061927e+05) (4.4194173824334396e-02,1.5536904580994921e+06) (2.2097086912169481e-02,2.4291271890978795e+07) };
\addlegendentry{$k = 1$}
\draw[draw=none] (axis cs:3.5355339059421481e-01,6.2778766495655850e+02) -- (axis cs:1.7677669529722873e-01,7.4832087077656770e+03) node[midway, above, text=red, fill=white, fill opacity=0.8, text opacity=1, inner sep=1pt, font=\scriptsize] {3.58};
\draw[draw=none] (axis cs:1.7677669529722873e-01,7.4832087077656770e+03) -- (axis cs:8.8388347648650084e-02,1.0303827673061927e+05) node[midway, above, text=red, fill=white, fill opacity=0.8, text opacity=1, inner sep=1pt, font=\scriptsize] {3.78};
\draw[draw=none] (axis cs:8.8388347648650084e-02,1.0303827673061927e+05) -- (axis cs:4.4194173824334396e-02,1.5536904580994921e+06) node[midway, above, text=red, fill=white, fill opacity=0.8, text opacity=1, inner sep=1pt, font=\scriptsize] {3.91};
\draw[draw=none] (axis cs:4.4194173824334396e-02,1.5536904580994921e+06) -- (axis cs:2.2097086912169481e-02,2.4291271890978795e+07) node[midway, above, text=red, fill=white, fill opacity=0.8, text opacity=1, inner sep=1pt, font=\scriptsize] {3.97};
\addplot+[brown, mark=o, thick] coordinates { (3.5355339059421481e-01,1.0787527737688988e+04) (1.7677669529722873e-01,1.1799388592414834e+05) (8.8388347648650084e-02,1.6489321751460987e+06) (4.4194173824334396e-02,2.5107474745989542e+07) (2.2097086912169481e-02,3.9391951673298842e+08) };
\addlegendentry{$k = 2$}
\draw[draw=none] (axis cs:3.5355339059421481e-01,1.0787527737688988e+04) -- (axis cs:1.7677669529722873e-01,1.1799388592414834e+05) node[midway, above, text=brown, fill=white, fill opacity=0.8, text opacity=1, inner sep=1pt, font=\scriptsize] {3.45};
\draw[draw=none] (axis cs:1.7677669529722873e-01,1.1799388592414834e+05) -- (axis cs:8.8388347648650084e-02,1.6489321751460987e+06) node[midway, above, text=brown, fill=white, fill opacity=0.8, text opacity=1, inner sep=1pt, font=\scriptsize] {3.80};
\draw[draw=none] (axis cs:8.8388347648650084e-02,1.6489321751460987e+06) -- (axis cs:4.4194173824334396e-02,2.5107474745989542e+07) node[midway, above, text=brown, fill=white, fill opacity=0.8, text opacity=1, inner sep=1pt, font=\scriptsize] {3.93};
\draw[draw=none] (axis cs:4.4194173824334396e-02,2.5107474745989542e+07) -- (axis cs:2.2097086912169481e-02,3.9391951673298842e+08) node[midway, above, text=brown, fill=white, fill opacity=0.8, text opacity=1, inner sep=1pt, font=\scriptsize] {3.97};
\addplot+[teal, mark=triangle*, thick] coordinates { (3.5355339059421481e-01,6.6467934434429830e+04) (1.7677669529722873e-01,7.0737527151771786e+05) (8.8388347648650084e-02,9.5178351466324534e+06) (4.4194173824334396e-02,1.4300238088929069e+08) (2.2097086912169481e-02,2.2351861587871308e+09) };
\addlegendentry{$k = 3$}
\draw[draw=none] (axis cs:3.5355339059421481e-01,6.6467934434429830e+04) -- (axis cs:1.7677669529722873e-01,7.0737527151771786e+05) node[midway, above, text=teal, fill=white, fill opacity=0.8, text opacity=1, inner sep=1pt, font=\scriptsize] {3.41};
\draw[draw=none] (axis cs:1.7677669529722873e-01,7.0737527151771786e+05) -- (axis cs:8.8388347648650084e-02,9.5178351466324534e+06) node[midway, above, text=teal, fill=white, fill opacity=0.8, text opacity=1, inner sep=1pt, font=\scriptsize] {3.75};
\draw[draw=none] (axis cs:8.8388347648650084e-02,9.5178351466324534e+06) -- (axis cs:4.4194173824334396e-02,1.4300238088929069e+08) node[midway, above, text=teal, fill=white, fill opacity=0.8, text opacity=1, inner sep=1pt, font=\scriptsize] {3.91};
\draw[draw=none] (axis cs:4.4194173824334396e-02,1.4300238088929069e+08) -- (axis cs:2.2097086912169481e-02,2.2351861587871308e+09) node[midway, above, text=teal, fill=white, fill opacity=0.8, text opacity=1, inner sep=1pt, font=\scriptsize] {3.97};
\end{loglogaxis}
\end{tikzpicture}
    \end{minipage}
    \begin{minipage}{0.48\textwidth}
\noindent\begin{tikzpicture}
\begin{loglogaxis}[
  cycle list shift=1,
  ymin=1e1,
  ymax=1e11,
  enlarge y limits=false,
  width=\linewidth, height=0.78\linewidth,
  xtick=data, xticklabel={\begingroup\pgfkeys{/pgf/fpu=true}\pgfmathparse{exp(\tick)}\pgfmathprintnumber[sci,precision=2]{\pgfmathresult}\endgroup},
  title style={font=\small}, label style={font=\scriptsize},
  tick label style={font=\scriptsize}, grid=both,
  major grid style={gray!40}, minor grid style={gray!20},
  legend pos=north east, legend style={font=\scriptsize}, legend cell align={left}]
\addplot+[red, mark=square*, thick] coordinates { (3.5355339059421481e-01,2.2327392743010705e+03) (1.7677669529722873e-01,3.5691455502418503e+04) (8.8388347648650084e-02,5.4370301011595968e+05) (4.4194173824334396e-02,8.4504720403295234e+06) (2.2097086912169481e-02,1.3320508254771371e+08) };
\addlegendentry{$k = 1$}
\draw[draw=none] (axis cs:3.5355339059421481e-01,2.2327392743010705e+03) -- (axis cs:1.7677669529722873e-01,3.5691455502418503e+04) node[midway, above, text=red, fill=white, fill opacity=0.8, text opacity=1, inner sep=1pt, font=\scriptsize] {4.00};
\draw[draw=none] (axis cs:1.7677669529722873e-01,3.5691455502418503e+04) -- (axis cs:8.8388347648650084e-02,5.4370301011595968e+05) node[midway, above, text=red, fill=white, fill opacity=0.8, text opacity=1, inner sep=1pt, font=\scriptsize] {3.93};
\draw[draw=none] (axis cs:8.8388347648650084e-02,5.4370301011595968e+05) -- (axis cs:4.4194173824334396e-02,8.4504720403295234e+06) node[midway, above, text=red, fill=white, fill opacity=0.8, text opacity=1, inner sep=1pt, font=\scriptsize] {3.96};
\draw[draw=none] (axis cs:4.4194173824334396e-02,8.4504720403295234e+06) -- (axis cs:2.2097086912169481e-02,1.3320508254771371e+08) node[midway, above, text=red, fill=white, fill opacity=0.8, text opacity=1, inner sep=1pt, font=\scriptsize] {3.98};
\addplot+[brown, mark=o, thick] coordinates { (3.5355339059421481e-01,3.9428364761624609e+04) (1.7677669529722873e-01,5.5710577496480395e+05) (8.8388347648650084e-02,8.3740522307885988e+06) (4.4194173824334396e-02,1.3005333945413488e+08) (2.2097086912169481e-02,2.0510115999425354e+09) };
\addlegendentry{$k = 2$}
\draw[draw=none] (axis cs:3.5355339059421481e-01,3.9428364761624609e+04) -- (axis cs:1.7677669529722873e-01,5.5710577496480395e+05) node[midway, above, text=brown, fill=white, fill opacity=0.8, text opacity=1, inner sep=1pt, font=\scriptsize] {3.82};
\draw[draw=none] (axis cs:1.7677669529722873e-01,5.5710577496480395e+05) -- (axis cs:8.8388347648650084e-02,8.3740522307885988e+06) node[midway, above, text=brown, fill=white, fill opacity=0.8, text opacity=1, inner sep=1pt, font=\scriptsize] {3.91};
\draw[draw=none] (axis cs:8.8388347648650084e-02,8.3740522307885988e+06) -- (axis cs:4.4194173824334396e-02,1.3005333945413488e+08) node[midway, above, text=brown, fill=white, fill opacity=0.8, text opacity=1, inner sep=1pt, font=\scriptsize] {3.96};
\draw[draw=none] (axis cs:4.4194173824334396e-02,1.3005333945413488e+08) -- (axis cs:2.2097086912169481e-02,2.0510115999425354e+09) node[midway, above, text=brown, fill=white, fill opacity=0.8, text opacity=1, inner sep=1pt, font=\scriptsize] {3.98};
\addplot+[teal, mark=triangle*, thick] coordinates { (3.5355339059421481e-01,3.9501933060741733e+05) (1.7677669529722873e-01,5.5407831727638720e+06) (8.8388347648650084e-02,8.2976658800759092e+07) (4.4194173824334396e-02,1.2869726189588392e+09) (2.2097086912169481e-02,2.0288666393415962e+10) };
\addlegendentry{$k = 3$}
\draw[draw=none] (axis cs:3.5355339059421481e-01,3.9501933060741733e+05) -- (axis cs:1.7677669529722873e-01,5.5407831727638720e+06) node[midway, above, text=teal, fill=white, fill opacity=0.8, text opacity=1, inner sep=1pt, font=\scriptsize] {3.81};
\draw[draw=none] (axis cs:1.7677669529722873e-01,5.5407831727638720e+06) -- (axis cs:8.8388347648650084e-02,8.2976658800759092e+07) node[midway, above, text=teal, fill=white, fill opacity=0.8, text opacity=1, inner sep=1pt, font=\scriptsize] {3.90};
\draw[draw=none] (axis cs:8.8388347648650084e-02,8.2976658800759092e+07) -- (axis cs:4.4194173824334396e-02,1.2869726189588392e+09) node[midway, above, text=teal, fill=white, fill opacity=0.8, text opacity=1, inner sep=1pt, font=\scriptsize] {3.96};
\draw[draw=none] (axis cs:4.4194173824334396e-02,1.2869726189588392e+09) -- (axis cs:2.2097086912169481e-02,2.0288666393415962e+10) node[midway, above, text=teal, fill=white, fill opacity=0.8, text opacity=1, inner sep=1pt, font=\scriptsize] {3.98};
\end{loglogaxis}
\end{tikzpicture}
    \end{minipage}
  }
  \caption{Spectral-norm condition numbers of the statically condensed systems as a function of the mesh size $h$ for the potential-based (left column) and Hessian-based (right column) stabilizations. The numbers along the curves are the observed growth exponents $p$ in $h^{-p}$.}
  \label{fig:condition-numbers}
\end{figure}
\endgroup

\FloatBarrier


\section{Analysis}
\label{sec:analysis}

\subsection{Properties of the local operators}

For every $T\in\Th$, $\underline v_T\in\Htwo{k}{T}$,
and $E\in\ET$, the identities below are obtained by testing the defining relations, using Remark \ref{rem:validity.formulae}, integrating by parts, and observing that polynomials of degree $\le \ell$ on $E$ are uniquely determined by their projection of degree $\le \ell-2$ and their endpoint values:
for $k\ge1$,
\begin{equation}
  \begin{gathered}
  \text{%
    $(G_E^{k+1}\underline v_E)'= H_E^k\underline v_E$
    and
    $G_E^{k+1}\underline v_E(x_V) = G_{v,V}$
    for all $V \in \VE$
  },
  \\
  \text{%
    $(\gamma_E^{k+2}\underline v_E)' = G_E^{k+1}\underline v_E\cdot t_E$
    and $\gamma_E^{k+2}\underline v_E(x_V) =v_V$  for all $V \in \VE$.
  }
  \end{gathered}
  \label{eq:edge-operators:characterization}
\end{equation}
For $k=0$, the first line and the endpoint values of the trace
remain valid, but
$(\gamma_E^2\underline v_E)'=G_E^1\underline v_E\cdot t_E$
need not hold for arbitrary discrete data.
Moreover, for every $k\ge0$, the $L^2$-orthogonal projections of the edge reconstructions coincide with the corresponding polynomial components:
\begin{equation}
  \lproj{k-1}{E}
    (G_E^{k+1}\underline v_E\cdot n_E)=G_{v,E,n},
  \qquad
  \lproj{k-2}{E}\gamma_E^{k+2}\underline v_E=v_E.
  \label{eq:edge-operators-dof-recovery}
\end{equation}

All the reconstructions are polynomially consistent up to degree $k+2$, i.e., for every $T\in\Th$ and every $v\in\Poly{k+2}(T)$, we have
\[
\begin{gathered}
  \text{%
    $\gamma_E^{k+2}\IHtwo{k}{E}v=v_{|E}$,
    $G_E^{k+1}\IHtwo{k}{E}v=(\nabla v)_{|E}$,
    and $H_E^k\IHtwo{k}{E}v=(\Hess v)_{|E}\,t_E$
    for all $E \in \ET$,
  }
  \\
  \text{%
    $H_T^k\IHtwo{k}{T}v=\Hess v$
    and $p_T^{k+2}\IHtwo{k}{T}v=v$.
  }
\end{gathered}
\]

Polynomial consistency, together with boundedness, results in optimal approximation properties for smooth functions; see, e.g., \cite[Lemma 1.43]{Di-Pietro.Droniou:20}.
Specifically, for every $k\ge0$, $T\in\Th$, $v\in H^{k+3}(T)$, and $E\in\ET$, the edge reconstructions satisfy
\begin{equation}
  h_T^{\frac12}
  \norm{L^2(E)}{\gamma_E^{k+2}\IHtwo{k}{E}v-v}
  +h_T^{\frac32}
  \norm{L^2(E;\Real^2)}
    {G_E^{k+1}\IHtwo{k}{E}v-\nabla v}
  \lesssim h_T^{k+3}\seminorm{H^{k+3}(T)}{v}.
  \label{eq:trace-operators:approximation}
\end{equation}
Moving to the element reconstructions, for all integers $m\in\{0,\ldots,k+1\}$, we have
\begin{equation}
  \seminorm{H^m(T;\Real^{2\times2})}
    {H_T^k\IHtwo{k}{T}v-\Hess v}
  \lesssim
  h_T^{k+1-m}\seminorm{H^{k+3}(T)}{v}
  \label{eq:hessian:approximation}
\end{equation}
and, for all integers $m\in\{0,\ldots,k+3\}$,
\begin{equation}
  \seminorm{H^m(T)}{p_T^{k+2}\IHtwo{k}{T}v-v}
  \lesssim
  h_T^{k+3-m}\seminorm{H^{k+3}(T)}{v}.
  \label{eq:potential:approximation}
\end{equation}
Trace approximation properties follow, as in \cite[Theorem 1.45]{Di-Pietro.Droniou:20}, from continuous trace inequalities and the above bounds for two consecutive values of $m$.

Combining the approximation properties \eqref{eq:trace-operators:approximation} of the trace operators and \eqref{eq:potential:approximation} of the potential yields
\begin{equation}
  s_T(\IHtwo{k}{T}v, \IHtwo{k}{T}v)^{\frac12}
  \lesssim h_T^{k+1} \seminorm{H^{k+3}(T)}{v}.
  \label{eq:stabilization:consistency}
\end{equation}

\subsection{Stability}

\begin{proposition}[Properties of the element potential reconstruction]
  \label{prop:potential-low-degree-extension}
  For every $\underline{v}_T\in\Htwo{k}{T}$,
  \begin{equation}
    \pi_T^{k-2}p_T^{k+2}\underline{v}_T=v_T.
    \label{eq:potential-element-moments}
  \end{equation}
  Moreover, for every $\tau\in\Poly{k}(T;\Real^{2\times2})$,
  \begin{equation}
    \begin{aligned}
      \int_Tp_T^{k+2}\underline{v}_T\,\dd\tau
      &=
      \int_T H_T^{k}\underline{v}_T:\tau
      -\sum_{E\in\ET} \omega_{TE} \int_E G_E^{k+1} \underline{v}_E \cdot (\tau n_E)
      \\
      &\quad
      + \sum_{E\in\ET} \omega_{TE} \int_E \gamma_E^{k+2} \underline{v}_E \, \operatorname{div}\tau\cdot n_E.
    \end{aligned}
    \label{eq:potential-low-degree-extension}
  \end{equation}
\end{proposition}

\begin{proof}
For $k=0,1$, the element-moment identity \eqref{eq:potential-element-moments} is void.
  Let $k\ge2$ and $q\in\Poly{k-2}(T)$.
  Since $\dd:\divdivKernelComplement{k}(T)\longrightarrow\Poly{k-2}(T)$ is an isomorphism, there is a unique $\tau_q \in \divdivKernelComplement{k}(T)$ such that $\dd\tau_q=q$.
  In particular, $\tau_q \in \divdivKernelComplement{k+4}(T)$, so it is an admissible test function in the definition \eqref{eq:hessian-potential} of the Hessian potential.
  Using the definition \eqref{eq:element-hessian} of the element Hessian with the same test function on the right-hand side gives
  \[
    \int_T p_T^{k+2}\underline{v}_T \, q
    \overset{\eqref{eq:hessian-potential},\eqref{eq:element-hessian}}{=}
    \int_T v_T \, q.
  \]
  Since this holds for every $q\in\Poly{k-2}(T)$, it proves the element-moment identity \eqref{eq:potential-element-moments}.
Now let $\tau\in\Poly{k}(T;\Real^{2\times2})$.
  Then $\dd\tau\in\Poly{k-2}(T)$, and the element-moment identity \eqref{eq:potential-element-moments} yields
  \[
    \int_Tv_T\,\dd\tau
    =\int_Tp_T^{k+2}\underline{v}_T\,\dd\tau.
  \]
  Replacing the element term in the definition \eqref{eq:element-hessian} of the element Hessian by this expression and rearranging proves the extended potential identity \eqref{eq:potential-low-degree-extension}.
\end{proof}

\begin{proposition}[Estimate of the Hessian of the potential reconstruction]
  For all $T \in \Th$ and all $\underline{v}_T \in \Htwo{k}{T}$, we have
  \begin{equation}\label{eq:estimate-hess-potential}
    \norm{L^2(T; \Real^{2\times 2})}{\Hess p_T^{k+2} \underline{v}_T}
    \lesssim \norm{a,T}{\underline{v}_T}
    \eqcolon a_T(\underline{v}_T, \underline{v}_T)^{\frac12}.
  \end{equation}
\end{proposition}

\begin{proof}
  Applying the continuous integration-by-parts identity \eqref{eq:ibp-hessian} with $v=p_T^{k+2}\underline v_T$ to the left-hand side of \eqref{eq:potential-low-degree-extension} and rearranging, we obtain, for all $\tau \in \Poly{k}(T; \mathbb{S})$,
  \[
  \begin{aligned}
    \int_T \Hess p_T^{k+2} \underline{v}_T : \tau
    &= \int_T H_T^{k} \underline{v}_T : \tau
    - \sum_{E \in \ET} \omega_{TE} \int_E (G_E^{k+1}\underline{v}_E - \nabla p_T^{k+2} \underline{v}_T) \cdot (\tau n_E)
    \\
    &\quad
    + \sum_{E\in\ET} \omega_{TE} \int_E (\gamma_E^{k+2} \underline{v}_E - p_T^{k+2} \underline{v}_T) \, \operatorname{div}\tau\cdot n_E.
  \end{aligned}
  \]
  The estimate follows by selecting $\tau = \Hess p_T^{k+2} \underline{v}_T \in \Poly{k}(T; \mathbb{S})$ and applying weighted Cauchy--Schwarz and trace inequalities.
\end{proof}

\begin{lemma}[Discrete Poincar\'e inequality]
    For all $\underline{v}_h \in \Htwo{k}{h,0}$, we have
    \begin{equation}\label{eq:discrete-poincare}
      \norm{L^2(\Omega)}{p_h^{k+2} \underline{v}_h}
      \lesssim \norm{a,h}{\underline{v}_h}.
    \end{equation}
  \end{lemma}

  \begin{proof}
    Let $\underline{v}_h \in \Htwo{k}{h,0}$ and observe that, by the homogeneous boundary conditions, both $\gamma_E^{k+2} \underline{v}_E$ and $G_E^{k+1} \underline{v}_E$ are zero on boundary edges.
    Apply the discrete Poincar\'e inequality in HHO spaces \cite[Lemma 2.15]{Di-Pietro.Droniou:20} to the vector
    $
    ( (p_T^{k+2} \underline{v}_T)_{T \in \Th}, (\gamma_E^{k+2} \underline{v}_E)_{E \in \Eh})
    $
    to infer
    \[
    \norm{L^2(\Omega)}{p_h^{k+2} \underline{v}_h}
    \lesssim \left(
      \sum_{T \in \Th}
      \norm{L^2(T; \Real^2)}{\nabla p_T^{k+2} \underline{v}_T}^2
      + \sum_{T \in \Th} h_T^{-1} \sum_{E \in \ET} \norm{L^2(E)}{\gamma_E^{k+2} \underline{v}_E - p_T^{k+2} \underline{v}_T}^2
      \right)^{\frac12}.
    \]
    To estimate the first term, apply the Poincar\'e inequality in HHO spaces again to the vector
    $
    ( (\nabla p_T^{k+2} \underline{v}_T)_{T \in \Th}, (G_E^{k+1} \underline{v}_E)_{E \in \Eh})
    $
    to get
    \[
    \norm{L^2(\Omega)}{p_h^{k+2} \underline{v}_h}
    \lesssim \bigg(
      \begin{aligned}[t]
        &\sum_{T \in \Th} \norm{L^2(T; \Real^{2\times 2})}{\Hess p_T^{k+2} \underline{v}_T}^2
       + \sum_{T \in \Th} h_T^{-1} \sum_{E \in \ET} \norm{L^2(E; \Real^2)}{G_E^{k+1} \underline{v}_E - \nabla p_T^{k+2} \underline{v}_T}^2
        \\
        &+ \sum_{T \in \Th} h_T^{-1} \sum_{E \in \ET} \norm{L^2(E)}{\gamma_E^{k+2} \underline{v}_E - p_T^{k+2} \underline{v}_T}^2
        \bigg)^{\frac12}.
      \end{aligned}
      \]
      The conclusion follows by using \eqref{eq:estimate-hess-potential} to estimate the first term,
      observing that $h_T^{-1} \le \operatorname{diam}(\Omega)^2 h_T^{-3}$ for the last term,
      and recalling the definition \eqref{eq:discrete-energy-seminorm} of the $\norm{a,h}{\cdot}$ norm.
  \end{proof}

\begin{corollary}[Norm property]
  \label{cor:norm-property}
  The map $\norm{a,h}{\cdot}$ defines a norm on $\Htwo{k}{h,0}$.
\end{corollary}

\begin{proof}
  Since $\norm{a,h}{\cdot}$ is a seminorm, it suffices to prove that, for all $\underline{v}_h \in \Htwo{k}{h,0}$, $\norm{a,h}{\underline{v}_h} = 0$ implies $\underline{v}_h = \underline{0}$.
  Let $\underline{v}_h \in \Htwo{k}{h,0}$ be such that $\norm{a,h}{\underline{v}_h} = 0$.
  By the discrete Poincar\'e inequality \eqref{eq:discrete-poincare}, this implies $p_T^{k+2} \underline{v}_T = 0$ for all $T \in \Th$ and, by the vanishing stabilization, $\gamma_E^{k+2} \underline{v}_E = 0$ and $G_E^{k+1} \underline{v}_E = 0$ for all $E \in \Eh$.
  The projection property \eqref{eq:potential-element-moments} gives $v_T=0$ for all $T\in\Th$.
Finally, the endpoint relations in \eqref{eq:edge-operators:characterization}, valid also at $k=0$, and \eqref{eq:edge-operators-dof-recovery} give $v_E=0$, $G_{v,E,n}=0$, $v_V=0$, and $G_{v,V}=0$ on every edge and at every vertex.
  This concludes the proof.
\end{proof}

\subsection{Estimate of the consistency error}
\label{sec:consistency-error}

\begin{lemma}[Estimate of the consistency error]
  \label{lem:consistency-error}
  For all $w \in H^2(\Omega) \cap H^4(\Th)$ such that $\Delta^2w\in L^2(\Omega)$ in the distributional sense, we have
  \[
  \norm{a,h,*}{\mathfrak E_h(w;\cdot)}
  \lesssim \bigg[
    \sum_{T \in \Th} \bigg(
    \begin{aligned}[t]
      &\norm{L^2(T; \Real^{2 \times 2})}{\Hess w - H_T^{k} \IHtwo{k}{T} w}^2
      + h_T \norm{L^2(\partial T; \Real^{2 \times 2})}{\Hess w - H_T^{k} \IHtwo{k}{T} w}^2
      \\
      &+ h_T^3 \norm{L^2(\partial T; \Real^2)}{\operatorname{div}(\Hess w - H_T^{k} \IHtwo{k}{T} w)}^2
      + s_T(\IHtwo{k}{T} w, \IHtwo{k}{T} w)
      \bigg)
      \bigg]^{\frac12},
    \end{aligned}
    \]
    where the consistency error $\mathfrak E_h(w;\cdot)$ and its dual norm are defined in \eqref{eq:consistency-functional} and \eqref{eq:consistency-functional-norm}, respectively.
\end{lemma}

\begin{proof}
  We start by noting that global $H^2$ regularity ensures that the traces of $\nabla w$ agree across interfaces, as do their tangential derivatives. Thus, the tangential component of $(\Hess w) n_E$ is continuous. Testing the distributional equation with smooth functions supported near an interior edge yields continuity of the normal-normal Hessian component and of $(\operatorname{div}\Hess w)\cdot n_E+\partial_{t_E}(t_E\cdot(\Hess w) n_E)$. Since the second term in this expression has already been proved to be continuous, we infer the continuity of $(\operatorname{div}\Hess w)\cdot n_E$.

  Now let $\underline{v}_h \in \Htwo{k}{h,0}$.
  We use \eqref{eq:potential-low-degree-extension} to reformulate the consistency term in $a_h$.
  Since the quantities $(\Hess w) n_E$ and $(\operatorname{div} \Hess w) \cdot n_E$ are single-valued at interfaces, and both $G_E^{k+1} \underline{v}_E$ and $\gamma_E^{k+2} \underline{v}_E$ vanish on boundary edges, we can insert these two reconstructions into the boundary terms to obtain
  \[
  \begin{aligned}
    \mathfrak E_h(w;\underline{v}_h)
    &= \sum_{T \in \Th} \int_T \dd (\Hess w - H_T^{k} \IHtwo{k}{T} w) \, p_T^{k+2} \underline{v}_T
    \\
    &\quad
    + \sum_{T \in \Th} \sum_{E \in \ET} \omega_{TE} \int_E (\Hess w - H_T^{k} \IHtwo{k}{T} w) n_E \cdot G_E^{k+1} \underline{v}_E
    \\
    &\quad
    - \sum_{T \in \Th} \sum_{E \in \ET} \omega_{TE} \int_E \operatorname{div} (\Hess w - H_T^{k} \IHtwo{k}{T} w) \cdot n_E \, \gamma_E^{k+2} \underline{v}_E
    \\
    &\quad
    - \sum_{T \in \Th} s_T(\IHtwo{k}{T} w, \underline{v}_T).
  \end{aligned}
  \]
  Applying \eqref{eq:ibp-hessian} for all $T \in \Th$ with $\tau = \Hess w - H_T^{k} \IHtwo{k}{T} w$ and $v = p_T^{k+2} \underline{v}_T$, we arrive at
  \[
  \begin{aligned}
    \mathfrak E_h(w;\underline{v}_h)
    &= \sum_{T \in \Th} \int_T (\Hess w - H_T^{k} \IHtwo{k}{T} w) : \Hess p_T^{k+2} \underline{v}_T
    \\
    &\quad
    + \sum_{T \in \Th} \sum_{E \in \ET} \omega_{TE} \int_E (\Hess w - H_T^{k} \IHtwo{k}{T} w) n_E \cdot (G_E^{k+1} \underline{v}_E - \nabla p_T^{k+2} \underline{v}_T)
    \\
    &\quad
    - \sum_{T \in \Th} \sum_{E \in \ET} \omega_{TE} \int_E \operatorname{div} (\Hess w - H_T^{k} \IHtwo{k}{T} w) \cdot n_E \, (\gamma_E^{k+2} \underline{v}_E - p_T^{k+2} \underline{v}_T)
    \\
    &\quad
    - \sum_{T \in \Th} s_T(\IHtwo{k}{T} w, \underline{v}_T).
  \end{aligned}
  \]
  The result follows by applying Cauchy--Schwarz inequalities and recalling the bound \eqref{eq:estimate-hess-potential} and the definition \mbox{\eqref{eq:discrete-energy-seminorm}} of $\norm{a,h}{\cdot}$.
\end{proof}

\begin{corollary}[Convergence rate of the consistency error for smooth functions]
  \label{cor:adjoint-consistency:convergence}
Let $w\in H^2(\Omega)\cap H^{k+3}(\Th)$ satisfy $\Delta^2w\in L^2(\Omega)$ in the distributional sense, and assume additionally $w\in H^4(\Th)$ if $k=0$. Then,
  \begin{equation}
    \norm{a,h,*}{\mathfrak E_h(w;\cdot)} \lesssim
  \begin{cases}
      \displaystyle
      h \seminorm{H^3(\Th)}{w} + h^2 \seminorm{H^4(\Th)}{w}
      & \text{if $k=0$,}
      \\[1ex]
      h^{k+1}\seminorm{H^{k+3}(\Th)}{w} & \text{if $k\ge1$}.
    \end{cases}
    \label{eq:adjoint-consistency:convergence}
  \end{equation}
\end{corollary}

\begin{proof}
For $k\ge1$, apply continuous trace inequalities to the boundary terms and use the approximation properties \eqref{eq:hessian:approximation} of $H_T^{k}$ and the consistency \eqref{eq:stabilization:consistency} of the stabilization bilinear form. For $k=0$, the same estimates bound the volume Hessian term, its boundary trace, and the stabilization by $\sum_T h_T^2\seminorm{H^3(T)}{w}^2$.
  Since $H_T^0\IHtwo{0}{T}w$ is constant, the remaining boundary term involves only $\operatorname{div}\Hess w$.
  The continuous trace inequality gives
  \begin{equation*}
    h_T^3\norm{L^2(\partial T;\Real^2)}{\operatorname{div}\Hess w}^2
    \lesssim h_T^2\seminorm{H^3(T)}{w}^2
      +h_T^4\seminorm{H^4(T)}{w}^2,
  \end{equation*}
  which proves the stated lowest-order estimate.
\end{proof}


\section{Application to the von K\'arm\'an equations}
\label{sec:von-karman}
\def\vkPotentialCriticalLoad{59.945225}
\def\vkPotentialCriticalKappa{6.073721}
\def\vkPotentialElements{18226}
\def\vkPotentialActiveDofs{54310}
\def\vkPotentialLastChangePercent{0.002811}
\def\vkPotentialPreviousChangePercent{0.006308}
\def\vkPotentialStableDigits{4}
\def\vkPotentialMaxDeflection{2.555\times10^{-5}}
\def\vkPotentialNearLoad{60}
\def\vkPotentialNearDistance{0.0548}
\def\vkPotentialNearLtwoAmplitude{0.19072}
\def\vkPotentialNearModalFraction{0.99994}
\def\vkPotentialNewtonMax{3}
\def\vkPotentialLargerBasisDifference{7.530\times10^{-10}}
\def\vkPotentialShiftDifference{1.182\times10^{-9}}
\def\vkPotentialSymmetricFirstDifference{2.009\times10^{-9}}
\def\vkPotentialSymmetricMaxDifference{3.593\times10^{-9}}
\def\vkPotentialSpectrumFirstMin{59.91239}
\def\vkPotentialSpectrumFirstMax{60.30638}
\def\vkPotentialResponsePlotMin{0.000751194206}
\def\vkPotentialResponsePlotMax{0.238395578103}
\def\vkPotentialResponseLabelY{0.013382129012}
\def\vkPotentialTriWOrder{2.56}
\def\vkPotentialQuadWOrder{2.91}
\def\vkPotentialTriPhiOrder{2.60}
\def\vkPotentialQuadPhiOrder{2.97}
\def\vkPotentialManufacturedResidual{6.398\times10^{-10}}
\def\vkPotentialTriWOrderZero{0.97}
\def\vkPotentialTriWOrderOne{2.56}
\def\vkPotentialTriWOrderTwo{3.00}
\def\vkPotentialQuadWOrderZero{0.99}
\def\vkPotentialQuadWOrderOne{2.91}
\def\vkPotentialQuadWOrderTwo{3.02}
\def\vkPotentialTriPhiOrderZero{0.98}
\def\vkPotentialTriPhiOrderOne{2.60}
\def\vkPotentialTriPhiOrderTwo{2.98}
\def\vkPotentialQuadPhiOrderZero{1.00}
\def\vkPotentialQuadPhiOrderOne{2.97}
\def\vkPotentialQuadPhiOrderTwo{3.00}
\def\vkPotentialManufacturedResidualZero{2.744\times10^{-10}}
\def\vkPotentialManufacturedResidualOne{6.398\times10^{-10}}
\def\vkPotentialManufacturedResidualTwo{2.928\times10^{-6}}
\def\vkPotentialManufacturedMaxResidual{5.606\times10^{-6}}
\def\vkPotentialManufacturedMaxTolerance{1\times10^{-5}}
\def\vkPotentialManufacturedInitializationDifference{1.005\times10^{-8}}

In the following, we extend the DDR method to the classical von K\'arm\'an model \cite{Karman:10}.
This model describes the response of a thin plate in two dimensions and accounts for geometric nonlinearities through a nonlinear strain--displacement relation (see \cite[Sections~5.1--5.7]{Ciarlet:21} for its derivation from three-dimensional elasticity).

The purpose of this section is to describe the model and its implementation, and to present numerical studies of the flat-branch spectrum under compression and of a small forced nonlinear response in the spirit of \cite{Reinhart:82}; see also \cite{Pichi.Rozza:19} concerning the construction of parametrized bifurcation diagrams.
A complete well-posedness and convergence analysis of the nonlinear scheme is left for future work.

In all numerical computations in this section, we use the better-performing potential-based stabilization \eqref{eq:potential-stabilization}.

\subsection{Continuous von K\'arm\'an problem}

Consider a plate whose configuration at rest coincides with the set $\Omega\subset\Real^2$ and whose boundary is split into an outer component $\Gamma_{\rm e}$ and, possibly, the boundary $\Gamma_{\rm h}$ of a hole.
The unknowns are the transverse deflection $w:\Omega\to\Real$ and the Airy stress function $\phi:\Omega\to\Real$.
For sufficiently regular scalar functions $r$ and $s$, set
\begin{equation*}
  [r,s]
  \coloneqq
  \operatorname{tr}(\Hess r)\operatorname{tr}(\Hess s)
  -\Hess r:\Hess s.
\end{equation*}
This is the symmetric von K\'arm\'an bracket, and the associated membrane stress tensor is $ N(\phi)\coloneqq\operatorname{tr}(\Hess\phi) I_2 -\Hess\phi$, where $I_2$ is the $2\times2$ identity matrix.
The following holds:
\begin{equation}
    [r,s] = N(r):\Hess s,\qquad
    \operatorname{div} N(r) = 0.
  \label{eq:von-karman-stress-identities}
\end{equation}
The divergence identity follows from equality of mixed derivatives and expresses the equilibrium of an Airy-generated membrane stress in the absence of in-plane body forces.

Given source terms $q : \Omega \to \Real$ and $g : \Omega \to \Real$, the equations considered here take the nondimensional form
\begin{equation}
  \label{eq:von-karman:strong}
  \begin{alignedat}{2}
    \Delta^2 w-[\phi,w]&=q &\qquad& \text{in $\Omega$,}
    \\
    \Delta^2\phi+\tfrac12[w,w]&=g &\qquad& \text{in $\Omega$.}
  \end{alignedat}
\end{equation}
The outer boundary is clamped for the deflection, whereas the hole is left mechanically free.
The Airy data are prescribed on every boundary component.

To state the weak problem, let
\[
V_{\rm e}
\coloneqq
\{v\in H^2(\Omega):v=\partial_n v=0\text{ on }\Gamma_{\rm e}\},\qquad
V_0
\coloneqq H_0^2(\Omega),
\]
and let $\phi_D\in H^2(\Omega)$ be an extension of the prescribed Airy boundary value and normal derivative.
The variational problem is: find $(w,\phi)\in V_{\rm e}\times(\phi_D+V_0)$ such that
\begin{equation}
  \begin{alignedat}{2}
    a(w,v)-\int_\Omega[\phi,w]v
    &=\int_\Omega q  v
    &\qquad&\forall v\in V_{\rm e},
    \\
    a(\phi,\psi)+\frac12\int_\Omega[w,w]\psi
    &=\int_\Omega g \, \psi
    &\qquad&\forall\psi\in V_0.
  \end{alignedat}
  \label{eq:von-karman:weak}
\end{equation}
Here, $a$ is the Hessian bilinear form already introduced in \eqref{eq:continuous-weak}, which, for isotropic bending, corresponds to unit flexural rigidity and a Poisson ratio of zero.
The geometric form appearing in the first equation has an additional symmetry that will be relevant in the spectral study.
Indeed, if $ N(\phi) n=0$ on $\Gamma_{\rm h}$, then an integration by parts, together with $v=0$ on $\Gamma_{\rm e}$ and \eqref{eq:von-karman-stress-identities}, gives
\begin{equation}
  \int_\Omega[\phi,z]v
  =-\int_\Omega  N(\phi)\nabla z\cdot\nabla v
  =\int_\Omega[\phi,v]z
  \qquad\forall (z,v)\in V_{\rm e}^2.
  \label{eq:von-karman-geometric-symmetry}
\end{equation}
Thus, the continuous linearized transverse operator defined in \eqref{eq:von-karman-linearisation} below is self-adjoint.

\subsection{DDR discretization of the coupled problem}

The DDR scheme for \eqref{eq:von-karman:weak} uses two copies of the space $\Htwo{k}{h}$ defined in \eqref{eq:global-space}, one for each scalar unknown.
All the local operators are likewise the reconstructions introduced above, applied to the appropriate copy.
Accounting for boundary conditions, we define the following subspaces of $\Htwo{k}{h}$:
\begin{equation*}
  \begin{aligned}
    \underline W_{h,{\rm e}}^k
    &\coloneqq
    \left\{
    \underline v_h\in\Htwo{k}{h}:
    \text{all edge and vertex components vanish on }\Gamma_{\rm e}
    \right\},\\
    \underline\Phi_{h,D}^k
    &\coloneqq
    \left\{
    \underline\phi_h\in\Htwo{k}{h}:
    \text{the prescribed Airy boundary data are interpolated on }
    \partial\Omega
    \right\}.
  \end{aligned}
\end{equation*}
For $(\underline r_h,\underline s_h) \in \Htwo{k}{h} \times \Htwo{k}{h}$, define the global discrete bracket $[\underline r_h,\underline s_h]_h\in\Poly{2k}(\Th)$ elementwise by
\begin{equation*}
  \left([\underline r_h,\underline s_h]_h\right)_{|T}
  \coloneqq
  \operatorname{tr}( H_T^{k}\underline r_T)
  \operatorname{tr}( H_T^{k}\underline s_T)
  - H_T^{k}\underline r_T:
    H_T^{k}\underline s_T
  \qquad\forall T\in\Th.
\end{equation*}
The DDR discretization of \eqref{eq:von-karman:weak} reads: find $(\underline w_h,\underline\phi_h)\in \underline W_{h,{\rm e}}^k
\times\underline\Phi_{h,D}^k$ such that
\begin{subequations}
  \label{eq:discrete-von-karman}
  \begin{alignat}{2}
    a_h(\underline w_h,\underline v_h)
    -\int_\Omega
      [\underline\phi_h,\underline w_h]_h \, p_h^{k+2}\underline v_h
    &= \int_\Omega q \, p_h^{k+2} \underline{v}_h
    &\qquad&\forall\underline v_h\in    \underline W_{h,{\rm e}}^k,
    \label{eq:discrete-von-karman-deflection}\\
    a_h(\underline\phi_h,\underline\psi_h)
    +\frac12\int_\Omega
      [\underline w_h,\underline w_h]_h \, p_h^{k+2}\underline\psi_h
    &= \int_\Omega g \, p_h^{k+2} \underline{\psi}_h
    &\qquad&\forall\underline\psi_h\in\Htwo{k}{h,0}.
    \label{eq:discrete-von-karman-compatibility}
  \end{alignat}
\end{subequations}
Although no convergence analysis is carried out here, using the Hessian-based potential $p_h^{k+2}$ (cf. \eqref{eq:global-hessian-potential}) in the terms involving the von K\'arm\'an bracket is consistent with the treatment of the volume load.

\subsection{Convergence test for the von K\'arm\'an problem}
\label{sec:convergence-test-von-karman}

We start by assessing the accuracy of the coupled method using the smooth benchmark of \cite{Shylaja.Kumar:24}, for which $\Gamma_{\rm h} = \emptyset$; see also \cite[Section~5.1]{Lovadina.Mora.ea:21} for a variant.
On the unit square $\Omega=(0,1)^2$, we prescribe
\begin{equation}
  w(x,y)=x^2(1-x)^2y^2(1-y)^2,
  \qquad
  \phi(x,y)=\sin^2(\pi x)\sin^2(\pi y).
  \label{eq:manufactured-von-karman}
\end{equation}
Both fields satisfy homogeneous clamped boundary conditions, namely $w=\partial_n w=\phi=\partial_n\phi=0$ on $\partial\Omega$.
The two source terms are manufactured from the exact solution according to \eqref{eq:von-karman:strong}.

We use $k=0,1,2$ and the triangular and Cartesian mesh families illustrated in \cref{fig:convergence-test-meshes}, with five successive mesh levels for each family.
For each field, we measure the discrete error with respect to its exact interpolant in the energy seminorm \eqref{eq:discrete-energy-seminorm}:
\begin{equation*}
  E_w(h)=\norm{a,h}{\underline w_h-\IHtwo{k}{h}w},
  \qquad
  E_\phi(h)=\norm{a,h}{\underline\phi_h-\IHtwo{k}{h}\phi}.
\end{equation*}
The nonlinear terms and error functionals are evaluated using quadrature of degree $20$.
For $k=0,1$, the Newton residual tolerance is $10^{-9}$, except on the finest Cartesian mesh, where it is $10^{-8}$.
For $k=2$, the absolute residual tolerances on mesh levels $0$ through $4$ are $10^{-9}$, $10^{-9}$, $5\times10^{-8}$, $10^{-6}$, and $10^{-5}$, respectively, to accommodate stagnation of the assembled residual on refined meshes.
The finest-mesh $k=2$ results were independently checked using the exact interpolant as the initial guess; for both energy errors, the relative differences from the results obtained with the linear initial guess are at most $\vkPotentialManufacturedInitializationDifference$.

\Cref{fig:von-karman-convergence} shows that both errors decrease monotonically with refinement. For $k=0$ and $k=2$, the observed orders approach $1$ and $3$, respectively, for both fields.
For $k=1$, on the finest triangular refinement, the observed orders are approximately $\vkPotentialTriWOrder$ for $w$ and $\vkPotentialTriPhiOrder$ for $\phi$; on Cartesian meshes, they reach $\vkPotentialQuadWOrder$ and $\vkPotentialQuadPhiOrder$, respectively.
These enhanced rates for the present smooth solution are consistent with the behavior observed at $k=1$ in the linear tests of \cref{sec:discrete-problem:numerical-validation}.

\begin{figure}[!htbp]
  \centering
  \input{figures/von_karman_convergence_test/combined/figure8.tex}
  \caption{Energy errors $E_w(h)$ (left column) and $E_\phi(h)$ (right column) for the manufactured solution \eqref{eq:manufactured-von-karman} for $k=0,1,2$.
    The top and bottom rows correspond to triangular and Cartesian meshes, respectively.
    The numbers along the curves are the observed orders between successive refinements.}
  \label{fig:von-karman-convergence}
\end{figure}

\subsection{Perforated-plate test case}
\label{sec:perforated-plate-test}

For the numerical experiments presented below, we restrict the discretization to $k=1$.
This choice is motivated by both numerical and geometric considerations.
On the numerical side, it keeps the cost of the repeated spectral and coupled nonlinear solves moderate.
It additionally results in a symmetric Hessian reconstruction; see \cref{rem:HT:symmetry}.
Concerning the geometry, the piecewise-linear approximation of curvilinear boundaries introduces an error that may limit the benefit of higher polynomial degrees.
An extension of the current method to curved edges appears possible, although nontrivial, following the ideas of \cite{Botti.Di-Pietro:18,Yemm:24}.
This topic is left for future work.

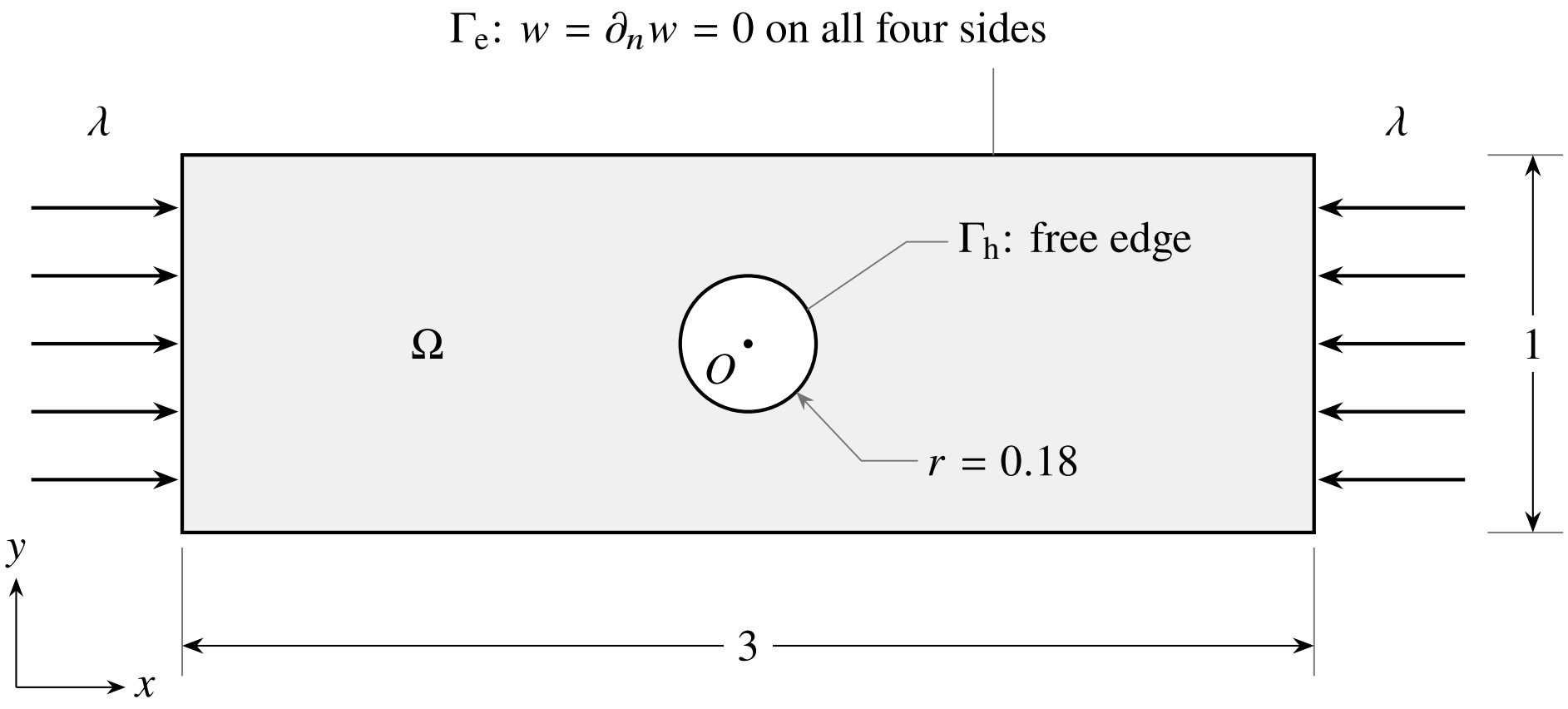
\begin{figure}[!htbp]
  \centering
  \begin{tikzpicture}[
    x=3cm,y=3cm,font=\small,
    boundary/.style={line width=0.8pt},
    load/.style={-{Stealth[length=2mm]},line width=0.8pt},
    dimension/.style={<->,>=Stealth,line width=0.4pt},
    guide/.style={line width=0.35pt,black!55}
  ]
    \path[fill=black!6,even odd rule]
      (-1.5,-0.5) rectangle (1.5,0.5) (0,0) circle[radius=0.18];
    \draw[boundary] (-1.5,-0.5) rectangle (1.5,0.5);
    \draw[boundary] (0,0) circle[radius=0.18];
    \fill (0,0) circle[radius=0.012];
    \node[below left,inner sep=2pt] at (0,0) {$O$};
    \node at (-0.85,0) {$\Omega$};

    \node at (0,0.83) {$\Gamma_{\rm e}$: $w=\partial_n w=0$ on all four sides};
    \draw[guide] (0.65,0.73) -- (0.65,0.5);
    \draw[guide] (0.156,0.09) -- (0.42,0.27) -- (0.53,0.27);
    \node[anchor=west,inner sep=2pt] at (0.53,0.27) {$\Gamma_{\rm h}$: free edge};
    \draw[guide,-{Stealth[length=1.5mm]}]
      (0.45,-0.31) -- (0.3,-0.31) -- (0.1273,-0.1273);
    \node[anchor=west,inner sep=2pt] at (0.45,-0.31) {$r=0.18$};

    \foreach \y in {-0.36,-0.18,0,0.18,0.36} {
      \draw[load] (-1.9,\y) -- (-1.51,\y);
      \draw[load] (1.9,\y) -- (1.51,\y);
    }
    \node at (-1.72,0.59) {$\lambda$};
    \node at (1.72,0.59) {$\lambda$};

    \draw[guide] (-1.5,-0.54) -- (-1.5,-0.88);
    \draw[guide] (1.5,-0.54) -- (1.5,-0.88);
    \draw[dimension] (-1.5,-0.8) -- node[fill=white,inner sep=3pt] {$3$} (1.5,-0.8);
    \draw[guide] (1.96,-0.5) -- (2.16,-0.5);
    \draw[guide] (1.96,0.5) -- (2.16,0.5);
    \draw[dimension] (2.08,-0.5) -- node[fill=white,inner sep=3pt] {$1$} (2.08,0.5);

    \draw[-{Stealth[length=1.5mm]}] (-1.94,-0.91) -- (-1.65,-0.91)
      node[right,inner sep=2pt] {$x$};
    \draw[-{Stealth[length=1.5mm]}] (-1.94,-0.91) -- (-1.94,-0.62)
      node[above,inner sep=2pt] {$y$};
  \end{tikzpicture}
  \caption{Perforated-plate test configuration with $O=(0,0)$.
    The arrows indicate nominal longitudinal compression of magnitude $\lambda$, represented by the outer Airy data \eqref{eq:perforated-plate-airy-data}.
    The hole is free for the deflection and traction-free in plane, with prescribed Airy data $\phi=0$ and $\nabla\phi=0$.}
  \label{fig:perforated-plate-configuration}
\end{figure}

We now consider a more realistic problem involving the deformation of a plate occupying the perforated rectangle
\begin{equation*}
  \Omega=
  \left((-\tfrac32,\tfrac32)\times(-\tfrac12,\tfrac12)\right)
  \setminus\overline{B_{0.18}(0)},
\end{equation*}
with $B_{0.18}(0)$ denoting the ball centered at the origin with radius $0.18$.
The geometry and boundary conditions are illustrated in \cref{fig:perforated-plate-configuration}.
The four straight outer sides form $\Gamma_{\rm e}$, and $\Gamma_{\rm h}=\partial B_{0.18}(0)$ is approximated by polygonal boundaries in the computation, since elements with straight edges are used.
As in \cref{sec:convergence-test-von-karman}, the deflection is clamped only on $\Gamma_{\rm e}$; no deflection DOF is fixed on $\Gamma_{\rm h}$.
We use the following Airy boundary datum:
\begin{equation}
  \phi_D(x,y)=-\frac{\lambda}{2}y^2,
  \qquad
  \nabla\phi_D(x,y)=(0,-\lambda y)^\top
  \quad\text{on }\Gamma_{\rm e},
  \label{eq:perforated-plate-airy-data}
\end{equation}
with compression parameter $\lambda>0$.
This choice generates the membrane stress $ N(\phi_D)=\left(\begin{smallmatrix}-\lambda&0\\0&0\end{smallmatrix}\right)$, corresponding to uniform longitudinal compression; the actual stress inside the perforated plate is nonuniform.
On the hole, we prescribe $\phi=0$ and $\nabla\phi=0$.

For the spectral study below, we set the loading terms $q=g=0$ in \eqref{eq:von-karman:strong}, so that a flat equilibrium exists for every compression parameter.
Seven conforming triangular meshes are obtained by remeshing the same Gmsh\footnote{\url{https://gmsh.info}} geometry, with target sizes $0.20(0.7)^j$, $j=0,\ldots,6$.
The meshes are independently generated, so that the polygonal approximation of the circular hole is also refined.

\subsubsection{Flat-branch spectrum and mesh refinement}
\label{sec:von-karman-results}

The flat branch has $w=0$ and an Airy function $\phi_\lambda^0$ solving the homogeneous biharmonic equation with the prescribed boundary data \eqref{eq:perforated-plate-airy-data}.
Since these data scale with $\lambda$, one has, by linearity,
\[
  \phi_\lambda^0=\lambda\phi_1^0.
\]
Buckling occurs when the plate under compression starts to deflect from its flat configuration.
It can be investigated through the linearization of the von K\'arm\'an equations along the flat branch, following the classical analysis of plate buckling in
\cite{Berger:67}.
To obtain the linearized problem, consider perturbations of the flat
branch of the form
\[
  w=\varepsilon z,
  \qquad
  \phi=\phi_\lambda^0+\varepsilon\chi,
  \qquad
  0<\varepsilon\ll1.
\]
After substituting the perturbations into the weak formulation
\eqref{eq:von-karman:weak} and neglecting higher-order terms in
$\varepsilon$, we obtain two decoupled equations. The Airy perturbation
satisfies a homogeneous biharmonic problem, whereas the transverse
linearization defines the operator
$L_\lambda:V_{\rm e}\to V_{\rm e}'$ by
\begin{equation}
  \langle L_\lambda z,v\rangle
  =
  a(z,v)
  -\lambda\int_\Omega[\phi_1^0,z]v
  \qquad\forall z,v\in V_{\rm e},
  \label{eq:von-karman-linearisation}
\end{equation}
where $\langle\cdot,\cdot\rangle$ denotes the duality pairing between $V_{\rm e}'$ and $V_{\rm e}$.
Buckling of the flat configuration can occur when the linearized
operator $L_\lambda$ has a nontrivial kernel. We therefore identify
as candidate buckling loads the positive values of $\lambda$ for which
\[
  \ker L_\lambda\neq\{0\},
\]
with the corresponding nonzero functions representing infinitesimal
transverse buckling modes.

For the buckling study, on each mesh we first compute the unit-load Airy solution $\underline\phi_{1,h}^0$ using \eqref{eq:discrete-von-karman-compatibility} with $\underline w_h=0$, the boundary data at $\lambda=1$, and no loading term.
We then set
\begin{equation*}
  g_h(\underline z_h,\underline v_h)
  \coloneqq
  \int_\Omega
    [\underline\phi_{1,h}^0,\underline z_h]_h\,
    p_h^{k+2}\underline v_h.
\end{equation*}
The discrete spectral problem consists of finding $\lambda_h \in \Real$ and $\underline z_h\in\underline W_{h,{\rm e}}^k\setminus\{\underline{0}\}$ such that
\begin{equation}
  a_h(\underline z_h,\underline v_h)
  =\lambda_h g_h(\underline z_h,\underline v_h)
  \qquad\forall\underline v_h\in\underline W_{h,{\rm e}}^k.
  \label{eq:von-karman-eigenproblem}
\end{equation}
The algebraic form of the above problem, after elimination of Dirichlet DOFs, is $K_h x=\lambda_h G_h x$.

Unlike its continuous counterpart in \eqref{eq:von-karman-geometric-symmetry}, the form $g_h$ is not exactly symmetric because its second argument is represented by the reconstructed potential $p_h^{k+2}\underline v_h$.
Its asymmetry is measured here by
\begin{equation}
  \alpha_h\coloneqq
  \frac{\norm{F}{G_h-G_h^\top}}{\norm{F}{G_h}},
  \label{eq:von-karman-geometric-asymmetry}
\end{equation}
where $\norm{F}{\cdot}$ denotes the Frobenius norm.

%

We compute six positive real eigenpairs on each mesh of the refined sequence described in Section \ref{sec:perforated-plate-test}.
Consistent with the preceding discussion, $G_h$ is neither assumed to be symmetric nor symmetrized.
An unrestarted Arnoldi iteration \cite[Chapter~6]{Saad:11} with a basis of dimension $128$ is applied to $K_h^{-1}G_h$, using a single sparse LU factorization of $K_h$.
For each eigenpair $(\lambda_h, x)$, we consider the relative residual
\begin{equation}
  \eta_h(\lambda_h, x)
  \coloneqq
  \frac{\norm{2}{K_h x-\lambda_hG_h x}}{\norm{2}{K_h x}+|\lambda_h|\norm{2}{G_h x}}.
  \label{eq:eigenpair-residual}
\end{equation}
The acceptance threshold is $10^{-8}$.
A Ritz value is classified as real when the magnitude of its imaginary part is at most $10^{-9}\max\{1,|\lambda_h|\}$.

On the finest mesh, the algebraic accuracy of the six loads is independently checked by recomputing them with a basis of dimension $192$ and a different starting vector; a shifted solve provides a further check of the first mode.
The larger-basis check changes the six loads by $\vkPotentialLargerBasisDifference$, and the shifted check changes the first by $\vkPotentialShiftDifference$.
As a separate structural check, on the finest mesh we replace $G_h$ by its symmetric part $\frac{G_h+G_h^\top}{2}$ and repeat the six-mode computation with the larger basis and an independent starting vector.
The relative change in the first load is $\vkPotentialSymmetricFirstDifference$, and the maximum relative change among the six loads is $\vkPotentialSymmetricMaxDifference$.
Thus, although the measured asymmetry is not at roundoff level, its influence on the reported part of the spectrum can be considered negligible.

Write $\lambda_{1,h}\le\cdots\le\lambda_{6,h}$ for the six computed positive loads.
For $i=1,\ldots,6$, we denote by $\underline z_{i,h}$ a corresponding discrete eigenvector satisfying \eqref{eq:von-karman-eigenproblem} with $\lambda_h=\lambda_{i,h}$.
We denote the first load by $\lambda_{\mathrm{cr},h}=\lambda_{1,h}$.
\Cref{tab:von-karman-results} reports the first load, the geometric asymmetry indicator, and the largest residual among the six modes.
The number of deflection DOFs after strong enforcement of the Dirichlet conditions is denoted by $n_{w,h}$.

\def\vkObservedEigenOrder{1.800}
\def\vkExtrapolatedCriticalLoad{59.94326}
\def\vkExtrapolatedCriticalKappa{6.07352}
\def\vkFinestExtrapolationDifference{0.00196}
\def\vkCoarsestGeometricAsymmetry{2.814\times10^{-2}}
\def\vkFinestGeometricAsymmetry{7.857\times10^{-3}}
\def\vkSymmetrizedFirstRelativeDifference{2.009e-09}
\def\vkSymmetrizedMaximumRelativeDifference{3.593e-09}

\begin{table}[!htbp]
  \centering
  \normalsize
  \begin{tabular}{@{}crrrrr@{}}
\toprule
Mesh & $h$ & $n_{w,h}$ & $\lambda_{\mathrm{cr},h}$ & $\alpha_h$ & $\max_{1\le i\le6}\eta_h(\lambda_{i,h},\boldsymbol x_{i,h})$ \\
\midrule
    0 & $2.46\times10^{-1}$ & 862 & 60.273551 & $2.8\times10^{-2}$ & $4.0\times10^{-13}$ \\
    1 & $1.72\times10^{-1}$ & 1676 & 60.066450 & $2.1\times10^{-2}$ & $8.0\times10^{-13}$ \\
    2 & $1.33\times10^{-1}$ & 3254 & 59.991776 & $1.8\times10^{-2}$ & $3.3\times10^{-12}$ \\
    3 & $8.94\times10^{-2}$ & 6548 & 59.962153 & $1.4\times10^{-2}$ & $1.1\times10^{-11}$ \\
    4 & $6.53\times10^{-2}$ & 13168 & 59.950692 & $1.2\times10^{-2}$ & $3.9\times10^{-11}$ \\
    5 & $4.40\times10^{-2}$ & 26648 & 59.946910 & $9.7\times10^{-3}$ & $1.6\times10^{-10}$ \\
    6 & $3.12\times10^{-2}$ & 54310 & 59.945225 & $7.9\times10^{-3}$ & $6.6\times10^{-10}$ \\
\bottomrule
\end{tabular}

  \caption{%
    Flat-branch eigenvalue computation on seven meshes.
    The loads are obtained directly from \eqref{eq:von-karman-eigenproblem};
    $\alpha_h$ is defined by \eqref{eq:von-karman-geometric-asymmetry} and $\eta_h$ by \eqref{eq:eigenpair-residual}.%
  }
  \label{tab:von-karman-results}
\end{table}

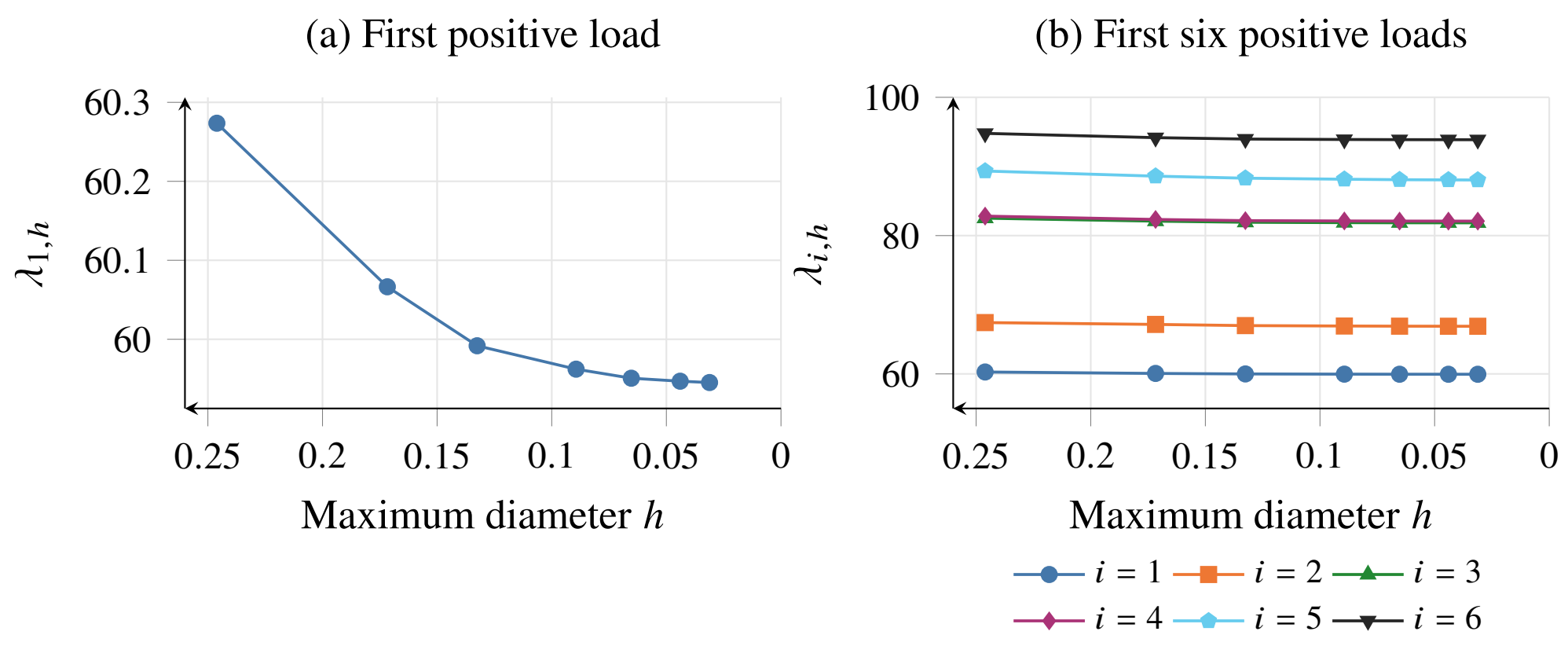
\begin{figure}[!htbp]
  \centering
  \begin{tikzpicture}
    \definecolor{vkblue}{HTML}{4477AA}
    \definecolor{vkorange}{HTML}{EE7733}
    \definecolor{vkgreen}{HTML}{228833}
    \definecolor{vkpurple}{HTML}{AA3377}
    \definecolor{vkcyan}{HTML}{66CCEE}
    \pgfplotsset{
      vkmode1/.style={color=vkblue,mark=*},
      vkmode2/.style={color=vkorange,mark=square*},
      vkmode3/.style={color=vkgreen,mark=triangle*},
      vkmode4/.style={color=vkpurple,mark=diamond*},
      vkmode5/.style={color=vkcyan,mark=pentagon*},
      vkmode6/.style={color=black!85,mark=triangle*,mark options={rotate=180}}
    }
    \begin{groupplot}[
      group style={group size=2 by 1,horizontal sep=1.45cm},
      width=0.44\textwidth,height=4.2cm,
      axis lines=left,tick align=outside,
      label style={font=\small},tick label style={font=\footnotesize},
      title style={font=\small},grid=major,grid style={black!10},
      xlabel={Maximum diameter $h$},xmin=0,xmax=0.26,
      x dir=reverse,
      xtick={0,0.05,0.1,0.15,0.2,0.25},
      xticklabel style={/pgf/number format/fixed,/pgf/number format/precision=2},
      every axis plot/.append style={line width=0.7pt,mark size=1.7pt}
    ]
      \nextgroupplot[
        title={(a) First positive load},
        ymin=\vkPotentialSpectrumFirstMin,ymax=\vkPotentialSpectrumFirstMax,
        ylabel={$\lambda_{1,h}$},
        scaled y ticks=false,
        yticklabel style={/pgf/number format/fixed,/pgf/number format/precision=2}
      ]
      \addplot[vkmode1] table[x=h,y=lambda1]
        {data/von_karman_buckling/eigenvalues.dat};
      \nextgroupplot[
        title={(b) First six positive loads},ymin=55,ymax=100,
        ylabel={$\lambda_{i,h}$},
        legend style={font=\scriptsize,draw=none,fill=none,
          at={(0.5,-0.43)},anchor=north,legend columns=3}
      ]
      \pgfplotsinvokeforeach{1,...,6}{
        \addplot[vkmode#1] table[x=h,y=lambda#1]
          {data/von_karman_buckling/eigenvalues.dat};
        \addlegendentry{$i=#1$}
      }
    \end{groupplot}
  \end{tikzpicture}
  \caption{%
    Refinement of the flat-branch spectrum for $k=1$: the first load on a magnified scale (a) and the first six loads (b).
  }
  \label{fig:von-karman-spectrum}
\end{figure}

\begin{figure}[!t]
  \centering
  \input{figures/von_karman_buckling_labels.tex}
  \begin{subfigure}{\textwidth}
    \centering
    \includegraphics[width=\textwidth]{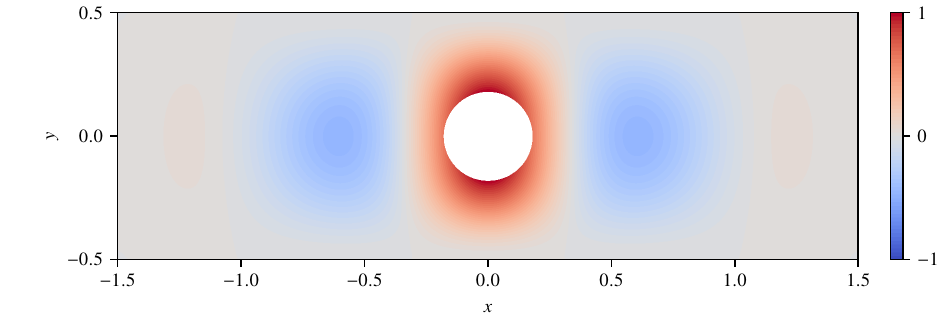}
    \caption{First eigenmode, $\lambda_{1,h}=\vkFirstModeLoad$.}
    \label{fig:von-karman-mode-one}
  \end{subfigure}\par\medskip
  \begin{subfigure}{\textwidth}
    \centering
    \includegraphics[width=\textwidth]{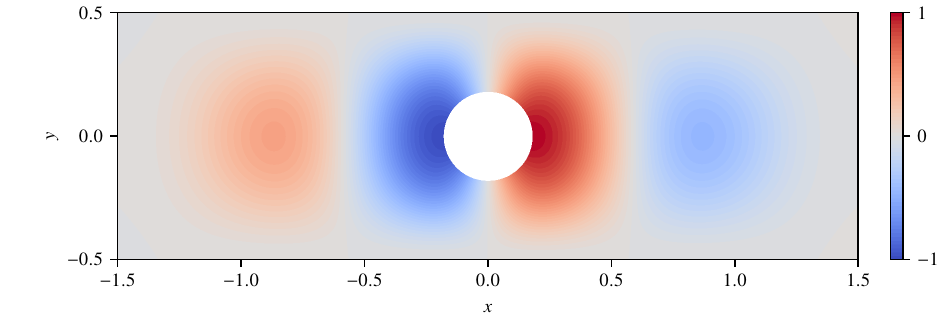}
    \caption{Second eigenmode, $\lambda_{2,h}=\vkSecondModeLoad$.}
    \label{fig:von-karman-mode-two}
  \end{subfigure}\par\medskip
  \begin{subfigure}{\textwidth}
    \centering
    \includegraphics[width=\textwidth]{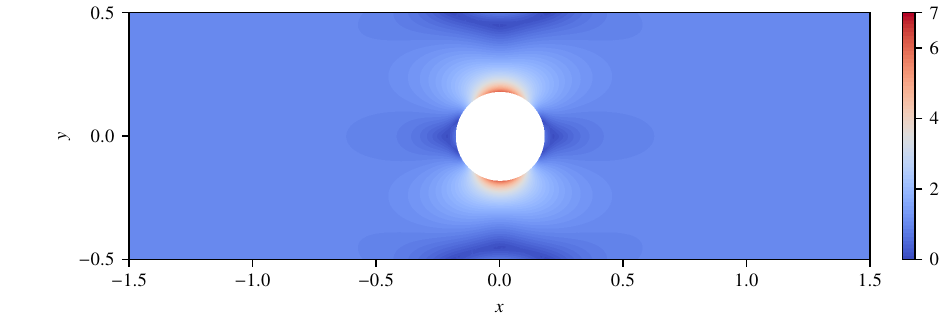}
    \caption{Spectral norm of the unit-load membrane stress, $\| N_h^0\|_2$.}
    \label{fig:von-karman-unit-stress}
  \end{subfigure}
  \caption{%
    Reconstructed fields on mesh \vkFieldMesh{} at degree $k=1$.
    The first two eigenmodes are separately normalized by their sampled maximum absolute deflections.
  }
  \label{fig:von-karman-buckling-fields}
\end{figure}

All six loads decrease upon mesh refinement, with shrinking successive absolute changes that suggest convergence; see \cref{fig:von-karman-spectrum}.
As reported in \cref{tab:von-karman-results}, on the seventh mesh, with $\vkPotentialElements$ triangles and $\vkPotentialActiveDofs$ active deflection DOFs, we obtain $\lambda_{\mathrm{cr},h}=\vkPotentialCriticalLoad$.
The last refinement changes the value of $\lambda_{\mathrm{cr},h}$ by about $\vkPotentialLastChangePercent\%$, compared with $\vkPotentialPreviousChangePercent\%$ for the preceding refinement;
the last two values of $\lambda_{\mathrm{cr},h}$ share $\vkPotentialStableDigits$ rounded significant digits.
The displayed residuals are small relative to the mesh-to-mesh changes, and the independent algebraic checks support attributing those changes primarily to discretization.
Finally, the asymmetry indicator $\alpha_h$ decreases from $\vkCoarsestGeometricAsymmetry$ to $\vkFinestGeometricAsymmetry$.

For spatial interpretation, we use the global reconstructed modes $p_h^{k+2}\underline z_{i,h}$ and define the unit-load membrane stress by
\[
   N_h^0|_T
  \coloneqq
  \operatorname{tr}( H_T^{k}\underline\phi_{1,T}^0) I_2
   - H_T^{k}\underline\phi_{1,T}^0.
\]
For display only, each mode in \cref{fig:von-karman-buckling-fields} is normalized by dividing it by its sampled maximum absolute value.
The first two modes illustrate distinct transverse patterns:
the first is even with respect to the vertical centerline, whereas the second changes sign across it.

The stress panel shows the concentration of membrane stress near the hole. Since the reconstructed stress along the discrete flat branch is $\lambda N_h^0$, the displayed unit-load stress magnitudes quantify the local amplification relative to the nominal compression.


\subsubsection{Small forced nonlinear response near the first critical load}
\label{sec:von-karman-nonlinear-response}

\begin{figure}[!htbp]
  \centering
  \includegraphics[width=\textwidth]{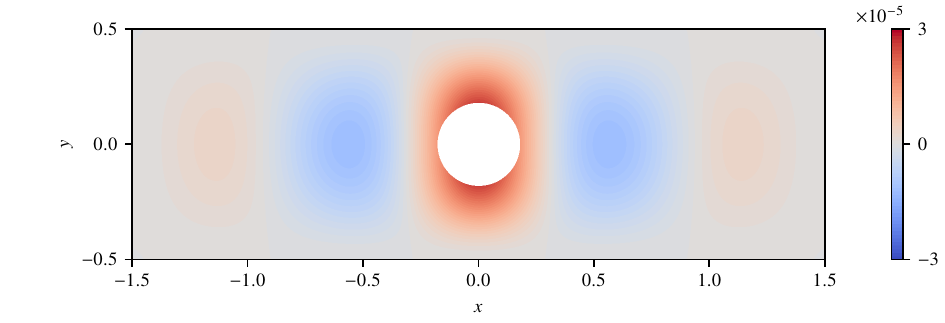}
  \caption{Reconstructed nonlinear deflection on mesh 6 with $k=1$, compression $\lambda=61$, and transverse forcing \eqref{eq:von-karman-small-forcing} with $\varepsilon=10^{-3}$.
    Unlike the eigenmodes in \cref{fig:von-karman-buckling-fields}, the deflection is not normalized; the colorbar shows its actual amplitude.}
  \label{fig:von-karman-nonlinear}
\end{figure}

To compare the nonlinear plate response with the first buckling mode, we solve the fully coupled nonlinear problem \eqref{eq:discrete-von-karman} on the finest mesh. We impose the same boundary conditions as in the spectral study, with values of $\lambda$ close to the first critical value, and apply a small transverse forcing to perturb the flat configuration:
\begin{equation}
  q_\varepsilon(x,y)=\varepsilon\sin\bigl(\pi(x+\tfrac32)\bigr)
  \sin\bigl(\pi(y+\tfrac12)\bigr),
  \qquad \varepsilon=10^{-3}.
  \label{eq:von-karman-small-forcing}
\end{equation}
\Cref{fig:von-karman-nonlinear} shows the reconstructed deflection without rescaling its amplitude.
Its sampled maximum absolute value is approximately $\vkPotentialMaxDeflection$.
The concentration near the hole and the two oppositely signed lateral lobes resemble the first eigenmode in \cref{fig:von-karman-mode-one}.
This is consistent with the proximity of $\lambda=61$ to $\lambda_{1,h}\approx\vkPotentialCriticalLoad$; moreover, the forcing is even in $x$, like the first mode, whereas the second mode is odd in $x$.

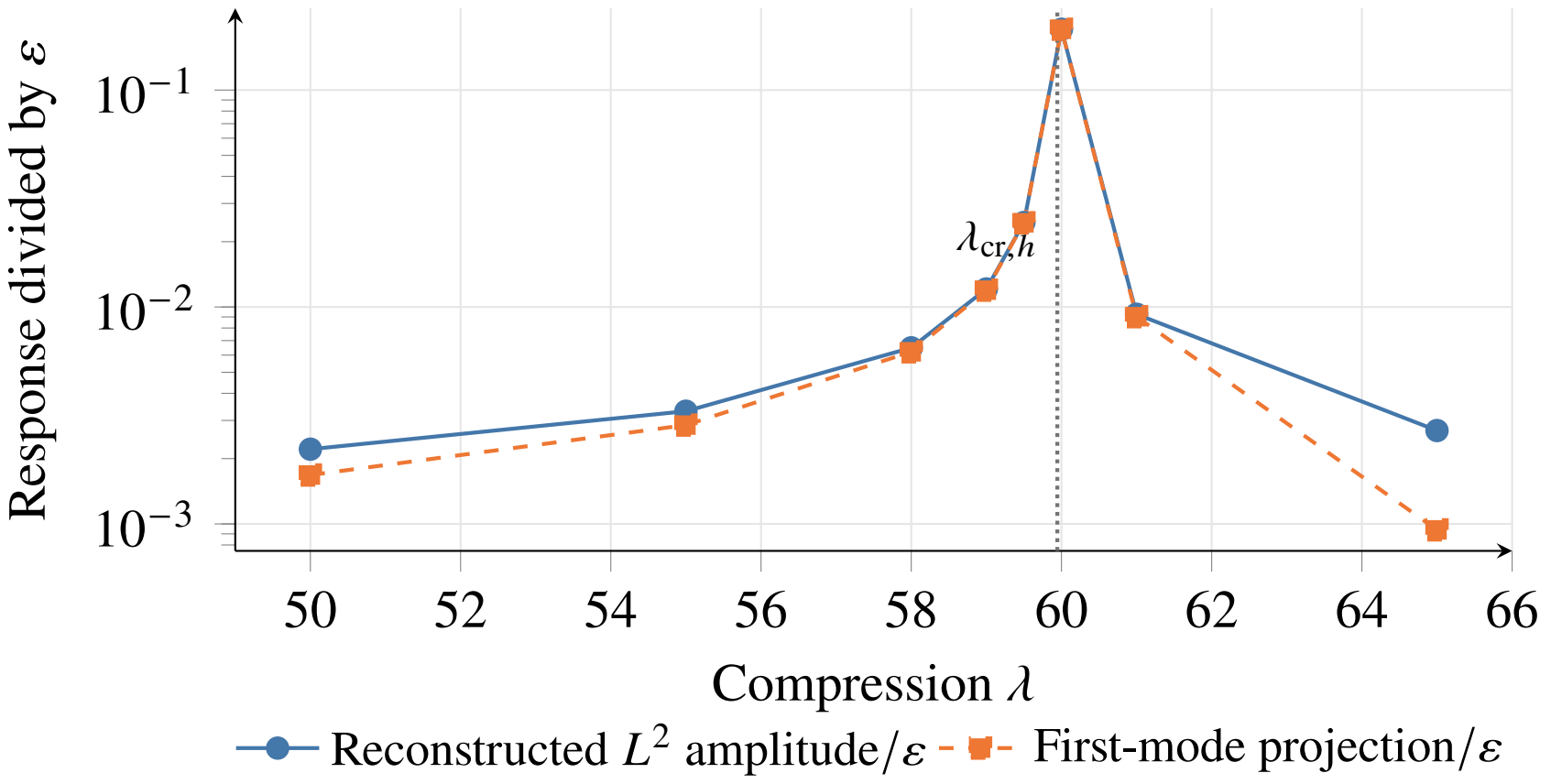
\begin{figure}[!htbp]
  \centering
  \begin{tikzpicture}
    \definecolor{vkblue}{HTML}{4477AA}
    \definecolor{vkorange}{HTML}{EE7733}
    \begin{axis}[
      width=0.72\textwidth,height=5.5cm,
      axis lines=left,tick align=outside,
      xlabel={Compression $\lambda$},
      ylabel={Response divided by $\varepsilon$},
      ymode=log,ymin=\vkPotentialResponsePlotMin,ymax=\vkPotentialResponsePlotMax,
      xmin=49,xmax=66,
      grid=major,grid style={black!10},
      legend style={font=\small,draw=none,fill=none,at={(0.5,-0.28)},anchor=north,legend columns=2},
      every axis plot/.append style={line width=0.8pt,mark size=2pt}
    ]
      \addplot[color=vkblue,mark=*]
        table[x=load,y=L2_over_epsilon]
        {data/von_karman_nonlinear_check/load_response.dat};
      \addlegendentry{Reconstructed $L^2$ amplitude$/\varepsilon$}
      \addplot[color=vkorange,mark=square*,dashed]
        table[x=load,y=mode1_over_epsilon]
        {data/von_karman_nonlinear_check/load_response.dat};
      \addlegendentry{First-mode projection$/\varepsilon$}
      \addplot[black!55,densely dotted,forget plot]
        coordinates {(\vkPotentialCriticalLoad,\vkPotentialResponsePlotMin) (\vkPotentialCriticalLoad,\vkPotentialResponsePlotMax)};
      \node[font=\footnotesize,anchor=south east]
        at (axis cs:\vkPotentialCriticalLoad,\vkPotentialResponseLabelY) {$\lambda_{\mathrm{cr},h}$};
    \end{axis}
  \end{tikzpicture}
  \caption{%
    Amplification of the small forced response around the first discrete critical load.
  }
  \label{fig:von-karman-load-response}
\end{figure}

To quantify this comparison, define
\begin{equation}
  \rho_1(\underline w_h)
  \coloneqq
  \frac{\left|\left(p_h^{k+2}\underline w_h,
          p_h^{k+2}\underline z_{1,h}\right)_{L^2(\Omega)}\right|}
       {\norm{L^2(\Omega)}{p_h^{k+2}\underline w_h}
        \norm{L^2(\Omega)}{p_h^{k+2}\underline z_{1,h}}},
  \label{eq:von-karman-modal-fraction}
\end{equation}
where $\underline{w}_h$ is the solution to the fully nonlinear model, while $\underline{z}_{1,h}$ is the first eigenmode.
We then repeat the nonlinear solve with $\varepsilon=10^{-3}$ for loads on both sides of $\lambda_{1,h}$.
The results in \cref{tab:von-karman-load-response,fig:von-karman-load-response} show a sharply amplified, almost purely first-mode response near the spectral threshold.
At $\lambda=\vkPotentialNearLoad$, only about $\vkPotentialNearDistance$ from the computed first load, one has
$\norm{L^2(\Omega)}{p_h^{k+2}\underline w_h}/\varepsilon\approx\vkPotentialNearLtwoAmplitude$ and
$\rho_1\approx\vkPotentialNearModalFraction$.
Moving away from the threshold reduces both the amplification and, eventually, the first-mode fraction.
All load-sweep computations converged within $\vkPotentialNewtonMax$ Newton corrections.

\begin{table}[!htbp]
  \centering
  \small
  \begin{tabular}{@{}rrrrr@{}}
\toprule
$\lambda$ & $L^2$ amplitude$/\varepsilon$ & First-mode projection$/\varepsilon$ & $\rho_1$ & Newton steps \\
\midrule
50.00 & 0.002211 & 0.001681 & 0.7602 & 2 \\
55.00 & 0.003308 & 0.002855 & 0.8631 & 2 \\
58.00 & 0.006485 & 0.006213 & 0.9582 & 2 \\
59.00 & 0.012110 & 0.011954 & 0.9872 & 2 \\
59.50 & 0.024521 & 0.024440 & 0.9968 & 2 \\
60.00 & 0.190716 & 0.190691 & 0.9999 & 3 \\
61.00 & 0.009260 & 0.009020 & 0.9742 & 2 \\
65.00 & 0.002694 & 0.000939 & 0.3486 & 2 \\
\bottomrule
\end{tabular}

  \caption{Small-forcing load sweep on mesh 6 with $k=1$ and $\varepsilon=10^{-3}$.
    The modal fraction $\rho_1$ is defined by \eqref{eq:von-karman-modal-fraction}.}
  \label{tab:von-karman-load-response}
\end{table}

\section*{Acknowledgements}

Funded by the European Union (ERC Synergy, NEMESIS, project number 101115663).
Views and opinions expressed are however those of the authors only and do not necessarily reflect those of the European Union or the European Research Council Executive Agency. Neither the European Union nor the granting authority can be held responsible for them.


\printbibliography

\end{document}